\documentclass[11pt]{amsart}
\usepackage[T1]{fontenc}
\usepackage{lmodern,amsmath,amssymb,amsthm,booktabs}
\usepackage[a4paper,margin=2.7cm]{geometry}
\usepackage{microtype,mathrsfs}
\usepackage[hypertexnames=false,hidelinks]{hyperref}
\hypersetup{pdftitle={An ambient metric framework for sharp compactness of the constant Q2N-curvature problem}}

\usepackage{pgfplots,float}
\pgfplotsset{compat=1.17}
\definecolor{sectorT}{HTML}{245C89}
\definecolor{sectorV}{HTML}{47745A}
\definecolor{sectorS}{HTML}{BC6628}

\numberwithin{equation}{section}
\newtheorem{theorem}{Theorem}[section]
\newtheorem{proposition}[theorem]{Proposition}
\newtheorem{lemma}[theorem]{Lemma}
\newtheorem{corollary}[theorem]{Corollary}
\theoremstyle{definition}
\newtheorem{definition}[theorem]{Definition}
\theoremstyle{remark}
\newtheorem{remark}[theorem]{Remark}

\DeclareMathOperator{\adj}{adj}
\DeclareMathOperator{\diag}{diag}
\DeclareMathOperator{\divg}{div}
\DeclareMathOperator{\Ric}{Ric}
\DeclareMathOperator{\Span}{span}
\DeclareMathOperator{\supp}{supp}
\DeclareMathOperator{\tr}{tr}
\DeclareMathOperator{\Vol}{Vol}

\newcommand{\ga}{\gamma}
\newcommand{\ep}{\epsilon}
\newcommand{\vep}{\varepsilon}
\newcommand{\vph}{\varphi}
\newcommand{\vrh}{\varrho}
\newcommand{\vth}{\vartheta}
\newcommand{\pp}{\partial}

\newcommand{\Q}{\mathbb{Q}}
\newcommand{\R}{\mathbb{R}}
\newcommand{\Ss}{\mathbb{S}}

\newcommand{\mca}{\mathcal{A}}
\newcommand{\mcb}{\mathcal{B}}
\newcommand{\mcc}{\mathcal{C}}
\newcommand{\mcd}{\mathcal{D}}
\newcommand{\mce}{\mathcal{E}}
\newcommand{\mcf}{\mathcal{F}}
\newcommand{\mch}{\mathcal{H}}
\newcommand{\mcj}{\mathcal{J}}
\newcommand{\mck}{\mathcal{K}}
\newcommand{\mcl}{\mathcal{L}}
\newcommand{\mcm}{\mathcal{M}}
\newcommand{\mcn}{\mathcal{N}}
\newcommand{\mco}{\mathcal{O}}
\newcommand{\mcp}{\mathcal{P}}
\newcommand{\mcq}{\mathcal{Q}}
\newcommand{\mcr}{\mathcal{R}}
\newcommand{\mcs}{\mathcal{S}}
\newcommand{\mct}{\mathcal{T}}
\newcommand{\mcu}{\mathcal{U}}
\newcommand{\mcv}{\mathcal{V}}
\newcommand{\mcy}{\mathcal{Y}}
\newcommand{\mcz}{\mathcal{Z}}

\newcommand{\mfB}{\mathfrak{B}}
\newcommand{\mfF}{\mathfrak{F}}
\newcommand{\mfJ}{\mathfrak{J}}
\newcommand{\mfd}{\mathfrak{d}}
\newcommand{\mfm}{\mathfrak{m}}

\newcommand{\mscD}{\mathscr{D}}
\newcommand{\mscF}{\mathscr{F}}
\newcommand{\mscH}{\mathscr{H}}
\newcommand{\mscP}{\mathscr{P}}

\newcommand{\msM}{\mathsf{M}}
\newcommand{\msP}{\mathsf{P}}
\newcommand{\msS}{\mathsf{S}}
\newcommand{\msT}{\mathsf{T}}
\newcommand{\msV}{\mathsf{V}}
\newcommand{\msZ}{\mathsf{Z}}
\newcommand{\msq}{\mathsf{q}}
\newcommand{\msv}{\mathsf{v}}

\newcommand{\bg}{\bar{g}}
\newcommand{\bph}{\bar{\phi}}
\newcommand{\hc}{\hat{c}}
\newcommand{\hg}{\hat{g}}
\newcommand{\whH}{\widehat{H}}
\newcommand{\wtbeta}{\widetilde{\beta}}
\newcommand{\wtv}{\widetilde{v}}

\newcommand{\mADM}{\mathrm{ADM}}
\newcommand{\mhom}{\mathrm{hom}}
\newcommand{\mlin}{\mathrm{lin}}
\newcommand{\TT}{\mathrm{TT}}
\newcommand{\mfull}{\mathrm{full}}
\newcommand{\msol}{\mathrm{sol}}
\newcommand{\msrc}{\mathrm{src}}

\allowdisplaybreaks

\makeatletter
\renewcommand{\l@subsection}{\@tocline{2}{0pt}{3pc}{3pc}{}}
\makeatother

\title[Compactness for the constant $Q_{2N}$-curvature problem]{An ambient metric framework \\ for sharp compactness of the constant \\ $Q_{2N}$-curvature problem}

\author{Liuwei Gong}
\address[Liuwei Gong]{Department of Mathematics, Chinese University of Hong Kong, Shatin, NT, Hong Kong}
\email{lwgong@math.cuhk.edu.hk}

\author{Seunghyeok Kim}
\address[Seunghyeok Kim]{Department of Mathematics and Research Institute for Natural Sciences, College of Natural Sciences, Hanyang University, 222 Wangsimni-ro Seongdong-gu, Seoul 04763, Republic of Korea}
\email{shkim0401@hanyang.ac.kr}
\email{shkim0401@gmail.com}

\author{Juncheng Wei}
\address[Juncheng Wei]{Department of Mathematics, Chinese University of Hong Kong, Shatin, NT, Hong Kong}
\email{wei@math.cuhk.edu.hk}
\date{November 4, 2026}

\begin{document}
\begin{abstract}
We develop a geometric method for studying $C^{2N}$-compactness of positive solutions to the constant $Q_{2N}$-curvature equation on closed Riemannian $n$-manifolds, for every integer $N\ge1$ and every dimension $n>2N$.

Assuming the GJMS operator is invertible, its Green function is positive, and the positive mass condition holds, we establish sharp dimension ranges for the compactness of the set of all positive solutions with prescribed positive constant $Q_{2N}$-curvature.
The largest dimension for compactness is $24$ for $N=1$ (\cite{KMS}), $2N+20$ for $N=2,3$ (\cite{GKW}), $2N+19$ for $N=4,5$, $2N+18$ for $N=6,7,8$, $2N+17$ for $9\le N\le17$, and $2N+16$ for all $N\ge18$. 
The constructions in \cite{Brendle,BrendleMarques,WeiZhao,GKW} for $N=1,2,3$ and in the companion paper \cite{GKW26} for $N\ge4$ give noncompact examples in every dimension exceeding the corresponding compactness bound.
We show that these $L^{\infty}$-unbounded examples can be chosen to satisfy the same operator and mass hypotheses, proving sharpness within this class.

Our method relates the Pohozaev quadratic forms to the Hessian $D^2\mscF(g_c)$, where $\mscF$ is the volume-normalized total $Q_{2N}$-curvature functional and $g_c$ is the standard round metric on $\Ss^n$.
The ambient metric argument of \cite{GKW26} supplies an explicit formula for $D^2\mscF(g_c)$ on transverse-traceless tensors.
The York decomposition projects metric perturbations onto this space.
The scalar--vector--tensor decomposition represents the Pohozaev quadratic forms by explicit finite-dimensional matrices.
We prove that these matrices are positive definite throughout the stated dimension range, a crucial step in establishing compactness.
\end{abstract}

\maketitle
\tableofcontents

\section{Introduction}\label{sec:introduction}
We study compactness for constant $Q$-curvature equations of arbitrary even order in a fixed conformal class on a closed manifold.
Our main contribution is a geometric method for computing variations of total curvature through ambient metric theory and applying the resulting formulas to conformal variational problems.

Let $(M^n,g)$ be a closed smooth Riemannian manifold and let $N\ge1$ be an integer with $n>2N$.
The order-$2N$ GJMS operator $P_{2N,g}$ on $(M,g)$ is a conformally covariant differential operator introduced by Graham--Jenne--Mason--Sparling \cite{GJMS}, whose principal part is $(-\Delta_g)^N$, where $\Delta_g:=\divg_g\nabla$.
We normalize the associated curvature by $Q_{2N,g}:=2P_{2N,g}1/(n-2N)$.
Conformal covariance implies that a positive function $u$ solves
\begin{equation}\label{eq:constant-q-equation}
P_{2N,g}u=\kappa u^{\frac{n+2N}{n-2N}} \quad \text{on } M
\end{equation}
if and only if $u^{4/(n-2N)}g$ has constant $Q_{2N}$-curvature $2\kappa/(n-2N)$.
We seek $C^{2N}(M)$-compactness of the full set of positive solutions for a fixed $\kappa>0$.

\medskip
The GJMS operators extend the second-order conformal Laplacian and the fourth-order Paneitz operator to higher even orders by means of the ambient metric \cite{GJMS}.
Juhl \cite{Juhl} expressed these operators and their curvatures through compositions of natural second-order operators and volume coefficients, and Fefferman--Graham \cite{FGJuhl} derived the formulas directly from the ambient construction.
These formulas become increasingly complicated as the order grows.
For compactness, the main difficulty is to identify the curvature terms that survive near a concentrating solution and determine their contribution to the Pohozaev identity.

\medskip
Refined bubble analysis in the Yamabe problem provides a model for evaluating these contributions.
Brendle's corrected test functions in \cite{Brendle07} enter the compactness argument of Khuri--Marques--Schoen \cite{KMS}, in which positivity of a Pohozaev quadratic form in local metric perturbations obstructs concentration.
Related reduced energies underlie the noncompactness constructions of Brendle and Brendle--Marques \cite{Brendle,BrendleMarques}.
Related methods have also been used for the fourth- and sixth-order constant $Q$-curvature problems \cite{GKW}, convergence of $Q$-curvature flow \cite{GLW}, isolated singularities of the Yamabe equation \cite{HJWXZI}, and the boundary Yamabe problem \cite{CM,GKMWBoundaryII}.

\medskip
\noindent\textbf{Main results.}
Our first main result gives the compactness bounds under the following operator and mass hypotheses.
We first define the locus on which positive mass is required.
For $n\ge2N+4$, let
\begin{equation}\label{eq:zero-weyl-locus}
\mcz:=\left\{\xi\in M:
\nabla^jW_g(\xi)=0\text{ for }0\le j\le\left\lfloor\frac{n-2N-4}{2}\right\rfloor\right\},
\end{equation}
where $W_g$ is the Weyl curvature tensor; for $2N+1\le n\le2N+3$, set $\mcz:=M$.
Fix a smooth family of conformal-normal representatives $g_\xi:=\Lambda_\xi^{4/(n-2N)}g$, with $\Lambda_\xi(\xi)=1$ and $d\Lambda_\xi(\xi)=0$, as specified in Section \ref{sec:main-mass}.
When $P_{2N,g}$ is invertible, we write $G_g$ for its Green function and $m_g(\xi)$ for its \emph{Green function mass} at $\xi\in\mcz$ in the chosen conformal-normal family.
This is the constant remaining after subtraction of the leading singularity and the nonradial terms through homogeneity zero; see Definition \ref{def:green-mass} and \eqref{eq:main-mass-expansion}.
The hypothesis is
\begin{equation}\label{eq:mass-hypothesis}
m_g(\xi)>0\quad\text{for every }\xi\in\mcz.
\tag{M}
\end{equation}

The dimension range of our compactness theorem is
\begin{equation}\label{eq:main-range}
2N<n\le
\begin{cases}
24 &\text{for }N=1,\\
2N+20 &\text{for }N=2,3,\\
2N+19 &\text{for }N=4,5,\\
2N+18 &\text{for }N=6,7,8,\\
2N+17 &\text{for }9\le N\le17,\\
2N+16 &\text{for }N\ge18.
\end{cases}
\end{equation}
\begin{theorem}
\label{thm:ac-main}
Let $(M^n,g)$ be a closed connected smooth manifold, and let $N\ge1$ and $n$ satisfy \eqref{eq:main-range}.
Assume
\[
\ker P_{2N,g}=0 \quad\text{and}\quad  G_g(x,y)>0\quad\text{for }x\ne y,
\]
and fix a conformal-normal family satisfying \eqref{eq:mass-hypothesis}.
Then, for each fixed $\kappa>0$, the set of positive solutions of \eqref{eq:constant-q-equation}
is uniformly bounded above and below by positive constants and is compact in $C^m(M)$ for every finite integer $m\ge0$.
\end{theorem}
For $N=1,2,3$, noncompactness in every dimension above \eqref{eq:main-range} follows from the constructions of \cite{Brendle,BrendleMarques,WeiZhao,GKW}.
For $N\ge4$, the companion paper \cite[Theorem 1]{GKW26} proves noncompactness throughout the complementary range using perturbations of the form $f(|x|^2)W_{ikj\ell}x^kx^\ell$, where $W$ is an algebraic Weyl tensor and $f$ is a polynomial of degree three.
These metrics can be chosen with a coercive operator and a strictly positive Green function.
Proposition \ref{prop:noncompactness-with-mass} proves that these examples can also be chosen to satisfy \eqref{eq:mass-hypothesis}.
Consequently, the dimension bounds in \eqref{eq:main-range} are sharp within the class satisfying the same operator and mass hypotheses.

Corollary \ref{cor:ac-jet-vanishing} places every limiting blow-up point in $\mcz$, so the mass hypothesis is needed only on that locus.
Under the operator hypotheses of Theorem \ref{thm:ac-main}, the smooth positive mass theorem, in the formulation of \cite[Theorem 1.1]{KWW}, implies \eqref{eq:mass-hypothesis} for $N=1$ whenever $(M,g)$ is not conformally diffeomorphic to the round sphere.
For $N=2$, the same conclusion holds if the Yamabe constant of $[g]$ is also positive \cite[Theorem A]{ALL}.
Corollary \ref{cor:low-order-compactness} records these consequences.
For a general background metric at higher orders, positivity of the Green function mass remains an assumption.

\medskip
\noindent\textbf{Geometric method.}
The method developed in this paper and its companion \cite{GKW26} uses total $Q_{2N}$-curvature as a geometric functional whose variations can be computed from the ambient Ricci-flatness equations.
At closed positive Einstein metrics, this calculation determines the fixed-volume second variation in transverse-traceless directions for arbitrary $N$ with $n>2N$, as stated in \cite[Theorem 2]{GKW26}.
The order $2N$ enters through the scalar ambient extension and the coefficients to be extracted, so the calculation treats all orders through the same differential equations and recurrence relations.
For the compactness and noncompactness problems, conformal covariance relates the quadratic term in the normalized energy $\mce_N$ in \eqref{eq:normalized-energy-quotient} to the second variation of the normalized total curvature functional $\mscF$ in \eqref{eq:normalized-total-curvature}.
York decomposition identifies the scalar correction term required by the linearized equation \eqref{eq:scalar-correction-equation}.
The scalar--vector--tensor decomposition then reduces the Hessian of $\mscF$ at the round metric to scalar radial pairings, from which we prove the positivity of the Pohozaev quadratic forms needed for compactness.

For a smooth trace-free metric variation $h$ at the round metric $g_c$ on the unit sphere $\Ss^n$, the York decomposition \cite[(3.1)]{York1974} reads $h=\Pi h+DV[h]$, where $D$ and $\Pi$ are defined in \eqref{eq:conformal-killing-operator} and \eqref{eq:york-projection}.
The tensor $\Pi h$ is \emph{transverse-traceless (TT)}: its trace and divergence with respect to $g_c$ both vanish.
In stereographic coordinates, we write
\begin{equation}\label{eq:main-notation1}
\Omega(y):=\frac{2}{1+|y|^2},\quad g_c=\Omega^2\delta,\quad\text{and}\quad U(y):=\Omega(y)^{\frac{n-2N}{2}}\quad\text{for }y\in\R^n,
\end{equation}
where $\delta$ is the Euclidean metric.
For $H:=\Omega^{-2}h$, the scalar correction term in \eqref{eq:scalar-correction-definition} is $\Psi[H]:=(n-2N)U\divg_{g_c}V[h]/(2n)$.
Setting $g_t:=\exp(tH)$, the first variation at $t=0$ of $(U+t\Psi[H])^{4/(n-2N)}g_t$ is
\[
h+\frac{4\Psi[H]}{(n-2N)U}g_c=\Pi h+\mcl_{V[h]}g_c.
\]
Thus the contribution of $\Psi[H]$ cancels the term proportional to $g_c$ in $DV[h]$.
Since $g_c$ is critical for $\mscF$, diffeomorphism invariance implies that $D^2\mscF(g_c)$ vanishes whenever one argument is a Lie derivative of $g_c$.
The conformal covariance identity \eqref{eq:energy-curvature-identity} therefore yields
\[
\left.\frac{d^2}{dt^2}\right|_{t=0}\mce_N(g_t,U+t\Psi[H])
=\frac{n-2N}{2}D^2\mscF(g_c)[\Pi h,\Pi h].
\]
Proposition \ref{thm:schur-york-identity-new} extends this identity by cutoff approximation to symmetric two-tensors $H$ whose components in Euclidean coordinates are polynomials of degree at most $d<(n-2N)/2$.
The vector equation determining $V[h]$ is independent of $N$, so the same projection $\Pi$ determines the scalar correction term at every order.
For the Yamabe flow, Brendle \cite{Brendle07} similarly obtains the scalar correction term from a vector field and relates the resulting energy to the second variation of the Einstein--Hilbert functional at the round metric.
The vector-field interpretation also appears in the fourth- and sixth-order problems \cite[Remark E.3]{GKW}.

For $n>2N$, \cite[Theorem 2]{GKW26} computes the second variation of $\mcq(g):=\int_M Q_{2N,g}\,dv_g$ in TT directions along fixed-volume metric variations at a closed Einstein metric $g_E$ normalized by $\Ric_{g_E}=(n-1)g_E$.
To compute the second variation of $\mcq$, let $g(t)$ be a metric path preserving the volume form, with $g(0)=g_E$, and let $g_\rho(t)$ be the corresponding family of metrics on $M$ in the normal form of the ambient metric, so $g_0(t)=g(t)$.
The proof differentiates the ambient Ricci-flatness equations in $t$ before extracting the Taylor coefficients in $\rho$ at $\rho=0$.
When $\left.\pp_t g(t)\right|_{t=0}$ is a TT eigentensor of the Lichnerowicz Laplacian $\Delta_L$, the first variation $\left.\pp_t g_\rho(t)\right|_{t=0}$ is a scalar function of $\rho$ times that eigentensor, and this function satisfies a scalar ordinary differential equation.
Once the first variation is known, the remaining second-order metric information needed to compute the Hessian of $\mcq$ is the scalar function
\[
\int_M\left.\frac{\pp^2}{\pp t^2}\right|_{t=0}
\log\frac{\det g_\rho(t)}{\det g_\rho(0)}\,dv_{g_E}.
\]
This integral satisfies a second scalar ordinary differential equation in $\rho$, so the full tensor $\left.\pp_t^2g_\rho(t)\right|_{t=0}$ need not be determined.
Substituting the resulting Taylor coefficients into the ambient formula for $\mcq$ represents its fixed-volume Hessian in TT directions by a degree-$N$ polynomial in $\Delta_L+2(n-1)$, with the sign convention of \cite[Theorem 2]{GKW26}.
This avoids deriving the full nonlinear $Q_{2N}$-curvature formula separately at each order.
In the critical dimension $n=2N\ge4$, Matsumoto obtained a related second-variation formula for total $Q_n$-curvature at an Einstein representative of a conformally Einstein manifold \cite[Theorem 0.2]{Matsumoto}.
He also obtained factorizations of tensor GJMS operators on TT tensors at Einstein metrics \cite[Theorem 4.2]{Matsumoto}.

The same ambient determinant equation computes the fixed-volume second variations of total renormalized volume coefficients on closed positive Einstein manifolds \cite[Theorem E.1]{GKW26}.
These second variations are strictly negative on nonzero TT tensors at the round metric $g_c$, including in the critical dimension; for $n>2N$, their relation to the second variation of total $Q_{2N}$-curvature is explicit in \cite[Proposition E.3]{GKW26}.
Since the higher-order Einstein tensor is the variational gradient of total $Q_{2N}$-curvature up to normalization, the second-variation formula also determines its linearization in TT directions \cite[Remark 1.2]{GKW26}.

\medskip
We apply the formula \eqref{eq:round-tt-hessian} for $D^2\mscF(g_c)$ to all homogeneous terms relevant to compactness in the conformal-normal metric expansion \eqref{eq:conformal-normal-metric-expansion}.
Proposition \ref{thm:finite-SVT-exhaustion-new} writes the polynomial decomposition of \cite[Corollary 6.10]{GKW} in scalar-, vector-, and tensor-derived families indexed by their angular harmonic degree on $\Ss^{n-1}$.
The projection $\Pi$ sends these tensors to TT tensors on $\Ss^n$; the scalar- and vector-derived tensors can also have nonzero TT projections, with radial components determined by the trace and divergence equations.
The projection $\Pi$ commutes with the $O(n)$ transformations fixing the stereographic axis, so it preserves angular type and harmonic degree but can mix radial profiles corresponding to different homogeneous degrees.
The same invariance of $D^2\mscF(g_c)$ eliminates pairings between distinct angular types and harmonic degrees by Schur's lemma.
Different radial powers within a fixed angular type and degree transform under the same representation of $O(n)$, so their mixed pairings need not vanish.

For each fixed angular harmonic of degree $\ell\ge2$, Proposition \ref{prop:weighted-reconstruction} reconstructs the corresponding TT tensors from one scalar radial function.
Under this reconstruction, the tensor operator $\mca=\nabla^*\nabla+2$ in \eqref{eq:round-tensor-operator}, which enters $D^2\mscF(g_c)$ through \eqref{eq:round-tt-hessian}, acts as the scalar differential operator $A_\ell$ in \eqref{eq:radial-operator}.
The same $A_\ell$ applies to all three angular types, while the tensor inner products reduce to scalar pairings involving the different operator polynomials $W_\sigma(A_\ell)$ in \eqref{eq:repair-norm-operators}.
Proposition \ref{alg:prop:onefunctional} expresses the energy matrix entries for every $\ell\ge2$ through $A_2$ by \eqref{alg:eq:universalenergy}.
By Lemma \ref{alg:lem:budget}, at most nine scalar pairings $m_0,\ldots,m_8$ in \eqref{alg:eq:scalar-pairings} determine these entries in the compactness range; their rational formulas in $N$ and $n$ are given by \eqref{alg:eq:partialfractions}.

\medskip
Under the dilation $H(y)\mapsto H(e^\tau y)$, the negative derivative at $\tau=0$ of the quadratic form with matrix $E^{\mhom}$ in \eqref{eq:repair-E} multiplies its entry pairing homogeneous degrees $d_i,d_j$ by $-(d_i+d_j)$
and is represented by the Pohozaev matrix $H^+$ in \eqref{eq:repair-Hplus}, as proved in Proposition \ref{prop:repair-dilation}.
Negativity of $D^2\mscF(g_c)$ does not by itself imply positive definiteness of $H^+$ in \eqref{eq:repair-Hplus}; Remark \ref{rem:symmetry-breaking} explains this distinction through the first moment of the curvature deviation.

We call $(n-2N)/2$ the \emph{critical degree}, at which the homogeneous energy integral has logarithmic growth.
Using the rational formulas, we prove that the Pohozaev matrices $H^+$ are positive definite for all homogeneous degrees below $(n-2N)/2$ throughout the dimension range \eqref{eq:main-range}.

Injectivity of $\Pi$ on the polynomial spaces, proved in Proposition \ref{prop:repair-injectivity}, then yields a positive quadratic bound in the original conformal-normal metric coefficients.
When the critical degree $(n-2N)/2$ is an integer at least two, Proposition \ref{prop:log-energy} and Lemma \ref{lem:ac-endpoint} establish the corresponding positivity for the leading logarithmic term.

\medskip
These positivity results yield our second main theorem, which provides the coercivity of the Pohozaev quadratic forms used in Theorem \ref{thm:ac-main}, including pairings between homogeneous metric terms of different degrees.
For each integer $d\ge2$, let $\mfJ_d$ be the space of symmetric matrix-valued homogeneous polynomials $H$ of degree $d$ on $\R^n$ satisfying $\tr H=0$ and $H(y)y=0$.

\begin{theorem}\label{thm:intro-algebra}
For every integer $N\ge1$ and dimension $n$ in \eqref{eq:main-range}, the Pohozaev quadratic form
\begin{equation}\label{eq:established-A}
H^+:=-(\mathbf D E^{\mhom}+E^{\mhom}\mathbf D)
\quad\text{on }\bigoplus_{2\le d<\frac{n-2N}{2}}\mfJ_d
\tag{A}
\end{equation}
is positive definite whenever this space is nonzero.
Here $E^{\mhom}$ is the energy matrix \eqref{eq:repair-E}, expressed in the basis indexed by the original homogeneous metric tensors, and $\mathbf D$ multiplies each homogeneous degree-$d$ component by $d$, with block matrices $D_\sigma$ in \eqref{eq:repair-S}.
For every fixed $N\in\mathbb N$, if the critical degree $d=(n-2N)/2$ is an integer at least two, the leading logarithmic coefficient in the oriented cutoff Pohozaev form is positive definite on $\mfJ_d$.
Here the cutoff is supported in $B_{2R}\subset\R^n$ and equals one on $B_R$, with $R\to\infty$; the coefficient is that of $\log R$ in Lemma \ref{lem:ac-endpoint}.
\end{theorem}
The matrix $E^{\mhom}$ uses the weighted radial integral in \eqref{eq:repair-E}; its relation to $D^2\mscF(g_c)$, including the positive volume and angular factors, is stated in \eqref{eq:radial-hessian-normalization}.
Below the critical degree, positivity of $H^+$ gives a negative quadratic contribution to the local Pohozaev integral \eqref{eq:ac-display-4-3}.

For $N=1$, Theorem \ref{thm:intro-algebra} corresponds to \cite[Proposition A.4]{KMS}, and Corollary \ref{cor:ac-jet-vanishing} recovers the Weyl-vanishing conclusion and rates of \cite[Theorem 6.1]{KMS} for the concentrating sequences considered here.
For general $N$, vanishing of the conformal-normal metric terms through degree $j+2$ is equivalent to vanishing of $W_g,\ldots,\nabla^jW_g$, by \eqref{eq:metric-weyl-jets}; this gives the range in \eqref{eq:zero-weyl-locus}.
The coefficient estimates also apply to local critical equations and smoothly varying metrics under the bubble bounds specified in Remark \ref{rem:local-scope}.

\medskip
The algebraic verification in Section \ref{sec:all-order} uses the same calculation of scalar pairings and $LDL^{\msT}$ pivot recurrence for every angular block, including $N=1,2,3$.
For each integer offset $n-2N\in\{5,\ldots,22\}$, the scalar pairings $m_0,\ldots,m_8$ in \eqref{alg:eq:scalar-pairings} determine every maximal block as a rational matrix in $N$.
There are $218$ such matrices and $430$ symbolic pivots.
The ancillary Python program verifies that these pivots are positive at all admissible orders in \eqref{alg:eq:order-intervals}.
At the critical homogeneous degree, Section \ref{sec:log-endpoint} proves positivity of the leading logarithmic coefficient in the cutoff Pohozaev form.
The same pairings also give the first large-order crossings in the tensor-, vector-, and scalar-derived sectors, including their first corrections as $N\to\infty$, in Theorem \ref{thm:three-sector-crossing}.
For sufficiently large $N\in\mathbb N$, the tensor crossing precedes the other two, and the first integer dimension with an indefinite block is $n=2N+17$, one above the compactness bound $2N+16$.
The companion paper \cite{GKW26} constructs an $L^\infty$-unbounded sequence of positive solutions of \eqref{eq:constant-q-equation} for one fixed smooth metric on $\Ss^n$, obtained by perturbing the round metric $g_c$.
The construction uses nondegenerate critical points of the reduced energy in the bubble center and scale variables.

\begin{remark}[Total curvature and the Pohozaev form]\label{rem:symmetry-breaking}
The Hessian $D^2\mscF(g_c)$ in \eqref{eq:round-tt-hessian} is $SO(n+1)$-invariant and, by Proposition \ref{prop:round-hessian-coercivity}, remains negative on TT tensors for every $n>2N$.
Choosing the oriented dilation axis breaks this symmetry to the subgroup $SO(n)$ fixing each pole; the resulting Pohozaev form can couple the radial copies of each angular representation.
Proposition \ref{prop:completion-global-TT} states the exact factor relating the global Pohozaev integral in \eqref{eq:completion-global-coefficient} to the Hessian of the \emph{normalized first moment of the curvature deviation}
\begin{equation}\label{eq:height-curvature-functional}
\mcm_{p_+}(g):=\Vol(g)^{-\frac{n-2N}{n}}
\int_{\Ss^n}Y_{p_+}\bigl(Q_{2N,g}-\overline Q_{2N,g}\bigr)\,dv_g.
\end{equation}
Here $g$ is a smooth metric on $\Ss^n$, and the bar denotes the average with respect to $dv_g$.
We fix $Y_{p_+}(x):=\langle x,p_+\rangle$ for $x\in \Ss^n\subset\R^{n+1}$, where $p_+$ is the pole of stereographic projection.
Thus $Y_{p_+}$ is a fixed first spherical harmonic, and the integral is the component in the direction $p_+$ of the first moment of the signed measure $(Q_{2N,g}-\overline Q_{2N,g})\,dv_g$.
The subtraction of the average curvature is necessary because $Y_{p_+}$ need not have zero average with respect to the varying measure $dv_g$.

The Hessian of $\mcm_{p_+}$ has no preferred sign on the full TT space.
Indeed, if a round rotation $\mathsf O$ exchanges the poles, then $Y_{p_+}\circ\mathsf O=-Y_{p_+}$, and, with $Y_{p_+}$ held fixed,
\[
D^2\mcm_{p_+}(g_c)[\mathsf O^*h,\mathsf O^*h]
=-D^2\mcm_{p_+}(g_c)[h,h].
\]
Our coercivity theorem concerns the images under $\Pi$ of a finite space of polynomial metric jets adapted to one pole; this projected space is not preserved by the rotation.
Below the critical degree, positivity of $H^+$ corresponds to negativity of $D^2\mcm_{p_+}(g_c)$ on the restricted space.
\end{remark}

\medskip
\noindent\textbf{Related results.}
When $N=1$, \eqref{eq:constant-q-equation} reduces to the Yamabe equation.
Schoen's approach \cite[Section 1]{Schoen} combines local blow-up estimates with the positive mass theorem.
Khuri--Marques--Schoen \cite[Theorem 1.1]{KMS} proved compactness through dimension $24$ on manifolds of positive Yamabe type outside the round conformal class, using the positive mass results described in their introduction.
Their argument controls all relevant homogeneous metric terms and proves vanishing of the Weyl tensor and its covariant derivatives at a possible blow-up point.
Brendle \cite{Brendle} constructed noncompact examples for $n\ge52$, and Brendle--Marques \cite{BrendleMarques} extended the construction to $25\le n\le51$.

For $N=2$, G. Li \cite[Theorem 1.2]{GangLi} proved compactness for $5\le n\le7$ under nonnegative scalar curvature and nonnegative, nontrivial $Q_4$-curvature, outside the round conformal class.
Li--Xiong \cite[Theorem 1.1]{LiXiong} proved compactness for $5\le n\le9$ under invertibility of the Paneitz operator, positivity of its Green function, and the corresponding positive-mass property.
Under the same operator assumptions, they also treated locally conformally flat manifolds of positive Yamabe type and, for $n\ge8$, manifolds with nowhere-vanishing Weyl tensor, excluding the round conformal class.
Gong--Kim--Wei \cite[Theorem 1.1]{GKW} extended the dimension range to $5\le n\le24$, outside the round conformal class, under those operator and mass hypotheses or under positivity of the Yamabe constant and nonnegative, nontrivial $Q_4$-curvature.
Wei--Zhao \cite{WeiZhao} constructed noncompact examples for every $n\ge25$.

For $N=3$, Gong--Kim--Wei \cite[Theorem 1.3]{GKW} proved compactness for $7\le n\le26$ outside the round conformal class, under invertibility, positivity of the Green function, and the sixth-order mass hypothesis.
They constructed noncompact examples for $n\ge27$ in \cite[Theorem 1.4]{GKW}.
Their fourth- and sixth-order arguments solve the linearized scalar equations through polynomial systems and decompose the Pohozaev form into the relevant tensor components.
The resulting Weyl-vanishing ranges are $0\le j\le\lfloor(n-8)/2\rfloor$ at order four and $0\le j\le\lfloor(n-10)/2\rfloor$ at order six, whenever these ranges are nonempty \cite[Theorem 1.7]{GKW}.
At order six, positivity of the Green function mass on a general background remains an assumption.
Our argument gives a new geometric proof of compactness for $N=1,2,3$ in these dimension ranges, under the operator and mass hypotheses of Theorem \ref{thm:ac-main}.

For arbitrary $N$, Mazumdar--Premoselli \cite[Theorem 1.1]{MP} established compactness in dimensions $2N+1\le n\le2N+5$ under invertibility and positivity of the Green function, assuming positive mass everywhere when $n\le2N+3$ and at points where the Weyl tensor vanishes when $2N+4\le n\le2N+5$.
They also treated locally conformally flat manifolds with positive mass and all dimensions $n\ge2N+4$ when the Weyl tensor is nowhere zero.
Their local estimates and the classification of entire positive flat solutions by Wei--Xu \cite[Theorem 1.3]{WX} provide the starting point for the blow-up analysis used here; the precise statements and hypotheses are recorded in Section \ref{sec:ac-tools}.

\medskip
\noindent\textbf{Positive mass at second and higher orders.}
For $N=1$, the Green function mass at the relevant poles is a positive multiple of the Arnowitt--Deser--Misner (ADM) mass of the scalar-flat metric obtained from the Green function; see Proposition \ref{prop:low-order-mass}.
Following the classical proofs of Schoen--Yau \cite{SY} and Witten \cite{Witten}, Bi--Hao--He--Shi--Zhu \cite{BHHZ} proved the smooth theorem through dimension $19$, and Brendle--Wang \cite[Theorem 1.5]{BW} established the all-dimensional result.
Khuri--Wang--Wang extended the theorem to singular metrics; we use the smooth positivity-and-rigidity statement of \cite[Theorem 1.1]{KWW} in Corollary \ref{cor:low-order-compactness}.

For $N=2$, earlier results for the Paneitz mass include \cite{HumbertRaulot,GurskyMalchiodi,HangYang}.
Under positive Yamabe constant and our operator hypotheses, Proposition \ref{prop:low-order-mass} identifies the mass at the required poles with a positive multiple of the fourth-order energy, using \cite[Proposition 8.3]{GKW}; positivity and rigidity follow from \cite[Theorem A]{ALL}.

For $N\ge3$, positivity of mass on a general background remains an assumption; conformal covariance and results under additional geometric hypotheses are treated in \cite{Michel,Ludewig}.
Proposition \ref{prop:positive-mass-ray} proves positivity near $g_c$ for every integer $N\ge1$: if $n>2N+3$ and $DW_{g_c}[H]\not\equiv0$, the mass along $g_t:=g_c\exp(tg_c^{-1}H)$ is bounded below by $ct^2$ for some $c>0$, uniformly on $\mcz(g_t)$ for small $t\ne0$.
The second-variation identity \eqref{eq:mass-ray-second-variation} proves this bound.

\medskip
\noindent\textbf{From coercivity to compactness.}
To apply the positive definiteness in Theorem \ref{thm:intro-algebra}, we compare the Pohozaev quadratic forms with the localized Pohozaev integral of a concentrating solution.
The coarse estimate \eqref{eq:ac-display-4-5}, together with the blow-up analysis in Section \ref{sec:ac-separation}, yields isolated simple blow-up points and quantitative vanishing of the metric coefficients, but its error is not yet small enough to determine the mass term.
Matching the bubble and its scalar correction term with the Green function expansion \eqref{eq:repair-Green-expansion} improves the annular error, as shown in Proposition \ref{prop:ac-matched-cap}.
Lemma \ref{lem:repair-angular-mass} shows that the nonradial terms do not contribute to the limiting Pohozaev integral.
The refined estimate \eqref{eq:ac-display-8-2} then contradicts the positive limit \eqref{eq:ac-display-8-5} implied by \eqref{eq:mass-hypothesis}.
In both limits, the solution index $i$ tends to infinity at a fixed coordinate radius $\rho>0$ before $\rho$ tends to zero.
Proposition \ref{prop:low-order-mass} identifies this Green function coefficient with the corresponding mass at orders two and four.

\medskip
\noindent\textbf{Background metrics.}
Theorem \ref{thm:ac-main} applies to general closed background metrics satisfying its operator and mass hypotheses.
The round sphere enters as the conformal compactification of the local bubble, so its use imposes no Einstein assumption on the original manifold.
The explicit harmonic decomposition and rational sign formulas use the round metric $g_c$, while the ambient second-variation formulas apply to closed positive Einstein metrics under their stated hypotheses.

\medskip
\noindent\textbf{Organization of the paper.}
Section \ref{sec:preliminaries} fixes the bubble and mass normalizations and recalls the ambient results used below.
Section \ref{sec:ambient-york} relates the second variation of $\mce_N$ along the scalar correction term to $D^2\mscF(g_c)$ and decomposes the homogeneous tensors $H\in\mfJ_d$.
Section \ref{sec:all-order} proves the exact matrix signs needed for compactness and, separately, the large-order crossings in Theorem \ref{thm:three-sector-crossing}.
Section \ref{sec:weighted-local-model} constructs the local corrections, treats the logarithmic endpoint, and proves the Green function expansion and mass identities.
Section \ref{sec:actual-comparison} proves that blow-up points are isolated simple and establishes quantitative vanishing of the metric coefficients.
Section \ref{sec:green-matching} matches the local profiles to the Green functions and completes the proof of Theorem \ref{thm:ac-main}.
Appendices \ref{app:weighted} and \ref{app:verification} contain the tensor reconstruction, rational scalar pairings, and exact-arithmetic verification.
Appendix \ref{app:noncompactness-mass} proves the mass assertion for the noncompactness examples.

\medskip
\noindent\textbf{Acknowledgement.}
L. Gong thanks Mingxiang Li, Yanyan Li and Wei Yuan for helpful discussions.
S. Kim thanks Haixia Chen for stimulating discussion.
L. Gong and J. Wei gratefully acknowledge the support of the Research Grants Council of Hong Kong (RGC) for the project titled \emph{On Fujita equation in the critical or supercritical regime}.
S. Kim acknowledges support from the Basic Science Research Program through the National Research Foundation of Korea (NRF), funded by the Ministry of Science and ICT (RS2025-00558417), as well as from the Samsung Science and Technology Foundation (Project Number SSTF-BA2601-01).

\medskip
\noindent\textbf{AI assistance.}
The authors formulated the method of relating compactness and noncompactness problems to the second variation of total $Q_{2N}$-curvature in TT directions and analyzing it through ambient metric theory independently of AI assistance in September 2025, following completion of \cite{GKW}.
OpenAI's Codex was subsequently used to assist with algebraic calculations, ambient recurrences and polynomial coefficient choices, develop exact-arithmetic verification programs, refine exposition, notation and proof details, and carry out additional proof and consistency checks.
The finite computations used in the paper are specified in Appendix \ref{app:verification}.
The authors are responsible for the final text and mathematical claims.

\section{Preliminaries}\label{sec:preliminaries}
We write $\operatorname{Scal}_g$ and $\operatorname{Rm}_g$ for scalar and Riemann curvature, respectively; plain $R$ denotes a Euclidean cutoff radius.

A polynomial symmetric two-tensor $H=(H_{ij})$ on $\R^n$ means a covariant two-tensor with $H_{ij}=H_{ji}$ whose components are real polynomials in the standard Euclidean coordinates.
For an integer $d\ge0$, the condition $\deg H\le d$ means that every component has total degree at most $d$; homogeneity of degree $d$ also refers to the components.
We use such tensors as directions of metric variation.
For such an $H$, the metric path $g_t:=\exp(tH)$ means $\delta\exp(t\delta^{-1}H)$, where the exponential on the right is the matrix exponential.

\subsection{Bubbles and Green function mass}\label{sec:main-mass}
For the standard bubble $U$ in \eqref{eq:main-notation1}, we define the dilation operator
\begin{equation}\label{eq:main-notation2}
\msZ := y\cdot\nabla+\frac{n-2N}{2}.
\end{equation}
Thus $\msZ U$ is the derivative at $\lambda=1$ of $\lambda^{(n-2N)/2}U(\lambda y)$.
We define
\begin{equation}\label{eq:bubble-linearization}
P_0:=(-\Delta)^N ,\quad 
\kappa_0:=P_{2N,g_c}1
=\frac{\Gamma\left(\frac{n}{2}+N\right)}{\Gamma\left(\frac{n}{2}-N\right)} ,\quad\text{and}\quad 
L_U:=(-\Delta)^N-\frac{n+2N}{n-2N}\kappa_0U^{\frac{4N}{n-2N}}.
\end{equation}
Here $L_U$ is the linearization at $v=U$ of the Euclidean equation $(-\Delta)^Nv=\kappa_0v^{(n+2N)/(n-2N)}$.
In the conformal energy space $\mch_N$ of functions on $\R^n$ defined in \eqref{eq:conformal-energy-space}, the kernel of $L_U$ is generated by the derivatives of the bubble under changes of center and scale; these are the \emph{Jacobi fields}.
If $u>0$ solves \eqref{eq:constant-q-equation} with coefficient $\kappa>0$, then
\[
\widetilde u:=\left(\frac{\kappa}{\kappa_0}\right)^{\frac{n-2N}{4N}}u
\quad\text{satisfies}\quad
P_{2N,g}\widetilde u=\kappa_0\widetilde u^{\frac{n+2N}{n-2N}}.
\]
We therefore take $\kappa=\kappa_0$ in the compactness proof.

For a sequence $(u_i)_{i\in\mathbb N}$ of these normalized solutions with $u_i(\xi_i)\to\infty$ at local maximum points $\xi_i$, we define the initial scale for each $i\in\mathbb N$ by
\begin{equation}\label{eq:bubble-height-scale}
\mu_i^{(0)}:=2u_i(\xi_i)^{-\frac{2}{n-2N}},
\quad\text{so that}\quad
(\mu_i^{(0)})^{-\frac{n-2N}{2}}U(0)=u_i(\xi_i).
\end{equation}
Thus, in normal coordinates centered at $\xi_i$, the bubble $(\mu_i^{(0)})^{-(n-2N)/2}U(x/\mu_i^{(0)})$ has the same value at $x=0$ as $u_i$ at $\xi_i$.
Lemma \ref{lem:ac-decomposition} then adjusts the center and scale so that the product $\chi_0(v_i-U-\psi_i)$ is orthogonal to the Jacobi fields in the energy inner product, as specified by \eqref{eq:ac-solution-projection}.
Here $\chi_0$ is the smooth radial cutoff function in that lemma, and $v_i$ and $\psi_i$ are defined in \eqref{eq:ac-local-bubble} and \eqref{eq:ac-display-2-1}.
The adjusted scale $\mu_i$ satisfies
\begin{equation}\label{eq:bubble-adjusted-scale}
\frac{\mu_i}{\mu_i^{(0)}}=1+o(1)\quad\text{as }i\to\infty.
\end{equation}

At each isolated blow-up point, Corollary \ref{cor:ac-jet-vanishing} fixes a coordinate radius $\rho_0>0$ for the orthogonality conditions.
For each $i\in\mathbb N$, we use $\chi_0$ with $R=\rho_0/\mu_i$: in the rescaled coordinates, $\chi_0(y)=1$ for $|y|\le2\rho_0/\mu_i$ and $\chi_0(y)=0$ for $|y|\ge5\rho_0/(2\mu_i)$.
We retain the chosen center and scale when replacing $\rho_0$ by a smaller coordinate radius $0<\rho<\rho_0$; the corresponding rescaled radius is $R_i=\rho/\mu_i$.
In the final mass comparison in Section \ref{sec:fine-completion}, we first let $i\to\infty$ with $\rho$ fixed, so that $R_i\to\infty$, and then let $\rho\downarrow0$.
The local Pohozaev integral \eqref{eq:ac-display-4-3} uses the operator difference $P_0-P_{2N,g}$.

\medskip
The conformal-normal convention specifying the finite part in \eqref{eq:mass-hypothesis} is as follows.
When $P_{2N,g}$ is invertible, $G_g(\xi,\cdot)$ denotes its Green function, normalized by $P_{2N,g}G_g(\xi,\cdot)=\delta_\xi$.
We choose a smooth family of conformal-normal metrics as in \cite[Theorem 5.1]{LP}:
\[
g_\xi:=\Lambda_\xi^{\frac{4}{n-2N}}g ,\quad\text{where } 
\Lambda_\xi>0 ,\quad  \Lambda_\xi(\xi)=1 ,\quad\text{and}\quad  d\Lambda_\xi(\xi)=0.
\]
In geodesic normal coordinates $x$ centered at $\xi$, the determinant is normalized by $\det g_\xi(x)=1+O(|x|^T)$ for a fixed integer $T>n+4$.
We write the metric as the matrix exponential of a symmetric tensor and expand its logarithm:
\begin{equation}\label{eq:conformal-normal-metric-expansion}
g_\xi(x)=\exp H_\xi(x) ,\quad\text{where } 
H_\xi(x)=\sum_{d=2}^{T-1}H_{\xi,d}(x)+O(|x|^T).
\end{equation}
Here $H_{\xi,d}$ is homogeneous of degree $d$, and normal coordinates and determinant normalization imply $H_{\xi,d}(x)x=0$ and $\tr H_{\xi,d}=0$.
We refer to this finite Taylor expansion as the \emph{conformal-normal metric jet}.
The family and the normalization order are fixed throughout the argument.
Equivalently, the set in \eqref{eq:zero-weyl-locus} is
\begin{equation}\label{eq:zero-jet-locus}
\mcz=\{\xi\in M:H_{\xi,d}=0\text{ for }2\le d\le\lfloor\frac{n-2N}{2}\rfloor\}.
\end{equation}
For $n\ge2N+4$, the equivalence follows from the relation between metric Taylor coefficients and covariant derivatives of the Weyl tensor \eqref{eq:metric-weyl-jets}.

At $\xi\in\mcz$, Theorem \ref{thm:repair-Green} proves the expansion
\begin{equation}\label{eq:main-mass-expansion}
\begin{split}
G_{g_\xi}(\xi,x)={}&\ga_{n,N}r^{2N-n}
+\sum_{d=d_0}^{n-2N}r^{d+2N-n}a_{\xi,d}(\theta)
+m_g(\xi)+O_{C^{2N}}(r^\vep),\\
d_0:={}&\max\{2,\lfloor\frac{n-2N}{2}\rfloor+1\} \quad\text{and}\quad 
\ga_{n,N}:=\frac{\Gamma\left(\frac{n-2N}{2}\right)}{4^N\pi^{\frac{n}{2}}\Gamma(N)}.
\end{split}
\end{equation}
Here $r:=|x|$, $\theta:=x/r$, each $a_{\xi,d}$ has zero integral over $\Ss^{n-1}$, and $0<\vep<\min\{1,2d_0-n+2N\}$.
The notation $O_{C^{2N}}(r^\vep)$ means that its derivatives of order $j$ are bounded by $Cr^{\vep-j}$ for $0\le j\le2N$.
\begin{definition}[Green function mass]\label{def:green-mass}
For $\xi\in\mcz$, the \emph{Green function mass} in the fixed conformal-normal family is the constant $m_g(\xi)$ in \eqref{eq:main-mass-expansion}: the finite part after subtraction of the radial singularity and all the displayed nonradial terms, including the bounded term when it is present.
Since the angular terms have zero spherical mean, this constant also satisfies
\begin{equation}\label{eq:mass-spherical-average}
m_g(\xi)
=\lim_{r\downarrow0}\left[
\frac{1}{|\Ss^{n-1}|}\int_{\Ss^{n-1}}
G_{g_\xi}(\xi,r\theta)\,d\theta-\ga_{n,N}r^{2N-n}\right].
\end{equation}
\end{definition}
When $2N+1\le n\le2N+3$, no angular terms remain.
The mass is then the constant term studied by Michel in \cite[Theorem 2.3]{Michel}, with conformal covariance as in \cite[Theorem 3.1]{Michel}.

For $N=3$, the mass normalization in \eqref{eq:main-mass-expansion} agrees with \cite[(10.44)]{GKW}.

\subsection{Ambient metric jets and the GJMS construction} \label{sec:ambient-preliminaries}
The \emph{ambient metric} is a metric in two auxiliary variables whose Ricci tensor vanishes to the order needed to define conformal invariants of $g$.
In the normal form recalled in \cite[Subsection 2.1]{GKW26}, its coefficients are encoded by a family $g_\rho$ with $g_0=g$.
The ambient Ricci equations determine this family, and powers of the ambient Laplacian acting on homogeneous extensions of functions determine the GJMS operators.
For $n>2N$, the Taylor coefficients $\pp_\rho^j g_\rho|_{\rho=0}$ through $j=N$ are uniquely determined before any even-dimensional obstruction occurs \cite[Theorem 2.9]{FGAmbient}.
Only these finite jets enter $P_{2N,g}$ and its metric variations.
They depend smoothly on a smooth variation of $g$, since the finite coefficient equations have nonzero dimensional denominators in this range.
With the normalization $P_{2N,g}1=(n-2N)Q_{2N,g}/2$ and leading term $(-\Delta_g)^N$, conformal covariance reads
\begin{equation}\label{eq:ambient-conformal-covariance}
P_{2N,e^{2\omega}g}u
=e^{-\left(\frac{n}{2}+N\right)\omega}P_{2N,g}\bigl(e^{\left(\frac{n}{2}-N\right)\omega}u\bigr)
\end{equation}
for every smooth real function $\omega$; see \cite[(2.6)]{GKW26}.

At a positive Einstein metric $g_E$ normalized by $\Ric_{g_E}=(n-1)g_E$, the normal-form family is $g_\rho=(1+\rho/2)^2g_E$ \cite[Lemma 2.2]{GKW26}.
For a TT eigentensor of the Lichnerowicz Laplacian $\Delta_L$, the first variation of this metric family is determined by the scalar equation \cite[(3.2)]{GKW26}.
The recurrence \cite[(3.4)]{GKW26} determines its Taylor coefficients in $\rho$; the recurrence \cite[(3.7)]{GKW26} determines those of the integrated second variation of $\log\det g_\rho$.
These coefficients enter the fixed-volume second variation of total $Q_{2N}$-curvature through \cite[(3.12)]{GKW26}.

\section{Second variation of total curvature and projection onto TT tensors}\label{sec:ambient-york}
We relate the second variation of the normalized energy $\mce_N$ in \eqref{eq:normalized-energy-quotient}, after including the scalar correction term $\Psi[H]$, to $D^2\mscF(g_c)$ by \eqref{eq:schur-york-identity-new}.
The cutoff estimates \eqref{eq:strict-york-cutoff-estimates} justify \eqref{eq:schur-york-identity-new} for the possibly singular tensors $h=\Omega^2H$ obtained from polynomial symmetric two-tensors $H$ on $\R^n$.

\subsection{Coercivity of the normalized total-curvature Hessian} \label{sec:round-hessian}
The normalized total curvature functional $\mscF$ is
\begin{equation}\label{eq:normalized-total-curvature}
\mscF(g):=\Vol(M,g)^{-\frac{n-2N}{n}}
\int_M Q_{2N,g}\,dv_g.
\end{equation}
The round metric $g_c$ is critical for $\mscF$: its variational gradient is rotation invariant and hence a multiple of $g_c$, and scale invariance forces this multiple to vanish.
The trace, divergence, and tensor pairings in the formula for $D^2\mscF(g_c)$ are taken with respect to $g_c$.
Let $\nabla^*$ be the formal $L^2(g_c)$ adjoint of the covariant derivative and set
\begin{equation}\label{eq:round-tensor-operator}
\mca:=\nabla^*\nabla+2.
\end{equation}
Thus $\nabla^*\nabla$ is the nonnegative rough Laplacian on tensors.
The eigenvalues of $\mca$ on TT tensors are $k(k+n-1)$, $k=2,3,\ldots$.
In the notation of \cite[Theorem 2]{GKW26}, we define
\begin{equation}\label{eq:ambient-hessian-coefficient}
p_{N,n}(x):=\mscP_N^{(n)}(-x).
\end{equation}
Here $\mscP_N^{(n)}$ is the polynomial representing the fixed-volume second variation of total $Q_{2N}$-curvature on TT tensors in the cited theorem.
The finite coefficient formula \cite[(3.12)]{GKW26} defines the same polynomial for real $n>2N$ and implies $p_{N,n}(0)=0$.
For every pair of smooth TT tensors $h_1,h_2$, \cite[Theorem 2]{GKW26} states
\begin{equation}\label{eq:round-tt-hessian}
D^2\mscF(g_c)[h_1,h_2]
=|\Ss^n|^{-\frac{n-2N}{n}}
\langle h_1,p_{N,n}(\mca)h_2\rangle_{L^2(g_c)}.
\end{equation}
Writing $DP_{2N,g}[h]$ for the first metric derivative in direction $h$, the same theorem also proves
\begin{equation}\label{eq:round-tt-first-variation}
DP_{2N,g_c}[h]1=0 \quad \text{for every smooth TT tensor }h.
\end{equation}
We write $H^N_{\TT}(\Ss^n)$ for the symmetric covariant two-tensors in $H^N(\Ss^n,g_c)$ whose trace and divergence with respect to $g_c$ vanish in the weak sense.

\begin{proposition}
\label{prop:round-hessian-coercivity}
For every integer $N\geq1$ and real $n>2N$, the polynomial $p_{N,n}$ satisfies
\begin{equation}\label{eq:hessian-polynomial-negativity}
p_{N,n}(x)\leq-\frac{((N-1)!)^2}{4(n-1)}
\binom{n-2}{N-1}x<0 \quad\text{for }x>0.
\end{equation}
For integer $n>2N$, there is $c_{N,n}>0$ such that
\begin{equation}\label{eq:round-hessian-coercivity}
-D^2\mscF(g_c)[h,h]\geq c_{N,n}\|h\|_{H^N}^2
\quad\text{for }h\in H^N_{\TT}(\Ss^n).
\end{equation}
\end{proposition}
\begin{proof}
We use the coefficient formula \cite[(3.12)]{GKW26} directly.
Set $a(\rho):=1+\rho/2$ and $\nu:=N-n/2<0$, and let $\psi_{-x}$ and $f_N$ be the polynomials through degree $N-1$ determined by \cite[(3.4)]{GKW26} and \cite[(3.10)]{GKW26}, respectively.
In this proof, primes denote derivatives in $\rho$, and $[\rho^j]$ denotes extraction of the coefficient of $\rho^j$.
For $x>0$, the recurrence for $\psi_{-x}$ has nonnegative multipliers and positive denominators $n-2m-2$ at every required index, so its coefficients are nonnegative.
The coefficients of $f_N$ are positive as well: the signs of $(1-N)_r$ and $(-\rho/2)^r$ cancel for $0\le r<N$.
Here $(\alpha)_r:=\prod_{j=0}^{r-1}(\alpha+j)$, with $(\alpha)_0:=1$.
Define
\[
J_x(\rho):=\frac{x}{2(n-1)}+
\int_0^\rho a(t)^2\bigl(\psi_{-x}'(t)\bigr)^2\,dt,
\qquad
B_N(\rho):=\frac12f_N(\rho)^2-\frac{\rho}{\nu}f_N(\rho)f_N'(\rho).
\]
The determinant equation \cite[(3.6)]{GKW26} and the coefficient formula yield
\[
p_{N,n}(x)=-2^{N-1}((N-1)!)^2
[\rho^{N-1}]\,a^{n-2}J_xB_N.
\]
All factors have nonnegative coefficients through order $N-1$ for real $n>2N$.
Retaining only $J_x(0)=x/[2(n-1)]$, $B_N(0)=1/2$, and $[\rho^{N-1}]a^{n-2}=2^{1-N}\binom{n-2}{N-1}$ proves \eqref{eq:hessian-polynomial-negativity}, including $N=1$.

For $N\ge2$, the leading coefficient in \cite[(3.14)]{GKW26} has sign $(-1)^{N+1}$ as a polynomial in the spectral variable.
After the substitution in \eqref{eq:ambient-hessian-coefficient}, it follows that
\[
-p_{N,n}(x)=\hc_{N,n}x^N+O(x^{N-1}) ,\quad\text{where }  \hc_{N,n}>0.
\]
For $N=1$, the same conclusion follows from $p_{1,n}(x)=-x/(4(n-1))$.
Combining this asymptotic with \eqref{eq:hessian-polynomial-negativity} proves $-p_{N,n}(x)\ge c(1+x)^N$ for $x\ge2(n+1)$.
The spectral decomposition of $\mca$ on TT tensors and elliptic norm equivalence applied to \eqref{eq:round-tt-hessian} prove \eqref{eq:round-hessian-coercivity} on finite spectral sums.
Density extends it to $H^N_{\TT}$.
\end{proof}

\subsection{Projection onto TT tensors and the scalar correction term}\label{sec:scalar-complement-new}
In the stereographic coordinates of \eqref{eq:main-notation1}, write $g_c=\Omega^2\delta$ and let $p_+$ be the projection pole.
For a polynomial symmetric two-tensor $H$ on $\R^n$, the tensor $h:=\Omega^2H$ expressed in this chart is defined on $\Ss^n\setminus\{p_+\}$ and may be singular at $p_+$.
To evaluate $D^2\mscF(g_c)$ on the TT projection $\Pi h$, we first establish continuity of $\Pi$ and convergence under cutoff.
We write $\vrh_+:=d_{g_c}(\cdot,p_+)$ for the round distance to the projection pole.
In this subsection, we use the unweighted $L^2$-based Sobolev spaces of functions and tensor fields on the round sphere $(\Ss^n,g_c)$.
The \emph{conformal Killing operator} removes the trace from the infinitesimal change of metric generated by a vector field:
\begin{equation}\label{eq:conformal-killing-operator}
\mcd_gV:=\mcl_Vg-\frac{2}{n}(\divg_gV)g.
\end{equation}
For the round metric, we write $D:=\mcd_{g_c}$, with formal $L^2(g_c)$ adjoint $D^*$, and set
\begin{equation}\label{eq:york-projection}
V[h]:=(D^*D)^\dagger D^*h \quad\text{and}\quad 
\Pi h:=h-DV[h].
\end{equation}
Here $(D^*D)^\dagger$ is the inverse of the elliptic vector operator $D^*D$ on the $L^2$ orthogonal complement of the conformal Killing fields $\ker D$, extended by zero on that kernel.
For a trace-free symmetric covariant two-tensor $h$, the divergence is the one-form with components $(\divg_{g_c}h)_k=\nabla^i h_{ik}$, where $\nabla$ is the Levi-Civita connection of $g_c$.
The adjoint $D^*h$ is the vector field obtained by multiplying this one-form by $-2$ and contracting it with the inverse metric $g_c^{-1}$:
\[
(D^*h)^j=-2(g_c)^{jk}\nabla^i h_{ik},
\qquad j\in\{1,\ldots,n\}.
\]
Thus $\Pi$ is the $L^2$ orthogonal projection onto TT tensors; \eqref{eq:york-projection} is the $n$-dimensional form of York's decomposition \cite[Section 3]{York1974}.
These definitions also apply to trace-free tensor distributions.
We write $S^2_0T^*\Ss^n$ for the bundle of symmetric covariant two-tensors with zero $g_c$-trace.

\begin{lemma}
\label{lem:sobolev-york-new}
For every real $r$, the maps
\[
V:H^r(S^2_0T^*\Ss^n)\to H^{r+1}(T\Ss^n) \quad\text{and}\quad  \Pi:H^r(S^2_0T^*\Ss^n)\to H^r(S^2_0T^*\Ss^n)
\]
are bounded.
Moreover $D^*\Pi h=0$ on the whole sphere in the weak sense, and $V[h]\perp\ker D$ when $h\in L^2$.
\end{lemma}
\begin{proof}
Elliptic theory on the closed sphere makes $(D^*D)^\dagger$ a pseudodifferential operator of order $-2$.
For every conformal Killing field $Y$, the definition of the formal adjoint gives $\langle D^*h,Y\rangle=\langle h,DY\rangle=0$.
Consequently $V$ has order $-1$, $\Pi$ has order zero, and $D^*DV[h]=D^*h$.
Their orders imply the asserted Sobolev bounds, and the definition of the generalized inverse ensures $V[h]\perp\ker D$.
\end{proof}

Let $d\ge0$ be an integer and let $H$ be a symmetric Euclidean trace-free polynomial tensor of degree at most $d$.
Writing $B_R:=\{x\in\R^n:|x|<R\}$, choose $\chi\in C_c^\infty(\R^n)$ equal to one on $B_1$ and zero outside $B_2$, and, for $R\ge1$, set
\[
H_R(x):=\chi\left(\frac{x}{R}\right)H(x) ,\quad 
h_R:=\Omega^2H_R ,\quad\text{and}\quad  h:=\Omega^2H.
\]
We regard $h_R$ as a smooth tensor on the whole sphere by extending it by zero near $p_+$, and $h$ initially as a tensor on $\Ss^n\setminus\{p_+\}$.
The estimates \eqref{eq:strict-york-cutoff-estimates} also hold for any cutoff family equal to one on $B_R$, supported in $B_{2R}$, and with derivatives of order $k$ bounded by $C_kR^{-k}$.

\begin{proposition}
\label{thm:strict-york-transmission-new}
Fix an integer $N\ge1$ and assume
\begin{equation}\label{eq:strict-degree-domain-new}
n>2d+2N.
\end{equation}
Then $h\in H^N$ and, for all $R\ge1$,
\begin{equation}\label{eq:strict-york-cutoff-estimates}
\|h_R-h\|_{H^N}\le C R^{d+N-\frac{n}{2}} \quad\text{and}\quad 
\|V[h_R]-V[h]\|_{H^{N+1}}
+\|\Pi h_R-\Pi h\|_{H^N}\le C R^{d+N-\frac{n}{2}}.
\end{equation}
The constants are independent of $R$.
In particular, $\Pi h$ is TT in the weak sense on the whole sphere.
\end{proposition}
\begin{proof}
The estimate for $h_R-h$ is \cite[(4.12)]{GKW26}, whose integrability hypothesis is exactly \eqref{eq:strict-degree-domain-new}.
Its proof uses the bounds $|\nabla^k h|\le C_k\vrh_+^{-d-k}$ near $p_+$ and the stated derivative bounds on the cutoff, so it applies to every cutoff family specified above.
Lemma \ref{lem:sobolev-york-new} then proves the estimates for $V$ and $\Pi$, including the weak TT condition on the whole sphere.
\end{proof}

\begin{corollary}
\label{cor:strict-hessian-transmission-new}
Let $\mcb_N:=D^2\mscF(g_c)$.
It extends to a continuous bilinear form on $H^N$, and
\[
\mcb_N(\Pi h_R,\Pi k_R)
\to\mcb_N(\Pi h,\Pi k) \quad \text{as }R\to\infty
\]
for symmetric trace-free polynomial tensors $H,K$ of degree at most $d$ satisfying \eqref{eq:strict-degree-domain-new},
with the quantitative norm bounds in \eqref{eq:strict-york-cutoff-estimates}.
Here $k_R$ and $k$ are defined from $K$ as $h_R$ and $h$ are from $H$.
The derivatives of these Hessian pairings also converge under $H(x)\mapsto H(\lambda x)$ and $K(x)\mapsto K(\lambda x)$ at $\lambda=1$, with the cutoffs dilated by the same factor.
\end{corollary}
\begin{proof}
The second variation of total $Q_{2N}$-curvature is continuous on $H^N$ by the total-order argument in \cite[Lemma 4.6]{GKW26}; differentiating the volume normalization adds only continuous lower-order terms.
Consequently
\begin{equation}\label{eq:round-hessian-continuity}
|\mcb_N(a,b)|\le C\|a\|_{H^N}\|b\|_{H^N}.
\end{equation}
Applying \eqref{eq:strict-york-cutoff-estimates} and \eqref{eq:round-hessian-continuity} to the difference of the two pairings proves the convergence.

Here dilation acts on tensor components by $H^\lambda(x):=H(\lambda x)$ for $\lambda>0$, so that $\pp_\lambda H^\lambda|_{\lambda=1}=EH$, where $E:=x\cdot\pp_x$.
The polynomial $EH$ still has degree at most $d$.
For example, we take $H_R^\lambda(x):=\chi(\lambda x/R)H(\lambda x)$ and set $h_R^\lambda:=\Omega^2H_R^\lambda$ and $h^\lambda:=\Omega^2H^\lambda$.
The additional derivative $(E\chi)(x/R)H(x)$ is supported on the cutoff annulus.
After compactification, its derivatives of order $k$ are $O(R^{d+k})$ on a spherical annulus of volume $O(R^{-n})$, so its $H^N$ norm is $O(R^{d+N-n/2})$.
Uniformly on any compact interval $0<a\le\lambda\le b<\infty$, both $h_R^\lambda-h^\lambda$ and their first $\lambda$ derivatives tend to zero in $H^N$ as $R\to\infty$.
Boundedness of $\Pi$ and $\mcb_N$ now gives convergence in $C^1$ of the quadratic forms.
If $H=\sum_jH_j$ is the decomposition into homogeneous terms, the derivative at $\lambda=1$ of the quadratic form multiplies the entry pairing $H_j$ with $H_k$ by $j+k$.
\end{proof}

\medskip
We define the scalar correction term from $\divg_{g_c}V[h]$ so that its contribution to the metric variation cancels the term proportional to $g_c$ in $DV[h]$.
For $h=\Omega^2H$ and $V[h]$ in \eqref{eq:york-projection}, we define
\begin{equation}\label{eq:scalar-correction-definition}
q[h]:=\frac{n-2N}{2n}\divg_{g_c}V[h] \quad\text{and}\quad 
\Psi[H]:=Uq[h].
\end{equation}
We write $\Psi[H]$ to display the dependence on $H$, and use $\Psi$ as an abbreviation only within proofs where $H$ is fixed.
In Euclidean coordinates, with $V=V[h]$, the same correction is $\Psi[H]=V\cdot\nabla U+\allowbreak(n-2N)U\divg_\delta V/(2n)$.
This identity follows from $\divg_{g_c}V=\divg_\delta V+nV(\log\Omega)$ and $U=\Omega^{(n-2N)/2}$.
For the scalar correction term $\Psi[H]$, we use the \emph{conformal energy space} of functions on $\R^n$
\begin{equation}\label{eq:conformal-energy-space}
\mch_N:=\{Uf:f\in H^N(\Ss^n)\} \quad\text{and}\quad 
\|Uf\|_{\mch_N}:=\|f\|_{H^N(\Ss^n)},
\end{equation}
where the factor $f$ in $Uf$ is the pullback of the spherical function $f$ to $\R^n$ by the inverse stereographic projection from $p_+$.
The norm $\|f\|_{H^N(\Ss^n)}$ is computed on the round sphere $(\Ss^n,g_c)$.

\begin{proposition}
\label{prop:scalar-corrector-unique-new}
Under \eqref{eq:strict-degree-domain-new}, $q[h_R]\to q[h]$ in $H^N$ as $R\to\infty$, at rate $O(R^{d+N-n/2})$.
The scalar correction term $\Psi[H]$ solves
\begin{equation}\label{eq:scalar-correction-equation}
L_U\Psi[H]=-DP_{2N,\delta}[H]U
\end{equation}
in the weak sense, where $L_U$ is the bubble linearization defined in \eqref{eq:bubble-linearization}.
It is the unique solution in $\mch_N$ whose round representative is orthogonal to the first spherical harmonics in $L^2(\Ss^n,g_c)$.
Its round representative is $q[h]=\Psi[H]/U$, and $\int_{\Ss^n}q[h]\,dv_{g_c}=0$.
\end{proposition}
\begin{proof}
The $H^{N+1}$ bound for $V$ in \eqref{eq:strict-york-cutoff-estimates} implies the asserted convergence.
First suppose that $H$ is smooth and compactly supported, and insert the York decomposition $h=\Pi h+\mcl_{V[h]}g_c-2(\divg_{g_c}V[h])g_c/n$ into the first variation of $P_{2N,g}1$.
Equation \eqref{eq:round-tt-first-variation} eliminates the TT term, and diffeomorphism invariance eliminates the Lie derivative term because $P_{2N,g_c}1=\kappa_0$ is constant.
Differentiating the conformal covariance law \eqref{eq:ambient-conformal-covariance} for the remaining trace term yields
\[
\left(P_{2N,g_c}-\frac{n+2N}{n-2N}\kappa_0\right)q[h]
=-DP_{2N,g_c}[h]1.
\]
Conformal covariance under stereographic projection proves \eqref{eq:scalar-correction-equation} for smooth compactly supported $H$, in particular for each approximating tensor $H_R$.
The first metric derivative on the right has order at most $2N$ in $h$, so the $H^N$ convergence allows passage to the limit in $H^{-N}(\Ss^n)$ and proves the equation for $h$ in the weak sense.
The mean of a divergence is zero.
If $Y$ is a first spherical harmonic, then $\nabla Y$ is conformal Killing, and therefore
\[
\int_{\Ss^n}q[h]Y\,dv_{g_c}
=-\frac{n-2N}{2n}\int_{\Ss^n}\langle V[h],\nabla Y\rangle\,dv_{g_c}=0.
\]
By \cite[Proposition 4.1]{GKW26}, the eigenvalue of $P_{2N,g_c}$ on spherical harmonics of degree $\ell\ge0$ is
\begin{equation}\label{eq:round-scalar-spectrum}
\prod_{j=0}^{N-1}\left(\ell+\frac{n}{2}+j\right)\left(\ell+\frac{n}{2}-j-1\right).
\end{equation}
These eigenvalues increase strictly for $\ell\ge0$, and the value at $\ell=1$ is $(n+2N)\kappa_0/(n-2N)$.
Thus the kernel of the linearized operator on the sphere consists precisely of first harmonics, proving uniqueness.
The same spectrum gives a bounded inverse from $H^{-N}$ to $H^N$ on their orthogonal complement.
\end{proof}

\medskip
For $N=1$, Brendle \cite{Brendle07} obtains the scalar correction term from a weighted divergence equation for a vector field and expresses the quadratic term in the resulting Yamabe energy through the remaining tensor.
Gong--Kim--Wei \cite[Remark 10.19]{GKW} explain the relation between the conformal Killing operator and their fourth- and sixth-order scalar correction terms.
Equation \eqref{eq:schur-york-identity-new} in Proposition \ref{thm:schur-york-identity-new} identifies the second derivative of $\mce_N(\exp(tH),U+t\Psi[H])$ at $t=0$ with $(n-2N)D^2\mscF(g_c)[\Pi h,\Pi h]/2$.
Since $U^{2n/(n-2N)}=\Omega^n$, we have $\int_{\R^n}U^{2n/(n-2N)}\,dx=|\Ss^n|$.
For smooth compactly supported metric variations, we define the energy near $(g,u)=(\delta,U)$ by
\begin{equation}\label{eq:normalized-energy-quotient}
\mce_N(g,u):=
\frac{\int_{\R^n}uP_{2N,g}u\,dv_g}
{\left(\int_{\R^n}u^{\frac{2n}{n-2N}}\,dv_g\right)^{\frac{n-2N}{n}}}.
\end{equation}
We choose the metric path $g_t:=\exp(tH)$, where the exponential is taken relative to $\delta$, so $dv_{g_t}=dx$ for trace-free $H$, and write $\mcf(H):=-DP_{2N,\delta}[H]U$.
Subscripts $g,u$ on $\mce_N$ denote derivatives at $(\delta,U)$, with metric variations in this exponential chart and additive scalar variations.
In \eqref{eq:schur-york-identity-new}, we pair functions on $\R^n$ using $dx$ and interpret pairings with distributions as the duality between $\mch_N$ and its dual.
Conformal covariance implies
\begin{equation}\label{eq:energy-curvature-identity}
\mscF\left(u^{\frac{4}{n-2N}}g\right)=\frac{2}{n-2N}\mce_N(g,u).
\end{equation}

\begin{proposition}
\label{thm:schur-york-identity-new}
Let $L_U^{-1}$ be the scalar inverse on the subspace whose round representatives are orthogonal to the first spherical harmonics, as specified in Proposition \ref{prop:scalar-corrector-unique-new}.
For a symmetric trace-free polynomial tensor $H$ of degree at most $d$ satisfying \eqref{eq:strict-degree-domain-new},
\begin{equation}\label{eq:schur-york-identity-new}
\mce_{N,gg}[H,H]
-\frac{2}{|\Ss^n|^{\frac{n-2N}{n}}}\langle\mcf(H),L_U^{-1}\mcf(H)\rangle
=\frac{n-2N}{2}\mcb_N(\Pi h,\Pi h).
\end{equation}
When the round representative $h=\Omega^2H$ is singular at $p_+$, all terms are defined by continuous extension from the smooth tensors $h_R$.
\end{proposition}
\begin{proof}
First let $H$ be smooth and compactly supported, and set $\Psi:=\Psi[H]$.
The identity $\int_{\Ss^n}q[h]\,dv_{g_c}=0$ implies $\int_{\R^n}U^{(n+2N)/(n-2N)}\Psi\,dx=0$.
Differentiating \eqref{eq:normalized-energy-quotient} gives
\[
\mce_{N,uu}[\Psi,\Psi]=\frac{2}{|\Ss^n|^{\frac{n-2N}{n}}}\langle\Psi,L_U\Psi\rangle \quad\text{and}\quad 
\mce_{N,gu}[H,\Psi]=-\frac{2}{|\Ss^n|^{\frac{n-2N}{n}}}\langle\mcf(H),\Psi\rangle.
\]
Equation \eqref{eq:scalar-correction-equation} shows that the second derivative of $\mce_N(g_t,U+t\Psi)$ at $t=0$ is the left-hand side of \eqref{eq:schur-york-identity-new}.
The derivative at $t=0$ of the metric $(U+t\Psi)^{4/(n-2N)}g_t$ is $h+4q[h]g_c/(n-2N)=\Pi h+\mcl_{V[h]}g_c$.
Criticality of $g_c$ and diffeomorphism invariance imply that $D^2\mscF(g_c)$ vanishes whenever one argument is a Lie derivative of $g_c$.
Equation \eqref{eq:energy-curvature-identity} therefore identifies this second derivative of $\mce_N$ with the right-hand side of \eqref{eq:schur-york-identity-new}.
Corollary \ref{cor:strict-hessian-transmission-new} and Proposition \ref{prop:scalar-corrector-unique-new} pass the identity to polynomial tensors by continuity.
\end{proof}

Section \ref{sec:weighted-local-model} extends the scalar correction term $\Psi[H]$ to weighted spaces beyond the energy range.
At $n=2d+2N$, the top derivative bound involves the divergent integral $\int_0^{\vrh_0}\vrh_+^{-1}\,d\vrh_+$, for any fixed $\vrh_0>0$, and does not imply $H^N$ convergence of the cutoffs.
Section \ref{sec:log-endpoint} treats this critical degree separately.

\subsection{Angular decomposition and homogeneous energy matrices} \label{sec:repair-raw-domain} 
We decompose each homogeneous tensor $H\in\mfJ_d$ by angular type and degree, then compute the projection $\Pi$ of each component.
For $n\ge5$ and an integer $d\ge2$, recall the space from \eqref{eq:established-A}:
\[
\mfJ_d=\{H_{ij}(x):H_{ij}=H_{ji},\ H\text{ polynomial and homogeneous of degree }d,
\ H_{ii}=0,\ x_iH_{ij}=0\}.
\]
In the polynomial decomposition, $(\divg_\delta H)_j:=\pp_iH_{ij}$ and $\divg_\delta V:=\pp_iV_i$ denote Euclidean divergence, with positive divergence convention.
Repeated Euclidean tensor indices are summed.
We set $\kappa_\ell:=\ell(\ell+n-2)$, the eigenvalue of $-\Delta_{\Ss^{n-1}}$ on scalar spherical harmonics of degree $\ell$.
For $x\ne0$, let $P_{ij}:=\delta_{ij}-x_ix_j/|x|^2$, the Euclidean orthogonal projection onto $x^\perp$.
For any symmetric covariant two-tensor $K=(K_{ab})$ at $x$, we define its \emph{tangential trace-free projection} by
\begin{equation}\label{eq:tangential-tracefree-projection}
(\msP K)_{ij}:=P_i{}^aP_j{}^bK_{ab}
-\frac{1}{n-1}P_{ij}P^{ab}K_{ab}.
\end{equation}
Indices in this formula are raised with the Euclidean metric.
Although $\msP$ is rational in $x$, the following expressions are polynomial tensors:
\begin{align*}
T_{\ell,j}[\msM]&:=|x|^{2j}\msM,\\
V_{\ell,j}[V]&:=|x|^{2j+2}\msP(\pp_aV_b+\pp_bV_a),\\
S_{\ell,j}[Y]&:=|x|^{2j+4}\msP(\pp_a\pp_bY) \quad\text{for }  j\ge0.
\end{align*}
Here $\msM$ is a componentwise harmonic tensor in $\mfJ_\ell$ with $\divg_\delta\msM=0$;
$V$ is a componentwise harmonic vector polynomial of degree $\ell$ with $x\cdot V=\divg_\delta V=0$; and $Y$ is a harmonic scalar polynomial of degree $\ell$.
In each case $\ell\ge2$.
We call these the \emph{tensor-derived}, \emph{vector-derived}, and \emph{scalar-derived} families, respectively,
and denote their images by $\mct_{\ell,j}$, $\mcv_{\ell,j}$, and $\mcs_{\ell,j}$.
They refine the three components in \cite[(6.16)]{GKW} by angular harmonic degree and radial power.
An \emph{algebraic Weyl tensor} is a four-tensor with the Riemann curvature symmetries, the first Bianchi identity, and zero Ricci contraction.

\begin{proposition}
\label{thm:finite-SVT-exhaustion-new}
Every $H\in\mfJ_d$ has a unique decomposition that is orthogonal after restriction to $\Ss^{n-1}$ in its $L^2$ tensor pairing:
\[
\mfJ_d=
\bigoplus_{\ell+2j=d}\mct_{\ell,j}
\ \oplus\!\bigoplus_{\ell+1+2j=d}\mcv_{\ell,j}
\ \oplus\!\bigoplus_{\ell+2+2j=d}\mcs_{\ell,j} \quad\text{for }  \ell\ge2.
\]
Moreover $\mfJ_2=\mct_{2,0}$, and this last space consists exactly of the tensors $H_{ij}=W_{ikj\ell}x^kx^\ell$ obtained from algebraic Weyl tensors $W$.
\end{proposition}

\begin{proof}
This is the orthogonal decomposition of \cite[Corollary 6.10]{GKW}, with its homogeneous components written as in \cite[Subsection D.1]{GKW}.
Writing $q$ for the integer index in that decomposition, the tensor, vector, and scalar formulas there coincide with the present ones under the index changes
\[
(\ell,j)=(d-2q,q) ,\quad 
(\ell,j)=(d-2q+1,q-1) ,\quad\text{and}\quad 
(\ell,j)=(d-2q,q-1),
\]
respectively.
The projection in \cite[(D.1)]{GKW} is \eqref{eq:tangential-tracefree-projection}, and the inner product is the same spherical $L^2$ pairing.
The degree-two assertion follows from the normal-coordinate curvature identity: a quadratic normal metric term is determined by its linearized curvature, whose Ricci contraction vanishes when the tensor is harmonic, divergence-free, and trace-free; see \cite[Section 4]{KMS}.
\end{proof}

For polynomial tensors satisfying \eqref{eq:strict-degree-domain-new}, the tensor-derived compactification is already TT, so $\Pi h=h$ and $q[h]=0$.
For a vector-derived tensor $h$, $D^*h$ is a radial multiple of a divergence-free tangential angular vector field and is therefore divergence-free on $\Ss^n$.
The round operator $D^*D$ preserves gradients; self-adjointness implies that its generalized inverse preserves the orthogonal complement of gradients.
Consequently $V[h]$ is divergence-free, and \eqref{eq:scalar-correction-definition} again gives $q[h]=0$.
These statements hold after multiplication of the polynomial tensors by smooth radial cutoff functions, and Proposition \ref{thm:strict-york-transmission-new} permits passage to the limit as $R\to\infty$.
The vector- and scalar-derived TT projections generally acquire radial components, as shown by \eqref{eq:reconstruction-ansatz}.

The projection $\Pi$ and Euclidean dilation commute with the orthogonal group $O(n)$, and $D^2\mscF(g_c)$ is invariant under this group.
For $n\ge5$, the angular harmonic representations indexed by tensor, vector, or scalar type and by degree $\ell$ are pairwise inequivalent.
Schur's lemma therefore eliminates pairings between distinct types or angular degrees.
Within a fixed angular representation, different radial powers can have nonzero pairings under $D^2\mscF(g_c)$, so all mixed pairings between them must be retained.
The family of polynomial symmetric two-tensors in \cite[(5.1)]{GKW26} is the tensor-derived family with $\ell=2$.

\medskip
For $n\ge4$, let $g=\exp H$ be in conformal normal coordinates at $\xi$, normalized by $\det g=1+O(|x|^{D+1})$ for an integer $D\ge2$,
and write $H=\sum_{d=2}^D H_d+O(|x|^{D+1})$ for its homogeneous Taylor expansion.
Hebey--Vaugon \cite{HebeyVaugon} proved that
\begin{equation}\label{eq:metric-weyl-jets}
H_2=\cdots=H_D=0
\quad\Longleftrightarrow\quad
\nabla^jW_g(\xi)=0\quad(0\le j\le D-2).
\end{equation}
For each integer $2\le d\le D$, the components of $\nabla^{d-2}W_g(\xi)$ are linear combinations of the coefficients of $H_d$, plus products of coefficients from $H_2,\ldots,H_{d-1}$.
Each product contains at least two factors, and the homogeneous degrees of the metric terms supplying these factors sum to $d$.
For example, when $D\ge4$, $\nabla^2W_g(\xi)$ is linear in the coefficients of $H_4$ plus terms quadratic in those of $H_2$.
Along the sequence in Corollary \ref{cor:ac-jet-vanishing}, estimate \eqref{eq:ac-display-5-2} therefore bounds these products by $O(\mu_i^{n-2N-d})$ for $2\le d\le(n-2N)/2$.
They decay faster than $\mu_i^{(n-2N)/2-d}$ below the critical degree and faster than $(\log(1/\mu_i))^{-1/2}$ at the critical degree.
Thus \eqref{eq:ac-display-5-2}--\eqref{eq:ac-display-5-3} also bound $\nabla^{d-2}W$ at the blow-up centers, computed in the same conformal-normal representatives; compare \cite[proof of Theorem 6.1]{KMS}.

\medskip
For $1\le j\le D-2$, under a smooth conformal change of $g$, the $j$th covariant derivative of the new Weyl tensor at $\xi$ is a nonzero multiple of $\nabla^jW_g(\xi)$ whenever $\nabla^rW_g(\xi)=0$ for $0\le r<j$.
Together with conformal covariance of $W$, induction on $j$ shows that the simultaneous vanishing in \eqref{eq:metric-weyl-jets} is independent of the representative of $[g]$.

\medskip
We next compute $\Pi$ within each angular block.
For the matrices representing $\Pi$ in the remainder of this subsection, we assume $n\ge5$ and $n>2N$, and take $D$ to be an integer satisfying $D<(n-2N)/2$, except where a weaker condition is explicitly stated.

We define the space of TT projections of the homogeneous metric terms by
\begin{equation}\label{eq:repair-domain}
\mct_D^{\mhom}
:=\left\{\Pi\left(\Omega^2\sum_{d=2}^D H_d\right):
H_d\in\mfJ_d\right\}.
\end{equation}
The superscript $\mhom$ refers to the original homogeneous tensors $H_d$, which index the projected space.
Although $H_d(0)=\pp H_d(0)=0$, the projection $\Pi$ need not preserve these identities.

\begin{proposition}
\label{prop:repair-injectivity}
The map in \eqref{eq:repair-domain} is injective.
On each fixed finite-dimensional space, the $H^N$ norm of the projected tensor is equivalent to a coefficient norm of $\sum_{d=2}^D H_d$.
\end{proposition}
\begin{proof}
It suffices to show that the normal-coordinate conditions exclude a nonzero homogeneous tensor in the image of the Euclidean conformal Killing operator $\mcd_\delta$.
Suppose $\Pi(\Omega^2H)=0$.
The vector field $V[h]$ is smooth near the stereographic origin by elliptic regularity.
Conformal covariance of \eqref{eq:conformal-killing-operator} gives $\mcd_{g_c}V=\Omega^2\mcd_\delta V$, hence $H=\mcd_\delta V$ there.
Taking homogeneous Taylor components, $H_d=\mcd_\delta V_{d+1}$.
We show that a homogeneous vector $V=V_{d+1}$ satisfying $x^i(\mcd_\delta V)_{ij}=0$ must vanish.
We set $r:=|x|$, $f:=x\cdot V$, and $b:=\divg V$.
Euler's identity gives
\[
dV=-\nabla f+\frac{2}{n} xb ,\quad\text{where } 
f=\frac{r^2b}{n(d+1)} \quad\text{and}\quad 
V=\frac{2xb}{n(d+1)}-\frac{r^2\nabla b}{nd(d+1)}.
\]
Taking divergence yields
\begin{equation}\label{eq:repair-no-normal-gauge}
r^2\Delta b=-(n-2)d(d-1)b.
\end{equation}
The Fischer eigenvalues of $r^2\Delta$ on degree-$d$ polynomials are $2a(2d-2a+n-2)\ge0$, $0\le a\le\lfloor d/2\rfloor$, whereas the coefficient in \eqref{eq:repair-no-normal-gauge} is negative.
Hence $b=0$, then $V=0$, and therefore every $H_d$ vanishes.
Injectivity and finite dimensionality give the claimed norm equivalence.
\end{proof}

The reconstruction in Proposition \ref{prop:weighted-reconstruction} relates the tensor operator $\mca$ in \eqref{eq:round-tensor-operator} to a scalar second-order operator $A_\ell$.
We express $A_\ell$ in a radial basis indexed by homogeneous polynomial degrees and compute the matrices representing $\Pi$ by triangular systems.
We set $z:=(1-r^2)/(1+r^2)$, where $r:=|x|$; thus $z=-Y_{p_+}$ for the height function in \eqref{eq:height-curvature-functional}.
For scalar functions of $z$, we write $\langle f,g\rangle_{w_\ell}:=\int_{-1}^1 f(z)g(z)w_\ell(z)\,dz$.
For $j\ge0$, we define
\begin{equation}\label{alg:eq:cyclic}
\Phi_j^{(\ell)}:=\frac{(1-z)^j}{(1+z)^{j+\ell}} \quad\text{and}\quad 
w_\ell:=(1-z^2)^{\frac{n}{2}+\ell-1} \quad\text{for }  \ell\ge2.
\end{equation}
We suppress the superscript $(\ell)$ on $\Phi_j^{(\ell)}$ while $\ell$ is fixed.
We define
\begin{equation}\label{eq:radial-operator}
A_\ell:=-(1-z^2)\pp_z^2+(n+2\ell)z\pp_z+\ell(\ell+n-1).
\end{equation}
Direct differentiation gives
\begin{equation}\label{eq:repair-A-chain}
\begin{aligned}
A_\ell\Phi_j={}&-\frac{j(n+2j+2\ell-2)}{2}\Phi_{j-1}
+\left[\frac{\ell(n-2)}{2}-2j(j+\ell)\right]\Phi_j\\
&+\frac{(j+\ell)(n-2j-2)}{2}\Phi_{j+1}.
\end{aligned}
\end{equation}
Integration with the weight $w_\ell$ gives the normalized moments and their ratio:
\begin{equation}\label{eq:repair-moments}
\mu_r^{(\ell)}:=\frac{\langle\Phi_i,\Phi_j\rangle_{w_\ell}}
{\langle\Phi_0,\Phi_0\rangle_{w_\ell}},\quad r=i+j,\qquad
\frac{\mu_{r+1}^{(\ell)}}{\mu_r^{(\ell)}}=\frac{n+2\ell+2r}{n-2\ell-2-2r}.
\end{equation}
The moment $\mu_r^{(\ell)}$ is finite for $n>2\ell+2r$, and its displayed ratio requires $n>2\ell+2r+2$.
For \eqref{eq:repair-moments}, substitute $v=(1-z)/2$; the beta parameters are $n/2+\ell+r$ and $n/2-\ell-r$, whose sum is $n$.

The index $\sigma\in\{T,V,S\}$ records whether the angular tensor is tensor-, vector-, or scalar-derived.
Proposition \ref{prop:weighted-reconstruction} expresses the tensor inner products through $c_{\sigma,\ell}W_\sigma(A_\ell)$, with the positive angular constants $c_{\sigma,\ell}$ defined in \eqref{eq:reconstruction-intertwining-norm} and
\begin{equation}\label{eq:repair-norm-operators}
W_T(x):=1 ,\quad  W_V(x):=x ,\quad\text{and}\quad 
W_S(x):=x(x+n-2).
\end{equation}
The operators $W_\sigma(A_\ell)$ determine the round $L^2$ tensor norm;
the factors $\wtbeta_\sigma$ in Proposition \ref{prop:weighted-height} arise when the tensor pairing includes multiplication by the stereographic height $z$.
Appendix \ref{app:weighted-reconstruction} proves both identities and the relation between $A_\ell$ and the round tensor operator $\mca$ under reconstruction.

The offsets $s_T:=0$, $s_V:=1$, and $s_S:=2$ record the extra homogeneous degrees in the three tensor constructions.
Fix $\sigma\in\{T,V,S\}$, $\ell\ge2$, and a nonnegative integer $M$ with $\ell+s_\sigma+2M\le D$.
Here $M$ is the largest retained radial index, distinct from the harmonic tensor $\msM$; we use the full basis $\Phi_0,\ldots,\Phi_M$ defined in \eqref{alg:eq:cyclic}.

\begin{proposition}
\label{prop:repair-C}
Fix a nonzero harmonic tensor $\msM$, vector polynomial $V$, or scalar polynomial $Y$ from the corresponding family in Proposition \ref{thm:finite-SVT-exhaustion-new}, and keep it fixed as $k$ varies.
For each $0\le k\le M$, define the homogeneous tensor of degree $d_k:=\ell+s_\sigma+2k$ by
\[
H_k:=\begin{cases}
T_{\ell,k}[\msM] &\text{for }\sigma=T,\\
V_{\ell,k}[V] &\text{for }\sigma=V,\\
S_{\ell,k}[Y] &\text{for }\sigma=S.
\end{cases}
\]
In the reconstruction \eqref{eq:reconstruction-ansatz} using this fixed angular datum, the scalar radial profile of $\Pi(\Omega^2H_k)$ lies in $\operatorname{span}\{\Phi_0,\ldots,\Phi_k\}$ and has a nonzero $\Phi_k$ coefficient.
Normalize this coefficient to one, and let the resulting coefficients in $\Phi_0,\ldots,\Phi_M$ form the $k$th column of an $(M+1)\times(M+1)$ matrix $C_\sigma$, with zero entries for $j>k$.
Then $C_\sigma$ is upper triangular with unit diagonal, $C_T$ is the identity matrix, and
\begin{equation}\label{eq:repair-CV}
(C_V)_{jk}=\prod_{r=j+1}^k
\frac{r(n+2\ell+2r-2)}{(\ell+r-1)(n-2r)} \quad\text{for }  j\le k.
\end{equation}
For $C_S=(p_{jk})$, put
\[
\mfd_j:=-\frac{j}{2}(n+2\ell+2j-2)\quad\text{and}\quad
u_j:=\frac{(j+\ell)(n-2j-2)}{2}.
\]
For each column $k$, set $p_{kk}:=1$, $y_k:=u_k$, and, for $j=k-1,k-2,\ldots,0$, define
\begin{equation}\label{eq:repair-CS}
y_j:=-\frac{\mfd_{j+1}}{u_{j+1}}y_{j+1} \quad\text{and}\quad 
p_{jk}:=\frac{y_j-\mfd_{j+1}p_{j+1,k}}{u_j}.
\end{equation}
All displayed denominators are nonzero under the degree condition $n>2N+2(\ell+s_\sigma+2M)$.
\end{proposition}
\begin{proof}
We characterize each projected column by pairing with reconstructed TT tensors.
The adjoint equations reduce to one bidiagonal system in the vector sector and two successive bidiagonal systems in the scalar sector.
Let $t:=1-z^2$ and $\mu_\ell:=(\ell-1)(\ell+n-1)$.
For $\sigma=V$, the tensor $H_k=V_{\ell,k}[V]$ has only the tangential coefficient $r^d$, with $d:=\ell+1+2k$.
Let $u(z)$ be the scalar radial profile of $\Pi(\Omega^2H_k)$ in \eqref{eq:reconstruction-ansatz}.
Pairing with a reconstructed TT tensor and integrating by parts in the round radial variable gives
\[
A_\ell u=-\mu_\ell d\,
\frac{(1-z)^k}{(1+z)^{k+\ell+1}}.
\]
The radial coefficient on the right follows from the identity
\[
t^{-1}\bigl(t\pp_z-(\ell-1)z\bigr)
\left((1-z)^{k+1}(1+z)^{-\ell-k}\right)
=-d\frac{(1-z)^k}{(1+z)^{k+\ell+1}}.
\]
We set $x:=r^2$ and let $\msq(x)$ be a polynomial.
Since $\Phi_j=2^{-\ell}(1+x)^\ell x^j$, \eqref{eq:repair-A-chain} is equivalent to
\[
A_\ell\bigl(2^{-\ell}(1+x)^\ell\msq\bigr)
=2^{-\ell}(1+x)^{\ell+1}\mcb_\ell\msq,
\]
where the bidiagonal operator is defined by
\[
\mcb_\ell x^j:=\mfd_jx^{j-1}+u_jx^j.
\]
Solving this bidiagonal equation for $j=k-1,k-2,\ldots,0$ gives $(C_V)_{jk}=-\mfd_{j+1}(C_V)_{j+1,k}/u_j$, hence \eqref{eq:repair-CV}.

For $\sigma=S$, take $H_k=S_{\ell,k}[Y]$, let $d:=\ell+2+2k$, and use $\kappa_\ell=\ell(\ell+n-2)$ and $\beta_\ell:=(n-2)(\kappa_\ell-n+1)/(n-1)$.
Using the mutually adjoint first-order operators $P$ and $P^\dagger$ in \eqref{eq:reconstruction-adjoint-operators}, integration by parts in the scalar pairing yields the source $t^{-2}((P^\dagger)^2-\kappa_\ell/(n-1))f$,
where $f:=(1-z)^{k+2}(1+z)^{-\ell-k}$ and $u$ now denotes the scalar reconstruction profile of this projection.
Here $P^\dagger f=-df$, so
\begin{equation}\label{eq:repair-scalar-source}
A_\ell(A_\ell+n-2)u
=\kappa_\ell\beta_\ell
\left(d^2-\frac{\kappa_\ell}{n-1}\right)
\frac{(1-z)^k}{(1+z)^{k+\ell+2}}.
\end{equation}
Direct operator multiplication gives
\[
A_\ell(A_\ell+n-2)
\bigl(2^{-\ell}(1+x)^\ell\msq\bigr)
=2^{-\ell}(1+x)^{\ell+2}\mcc_\ell\mcb_\ell\msq,
\]
where
\[
\mcc_\ell x^j:=\mfd_jx^{j-1}+u_{j+1}x^j.
\]
Solving the two bidiagonal systems successively gives \eqref{eq:repair-CS}.
Before imposing $p_{kk}=1$, the $k$th scalar projection column equals the column in \eqref{eq:repair-CS} multiplied by
\[
\frac{\kappa_\ell\beta_\ell\left(d^2-\frac{\kappa_\ell}{n-1}\right)}
{(k+\ell)(k+\ell+1)(n-2k-2)(n-2k-4)}.
\]
It is nonzero: $d\ge\ell+2$ and $\kappa_\ell/(n-1)\le\ell^2$ for $n\ge3$, while the other factors are strictly positive.
Before imposing a unit diagonal, the corresponding vector column multiplier is $-\mu_\ell d/[(k+\ell)(n-2k-2)]$ and is likewise nonzero.
Adjointness shows that these sources characterize the orthogonal projection $\Pi$, since the norm operators $W_\sigma(A_\ell)$ have trivial kernels on the angular components under consideration.
Proposition \ref{thm:strict-york-transmission-new} justifies the integrations by parts through cutoff convergence.
\end{proof}

\medskip
For the fixed angular datum in Proposition \ref{prop:repair-C}, the columns of $C_\sigma$ are the normalized scalar radial profiles of $\Pi(\Omega^2H_k)$ in the basis $\Phi_0,\ldots,\Phi_M$.
With $p_{N,n}$ from \eqref{eq:ambient-hessian-coefficient}, we define the energy matrix and its derivative under $H(x)\mapsto H(e^\tau x)$ by
\begin{align}
E_{ij}^{\Phi}&:=\langle\Phi_i,W_\sigma(A_\ell)p_{N,n}(A_\ell)\Phi_j\rangle_{w_\ell}\quad\text{and}\quad
E^{\mhom}:=C_\sigma^{\msT}E^\Phi C_\sigma,\label{eq:repair-E}\\
D_\sigma&:=\diag(\ell+s_\sigma+2j)_{j=0}^M\quad\text{and}\quad
S^{\mhom}:=D_\sigma E^{\mhom}+E^{\mhom}D_\sigma.
\label{eq:repair-S}
\end{align}
We define the Pohozaev matrix with the sign convention
\begin{equation}\label{eq:repair-Hplus}
H^+:=-S^{\mhom}.
\end{equation}
Let $\mcr_\sigma$ be the map in \eqref{eq:reconstruction-ansatz} that reconstructs a TT tensor from its scalar radial function.
For a unit $L^2$ angular harmonic and a coefficient vector $a$ with $\Pi h=\mcr_\sigma\sum_j(C_\sigma a)_j\Phi_j$, equations \eqref{eq:round-tt-hessian} and \eqref{eq:reconstruction-intertwining-norm} give
\begin{equation}\label{eq:radial-hessian-normalization}
D^2\mscF(g_c)[\Pi h,\Pi h]
=|\Ss^n|^{-\frac{n-2N}{n}}c_{\sigma,\ell}\,a^{\msT}E^{\mhom}a.
\end{equation}
Thus $E^{\mhom}$ omits exactly the positive volume factor and angular constant displayed here.
Independent nonzero column rescalings commute with $D_\sigma$ and act on \eqref{eq:repair-S} by congruence.

\begin{proposition}
\label{prop:repair-dilation}
The matrix $S^{\mhom}$ in \eqref{eq:repair-S} represents the derivative at $\tau=0$ of the quadratic form with matrix $E^{\mhom}$ under $H(x)\mapsto H(e^\tau x)$.
The generator of this dilation in the $\Phi$ basis is $C_\sigma D_\sigma C_\sigma^{-1}$.
\end{proposition}
\begin{proof}
Let $\vph_\tau$ be the round conformal diffeomorphism induced in stereographic coordinates by $x\mapsto e^\tau x$,
let $\omega_\tau$ satisfy $\vph_\tau^*g_c=e^{2\omega_\tau}g_c$, and define $T_\tau h:=e^{-2\omega_\tau}\vph_\tau^*h$.
Then $T_\tau(\Omega^2H)=\Omega^2H(e^\tau x)$.
Conformal covariance gives
\[
T_\tau\mcd_{g_c}V=\mcd_{g_c}(\vph_\tau^*V) \quad\text{and}\quad  \Pi T_\tau(I-\Pi)=0.
\]
Thus $\Pi T_\tau h$ depends only on $\Pi h$.
The infinitesimal generator of $T_\tau$ on trace-free tensors is
\[
\mcl_X-2z=\nabla_X
=\frac{1}{2}[\mca,z]-\frac{n}{2}z ,\quad\text{where }  X:=x\cdot\pp_x=-\nabla z.
\]
Here $I$ is the identity endomorphism, $[\mca,z]$ is the commutator with multiplication by $z$, and $\nabla X=zI$.
The last equality is the product rule for $\nabla^*\nabla$.
On the coefficients of the original homogeneous tensors the same action is exactly $\exp(\tau D_\sigma)$.
Differentiating $\exp(\tau D_\sigma)\allowbreak E^{\mhom}\allowbreak\exp(\tau D_\sigma)$ at $\tau=0$ proves \eqref{eq:repair-S}.
Multiplying the right-hand side of \eqref{eq:radial-hessian-normalization} by $(n-2N)/2$ gives the second derivative at $t=0$ of $\mce_N(\exp(tH),U+t\Psi[H])$, by \eqref{eq:schur-york-identity-new}.
Corollary \ref{cor:strict-hessian-transmission-new} justifies differentiation for compactified polynomial tensors satisfying \eqref{eq:strict-degree-domain-new}.
Conjugating the coefficient action by $C_\sigma$ proves the statement in the $\Phi$ basis.
\end{proof}

In \eqref{eq:repair-E}, the polynomial $W_\sigma(A_\ell)p_{N,n}(A_\ell)$ has degree $N+s_\sigma$ and must be applied before taking the pairing or restricting the indices.
Terms up to $\Phi_{j+N+s_\sigma}$ contribute, so the beta-integral convergence condition is
\[
n>2N+2(\ell+s_\sigma+2M).
\]

\section{Coercivity of the Pohozaev quadratic forms at all orders}\label{sec:all-order}
Equations \eqref{eq:round-tt-hessian} and \eqref{eq:schur-york-identity-new} evaluate the second variation of $\mce_N$ along the scalar correction term $\Psi[H]$ through the operator $p_{N,n}(\mca)$.
Differential identities between angular degrees reduce the entries of the associated Pohozaev matrices to radial pairings for $A_2$.
We construct all these matrices from the same scalar pairings and compute their $LDL^{\msT}$ pivots by one exact recurrence, including for $N=1,2,3$.
The pivot signs prove positive definiteness throughout \eqref{eq:main-range} and identify the first integer dimension where it fails; the determinant zeros locate the three large-order crossings.

\subsection{Angular reduction and rational formulas}\label{alg:sec:moment-functional}
We use the Pohozaev form $H^+$ in \eqref{eq:repair-Hplus}, with the matrices representing $\Pi$ and positive normalizations specified in Section \ref{sec:repair-raw-domain}.

\begin{theorem}\label{alg:thm:main}
Let $N\ge1$ and $n>2N$ be integers in the dimension range \eqref{eq:main-range}.
Every nonempty angular block
\[
(\sigma,\ell,M) ,\quad\text{where } 
2(\ell+s_\sigma+2M)<n-2N,
\]
is positive definite in the orientation \eqref{eq:repair-Hplus}.
Consequently the form is positive definite on the projected space
\[
\Pi\left(\Omega^2
\bigoplus_{d=2}^{\lfloor\frac{n-2N-1}{2}\rfloor}\mfJ_d\right).
\]
\end{theorem}

\begin{proposition}\label{alg:thm:uniform}
Let $N\ge1$ and $n$ be integers with $5\le n-2N\le22$.
Every block satisfying the degree inequality in Theorem \ref{alg:thm:main} is positive definite for the following orders:
\begin{equation}\label{alg:eq:order-intervals}
\begin{array}{c|cccccc}
n-2N&5,\ldots,16&17&18&19&20&21,22\\\hline
N&1,2,\ldots&1,\ldots,17&1,\ldots,8&1,\ldots,5&1,2,3&1.
\end{array}
\end{equation}
\end{proposition}

For fixed angular type $\sigma$ and degree $\ell$, a block is indexed by the homogeneous degrees $d_j=\ell+s_\sigma+2j$, $0\le j\le M$.
Taking the largest admissible $M$ gives the \emph{maximal} block; smaller ranges of radial indices give leading principal submatrices.

Starting from $B_{0j}=m_j$, the recurrence \eqref{alg:eq:pairrec} determines the pairing matrix.
The polynomial $G_{\sigma,\ell}$ and the matrices $U_{\sigma,\ell}$ then give $E^{\mhom}$ through \eqref{alg:eq:universalenergy}, followed by $H^+=-(D_\sigma E^{\mhom}+E^{\mhom}D_\sigma)$.
For the rational calculation we use the measure $d\mu$ in \eqref{alg:eq:normalized-pairing} and divide by $C_{N,n}$ in \eqref{alg:eq:scalar-pairings}; these positive normalizations preserve positive definiteness.
For a block $\mathsf H$ computed with these two normalizations, set $\mathsf H^{(0)}:=\mathsf H$ and compute its $LDL^{\msT}$ pivots by
\begin{equation}\label{alg:eq:ldl-pivots}
\pi_r:=\mathsf H^{(r)}_{rr},\qquad
\mathsf H^{(r+1)}_{ij}:=\mathsf H^{(r)}_{ij}
-\frac{\mathsf H^{(r)}_{ir}\mathsf H^{(r)}_{rj}}{\pi_r}
\quad (i,j>r).
\end{equation}
Here $0\le r\le M$, and the update is made whenever $r<M$ and $\pi_r\ne0$.
Writing $\Delta_r$ for the leading principal minor of size $r$ and setting $\Delta_0:=1$, we have $\pi_r=\Delta_{r+1}/\Delta_r$ whenever the preceding pivots are nonzero.
Thus positivity of all pivots is equivalent to positive definiteness.
For every fixed integer offset $n-2N\in\{5,\ldots,22\}$, all matrices and pivots are computed over $\Q(N)$ before any integer specialization.

The pairings $m_0,\ldots,m_8$ determine $218$ maximal rational matrix families, of size at most five, and their $430$ symbolic pivots.
The ancillary program verifies their positivity for the orders specified by \eqref{alg:eq:order-intervals}, as recorded in Section \ref{sec:matrix-positivity}.
This includes $N=1,2,3$ in the same construction and positivity proof.
The pairings through $m_8$ also give the large-order crossings in Theorem \ref{thm:three-sector-crossing}.
Section \ref{sec:log-endpoint} treats the critical homogeneous degree by its logarithmic coefficient.

\begin{lemma}\label{alg:lem:ladder}
For $\ell\ge2$, let $\pp:=\pp_z$ and $\lambda_\ell:=\ell(\ell+n-1)$.
Then
\begin{equation}\label{alg:eq:ladder}
\pp A_\ell=A_{\ell+1}\pp \quad\text{and}\quad 
\pp^*\pp=A_\ell-\lambda_\ell,
\end{equation}
where the adjoint is from the $w_{\ell+1}$ pairing to the $w_\ell$ pairing.
Moreover
\begin{equation}\label{alg:eq:ladderbasis}
\pp\Phi_j^{(\ell)}
=-j\Phi_{j-1}^{(\ell+1)}-(j+\ell)\Phi_j^{(\ell+1)}.
\end{equation}
Thus $\pp$ is an isomorphism from $\Span\{\Phi_j^{(\ell)}:0\le j\le M\}$ to $\Span\{\Phi_j^{(\ell+1)}:0\le j\le M\}$.
\end{lemma}

\begin{proof}
Direct differentiation proves \eqref{alg:eq:ladder}--\eqref{alg:eq:ladderbasis}, using $\pp^*=-(1-z^2)\pp_z+(n+2\ell)z$.
The matrix in \eqref{alg:eq:ladderbasis} is upper triangular with nonzero diagonal $-(j+\ell)$.
\end{proof}

We write $J_\ell$ for the matrix in \eqref{alg:eq:ladderbasis} and define
\[
U_{\sigma,\ell}:=(J_{\ell-1}\cdots J_2)^{-1}C_\sigma.
\]
The empty product is the identity, so $U_{\sigma,2}=C_\sigma$.
To invert $J_r$, solve $J_rx=y$ by the downward recurrence
\[
x_j=-\frac{y_j+(j+1)x_{j+1}}{j+r} ,\quad\text{where }  x_{M+1}=0,
\]
for $r=\ell-1,\ldots,2$.

\begin{proposition}\label{alg:prop:onefunctional}
We define
\begin{equation}\label{alg:eq:G}
G_{\sigma,\ell}(x):=W_\sigma(x)
\prod_{r=2}^{\ell-1}(x-r(r+n-1)).
\end{equation}
The bilinear form with matrix $E^{\mhom}$ in \eqref{eq:repair-E} is
\begin{equation}\label{alg:eq:universalenergy}
\langle U_{\sigma,\ell}a,
G_{\sigma,\ell}(A_2)p_{N,n}(A_2)U_{\sigma,\ell}b\rangle_{w_2},
\end{equation}
Here $a,b$ are coefficient vectors in the original homogeneous basis, and $U_{\sigma,\ell}a,U_{\sigma,\ell}b$ are interpreted as functions through their expansions in the $\Phi_j^{(2)}$ basis.
Differentiating this form under $H(x)\mapsto H(e^\tau x)$ uses the degree matrix $\diag(\ell+s_\sigma+2j)$.
\end{proposition}
\begin{proof}
The first identity in \eqref{alg:eq:ladder} moves every polynomial in $A_{r+1}$ through $\pp$, replacing $A_{r+1}$ by $A_r$.
Integration by parts with \eqref{alg:eq:ladder} then inserts $A_r-\lambda_r$.
Repeating this $\ell-2$ times proves \eqref{alg:eq:universalenergy}.
The endpoint terms at $z=-1$ vanish because their pole orders are bounded by the paired degree in \eqref{alg:eq:momentbound} below and $n>2N+2d_M$.
At $z=1$, the functions are regular and the weight vanishes to the required order.
These bounds also justify the integrations by parts as limits of integrals over truncated intervals.
\end{proof}

The Pohozaev matrix is obtained from \eqref{alg:eq:universalenergy} by multiplying its $ij$ entry by $-(d_i+d_j)$, as in \eqref{eq:repair-Hplus}.
The indices here refer to the original homogeneous tensors, even after the change of basis encoded by $U_{\sigma,\ell}$.

\medskip
We next express the common radial pairing by finitely many scalar quantities.
The presence of a block with $d_M\ge2$ forces $n>2N+4$; throughout this calculation of scalar pairings we work in that range.
We write $\msv:=n/2$ to display the poles of the rational formulas below and in Appendix \ref{app:rational-moments}.
We set $t:=r^2$.
Dividing the reconstruction functions $\Phi_j^{(2)}=(1+t)^2t^j/4$ by their common factor $(1+t)^2/4$ identifies their span with polynomials in $t$.
Conjugating $A_2$ by $(1+t)^2/4$ yields the radial operator in \cite[Lemma 5.2]{GKW26}, with
\begin{equation}\label{alg:eq:triA}
At^j=a_j^0t^{j-1}+b_j^0t^j+c_j^0t^{j+1},
\end{equation}
where
\[
a_j^0:=-\frac{j(n+2j+2)}{2} ,\quad 
b_j^0:=n-2-2j(j+2) ,\quad\text{and}\quad 
c_j^0:=\frac{(j+2)(n-2j-2)}{2}.
\]
Dividing the weighted integral by $\langle\Phi_0^{(2)},\Phi_0^{(2)}\rangle_{w_2}$ yields the probability measure
\begin{equation}\label{alg:eq:normalized-pairing}
d\mu(t):=\frac{t^{\msv+1}(1+t)^{-2\msv}}{B(\msv+2,\msv-2)}\,dt,\quad
\mu_r:=\int_0^\infty t^r\,d\mu(t),\quad
\frac{\mu_{r+1}}{\mu_r}=\frac{n+4+2r}{n-6-2r}.
\end{equation}
Here $B(a,b):=\int_0^\infty t^{a-1}(1+t)^{-a-b}\,dt$ is the \emph{beta function} for $a,b>0$, and $\langle f,g\rangle_\mu:=\int_0^\infty fg\,d\mu$.
The \emph{moments} satisfy $\mu_0=1$ and are finite for nonnegative integers $r<\msv-2$; the displayed ratio is used only for $r<\msv-3$.
We define the positive normalization and the \emph{scalar pairings} by
\begin{equation}\label{alg:eq:scalar-pairings}
C_{N,n}:=\frac{N}{2}(\msv-N+1)_{N-1}(\msv+1)_{N-1}>0 \quad\text{and}\quad  m_j:=\frac{\langle1,p_{N,n}(A)t^j\rangle_\mu}{C_{N,n}}.
\end{equation}
The factor $C_{N,n}$ is the normalization in \cite[(6.7)]{GKW26}.

\begin{lemma}\label{alg:lem:budget}
Every block $(\sigma,\ell,M)$ is determined by
\begin{equation}\label{alg:eq:momentbound}
m_0,\ldots,m_K ,\quad\text{where } 
K:=2M+\ell-2+s_\sigma=d_M-2.
\end{equation}
The same reconstruction over the entire range \eqref{eq:main-range} uses at most $m_0,\ldots,m_8$.
\end{lemma}
\begin{proof}
Self-adjointness transfers the operators to one side of the pairing, so the total polynomial degree bounds the number of scalar quantities required.
The columns of $U_{\sigma,\ell}$ have polynomial degree at most $M$.
The degree of \eqref{alg:eq:G} is $\ell-2+s_\sigma$.
The coefficients $c_j^0$ in \eqref{alg:eq:triA} are nonzero for the indices used here, so $t^i$ is a degree-$i$ polynomial in $A$ applied to $1$.
Moving this polynomial to the second argument shows that the highest scalar-pairing index required is $2M+\ell-2+s_\sigma$.
The degree inequality $2d_M<n-2N$ gives
\[
n-2N>2d_M=2K+4.
\]
Thus all required scalar pairings converge and the endpoint terms in the integrations by parts vanish.
Throughout \eqref{eq:main-range}, $d_M\le10$ and $K\le8$.
\end{proof}

The scalar pairings determine the full pairing matrix by recurrence.
We set
\[
B_{ij}:=C_{N,n}^{-1}\langle t^i,p_{N,n}(A)t^j\rangle_\mu.
\]
Starting from $B_{0j}=m_j$, self-adjointness gives the recurrence
\begin{equation}\label{alg:eq:pairrec}
B_{i+1,j}=\frac{1}{c_i^0}
\bigl[a_j^0 B_{i,j-1}+(b_j^0-b_i^0)B_{ij}
+c_j^0 B_{i,j+1}-a_i^0 B_{i-1,j}\bigr].
\end{equation}
Since $c_i^0\ne0$ for the indices used here, this determines successive rows of the pairing matrix.
Terms with a negative index are omitted.
Only entries with $i+j\le K$ are formed.
Apply the full polynomial $G_{\sigma,\ell}(A)$ before using this pairing;
restricting $A$ to the retained indices first would omit terms that pass through a higher polynomial degree before contributing to the block.

\begin{proposition}\label{alg:prop:moments}
Let $N\ge1$ be an integer.
For $j=0,\ldots,8$ and $n>2N+2j+4$, there are explicit rational functions $c_j(\msv)$ and $a_{j,i}(\msv)$, independent of $N$, such that
\begin{equation}\label{alg:eq:partialfractions}
m_j(N,n)=c_j(\msv)+\sum_{i=1}^{j+2}
\frac{a_{j,i}(\msv)}{(\msv-i)^2-N^2}.
\end{equation}
They are determined by systems of sizes at most ten.
\end{proposition}
For $0\le j\le3$, \eqref{alg:eq:partialfractions} agrees with the entries $e_{0j}$ in \cite[Theorem 6.5]{GKW26} on their common integrability range.
Appendix \ref{app:rational-moments} verifies the additional polynomial identities required through $j=8$ and derives \eqref{alg:eq:partialfractions}.
For a block requiring the pairings through $m_K$, its strict degree condition gives $\msv-i>N$ for $1\le i\le K+2$.
Thus the denominators $(\msv-i)^2-N^2$ are positive throughout the required range, and the formulas apply directly to $N=1,2,3$ as well.

\subsection{Positivity and the first loss of coercivity}\label{sec:matrix-positivity}
Fix an integer offset $a:=n-2N\in\{5,\ldots,22\}$.
With $n=2N+a$, equations \eqref{alg:eq:pairrec}, \eqref{alg:eq:G}, and \eqref{alg:eq:universalenergy} give every maximal block as a matrix over $\Q(N)$.
Compute all its pivots by \eqref{alg:eq:ldl-pivots}.
This produces the $218$ rational matrix families and $430$ symbolic pivots listed in Appendix \ref{app:ancillary}.

\begin{proof}[Proof of Proposition \ref{alg:thm:uniform}]
Exact rational verification in the ancillary program establishes that every matrix entry is defined and all $430$ pivots are positive for the corresponding integer orders in \eqref{alg:eq:order-intervals}.
The recurrence \eqref{alg:eq:ldl-pivots} therefore proves positive definiteness of every maximal block at these orders, including $N=1,2,3$ whenever admissible.
Every block with a smaller range of radial indices is a leading principal submatrix of a maximal block, so it is also positive definite.
Appendix \ref{app:ancillary} describes the common construction and how to reproduce the verification; the sign checks are implemented in the program.
\end{proof}

For example, when $n-2N=6$, the degree inequality permits a single one-dimensional block, whose normalized entry is
\[
-4m_0(N,2N+6)
=\frac{2(2N+3)(2N^2+15N+37)}{15(2N+1)}.
\]
This expression is positive for every $N\ge1$.

\begin{proof}[Proof of Theorem \ref{alg:thm:main}]
The integer orders in \eqref{alg:eq:order-intervals} cover exactly the cases in \eqref{eq:main-range} that admit a homogeneous degree $d\ge2$ with $2d<n-2N$.
Proposition \ref{alg:thm:uniform} proves positivity on every maximal angular block, and all smaller ranges of radial indices follow by restriction.
For $1\le n-2N\le4$ there is no homogeneous degree $d\ge2$ with $2d<n-2N$.
\end{proof}

Corollary \ref{cor:first-pohozaev-failure} uses the scalar pairings \eqref{alg:eq:partialfractions} and pivot recurrence \eqref{alg:eq:ldl-pivots}; Appendix \ref{app:ancillary} describes their implementation.

\begin{corollary}\label{cor:first-pohozaev-failure}
For each $N\in\mathbb N$, the first integer dimension at which a matrix representing the Pohozaev quadratic form below the critical degree fails to be positive definite is the first integer above the range \eqref{eq:main-range}.
\end{corollary}
\begin{proof}
Theorem \ref{alg:thm:main} proves positivity in every smaller admissible integer dimension.
For $N\ge4$, \cite[Theorem 6.8]{GKW26} proves indefiniteness at the next dimension for the tensor-derived block at angular degree two with homogeneous degrees $2,4,6,8$.
For $N=1,2,3$, construct the tensor-derived blocks at angular degree two from \eqref{alg:eq:partialfractions} and \eqref{alg:eq:universalenergy}.
The required scalar-pairing indices are at most $6,8,6$, respectively, and the exact pivots from \eqref{alg:eq:ldl-pivots} have the following signs:
\[
\begin{array}{c|c|c|c}
N&n&\text{degrees}&\text{pivot signs}\\\hline
1&25&2,4,6,8&+,+,+,-\\
2&25&2,4,6,8,10&+,+,+,+,-\\
3&27&2,4,6,8&+,+,+,-.
\end{array}
\]
All these pairings lie in the strict integrability range, so the calculation requires no exceptional low-order formula.
\end{proof}

\medskip
When $n-2N$ is even, the homogeneous degree $(n-2N)/2$ is excluded by the integrability condition \eqref{eq:strict-degree-domain-new}.
When this degree is at least two, Proposition \ref{prop:log-energy} computes the leading logarithmic term of $D^2\mscF(g_c)$ evaluated on cutoff tensors at fixed $N,n$.
Lemma \ref{lem:ac-endpoint} shows that, after dilation and the orientation \eqref{eq:repair-Hplus},
the logarithmic coefficient is positive definite on the images under $\Pi$ of the tensors in $\mfJ_{(n-2N)/2}$.
Together with Theorem \ref{alg:thm:main}, this proves Theorem \ref{thm:intro-algebra}.

\subsection{The first crossings at large order}\label{sec:three-sector-crossing}
We now describe the large-order loss of positivity in the tensor-, vector-, and scalar-derived sectors with angular degree $\ell=2$.
We write $\sigma=T,V,S$ for these sectors and set $s_T:=0$, $s_V:=1$, and $s_S:=2$.
Fix the radial indices $j=0,1,2,3$, whose homogeneous degrees are $2+s_\sigma+2j$.
Let $H_\sigma^+(N,n)$ be the resulting four-by-four matrix in the orientation \eqref{eq:repair-Hplus}.
We normalize the columns of $C_\sigma$ in Proposition \ref{prop:repair-C} to have diagonal entries equal to one and use the probability measure $d\mu$ and positive energy factor $C_{N,n}$ in \eqref{alg:eq:normalized-pairing} and \eqref{alg:eq:scalar-pairings}, respectively.
After division of the energy by $C_{N,n}$, the matrix is specified by \eqref{alg:eq:partialfractions}, \eqref{alg:eq:pairrec}, and \eqref{alg:eq:universalenergy} at $\ell=2$.
Its entries are defined by convergent integrals when $n>2N+16+2s_\sigma$.
Real values of $N,n$ below refer to the rational continuation of this finite matrix.

To expand these matrices as $N\to\infty$, we set $\vep:=N^{-1}$, use a local real variable $t$ in $n=2N+t$, and define
\begin{equation}\label{eq:three-sector-normalization}
\whH_\sigma(\vep,t)
:=\vep^{2+s_\sigma}
H_\sigma^+(\vep^{-1},2\vep^{-1}+t).
\end{equation}
The formula extends analytically to $\vep=0$ when $t>16+2s_\sigma$.

\begin{theorem}
\label{thm:three-sector-crossing}
For each $\sigma\in\{T,V,S\}$, the determinant of $\whH_\sigma(0,t)$ has a unique simple zero $a_\sigma$ in the interval displayed below.
We define
\begin{equation}\label{eq:three-sector-b}
b_\sigma:=-\frac{\pp_\vep\det\whH_\sigma(0,a_\sigma)}
{\pp_t\det\whH_\sigma(0,a_\sigma)}.
\end{equation}
These constants satisfy $b_\sigma>0$ and have the following rounded values:
\begin{center}
\begin{tabular}{cccc}\toprule
$\sigma$&interval for $a_\sigma$&$a_\sigma$&$b_\sigma$\\\midrule
$T$&$(16,17)$&16.2018714780&13.5128796221\\
$V$&$(18,19)$&18.6086755387&22.5041888032\\
$S$&$(20,21)$&20.8572849379&33.1661674079\\\bottomrule
\end{tabular}
\end{center}
For all sufficiently large integers $N$, the determinant of $H_\sigma^+(N,n)$ has exactly one zero with $n-2N$ in a fixed neighborhood of $a_\sigma$, and this zero satisfies
\begin{equation}\label{eq:three-sector-expansion}
n_{N,\sigma}=2N+a_\sigma+\frac{b_\sigma}{N}+O(N^{-2}).
\end{equation}
As $n$ increases through the zero within that neighborhood, the matrix changes from positive definite to having three positive eigenvalues
and one negative eigenvalue; at the zero it has three positive eigenvalues and a one-dimensional kernel.
Moreover, this is the first loss of positivity in the fixed $\ell=2$ sector when all radial indices satisfying $2(2+s_\sigma+2j)<n-2N$ are included; the largest such index is $M=3$ at the crossing.
\end{theorem}

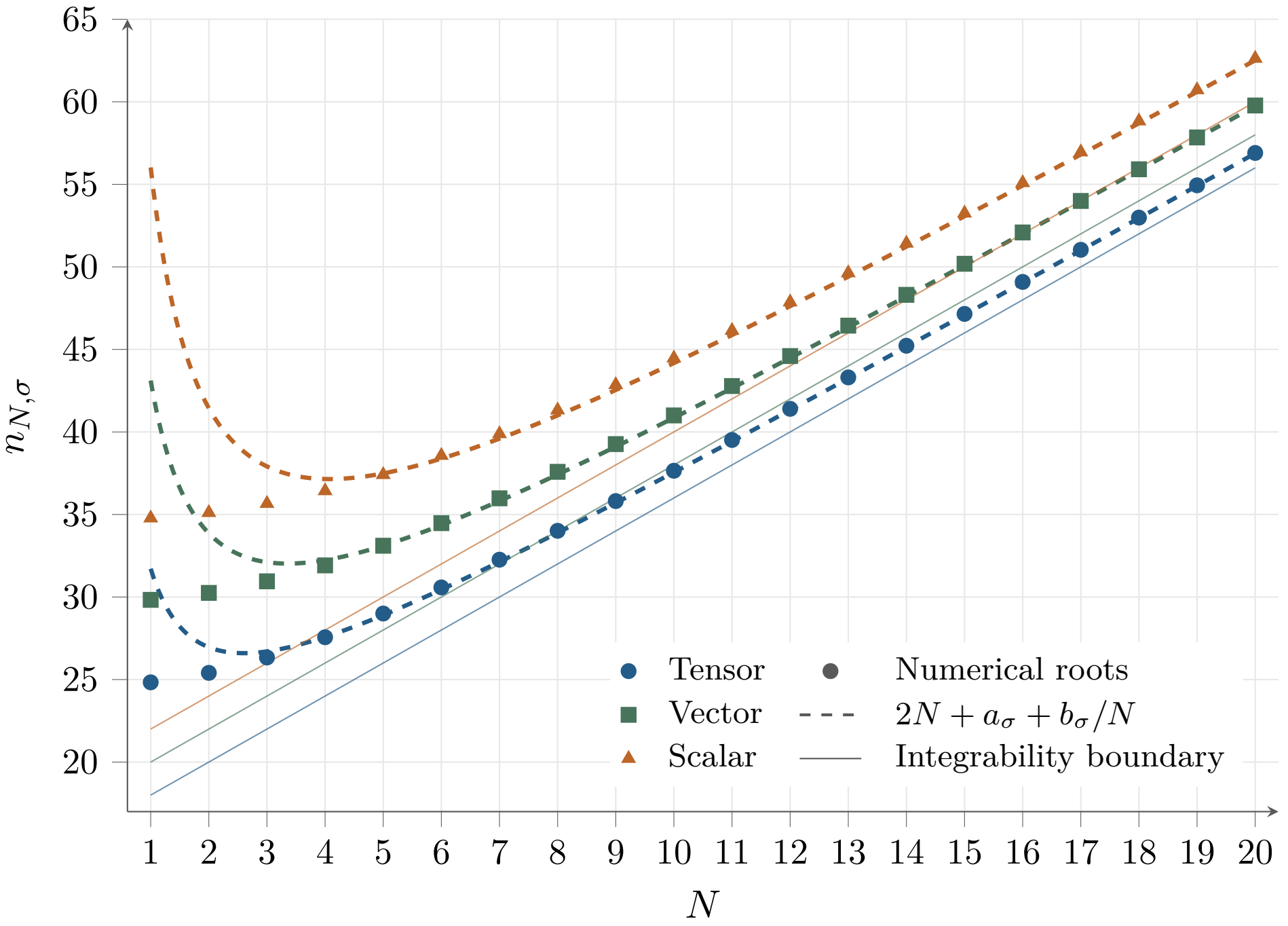
\begin{figure}[tbp]
\centering
\begin{tikzpicture}
\begin{axis}[
width=12.8cm,height=9.3cm,axis lines=left,axis line style={black!65},
tick align=outside,tick style={black!65},tick label style={font=\small},
label style={font=\normalsize},grid=major,grid style={black!9},
legend style={draw=none,fill=white,font=\footnotesize,cells={anchor=west},column sep=7pt},
legend columns=2,legend pos=south east,xlabel={$N$},ylabel={$n_{N,\sigma}$},
xmin=0.6,xmax=20.4,ymin=17,ymax=65,xtick={1,2,...,20},ytick={20,25,...,65}]
\addplot[sectorT,very thick,dashed,domain=1:20,samples=400,forget plot] {2*x+16.201871478047061988730081716198537044551933747673291029518009857+13.51287962214600311360673733738453207981115116299752974440316054/x};
\addplot[sectorT,thin,opacity=0.65,domain=1:20,samples=2,forget plot] {2*x+16};
\addplot[sectorT,only marks,mark=*,mark size=2pt,forget plot] coordinates {(1,24.83252084045323) (2,25.405566247922625) (3,26.337972340628898) (4,27.563606428075698) (5,29.0004788317551) (6,30.58181266523677) (7,32.26218000133571) (8,34.012068278645586) (9,35.81211178027384) (10,37.649268940339134) (11,39.51451342243104) (12,41.401438986594215) (13,43.30539560260797) (14,45.222939456755256) (15,47.151473615133035) (16,49.08900769287062) (17,51.03399366473955) (18,52.9852114915643) (19,54.94168801066647) (20,56.90263846045972)};
\addplot[sectorV,very thick,dashed,domain=1:20,samples=400,forget plot] {2*x+18.608675538677390379372179300525152035795341527472911482377369147+22.504188803156277224073109932937206731949003565555542726000834863/x};
\addplot[sectorV,thin,opacity=0.65,domain=1:20,samples=2,forget plot] {2*x+18};
\addplot[sectorV,only marks,mark=square*,mark size=2pt,forget plot] coordinates {(1,29.82918452105917) (2,30.249661795814685) (3,30.951260416459483) (4,31.916471345834506) (5,33.107149166662396) (6,34.47612745086279) (7,35.97999686034542) (8,37.58450875949836) (9,39.26442175731309) (10,41.00151772161804) (11,42.78262628487357) (12,44.59812186068904) (13,46.44087003662621) (14,48.30550487123983) (15,50.18793267078233) (16,52.08498797729986) (17,53.994192044986576) (18,55.913580975261624) (19,57.84158177121057) (20,59.77692178201997)};
\addplot[sectorS,very thick,dashed,domain=1:20,samples=400,forget plot] {2*x+20.857284937942843186452585432598149729335421430410689336902020748+33.166167407876854598439443659714823084478503857873117733540616806/x};
\addplot[sectorS,thin,opacity=0.65,domain=1:20,samples=2,forget plot] {2*x+20};
\addplot[sectorS,only marks,mark=triangle*,mark size=2pt,forget plot] coordinates {(1,34.782433975881204) (2,35.11038060333288) (3,35.66107653352499) (4,36.43231861438254) (5,37.41068837369959) (6,38.57185034736634) (7,39.88596429318871) (8,41.32364018713189) (9,42.859411627569) (10,44.47270459622968) (11,46.14744935414623) (12,47.87125193298429) (13,49.63456465919934) (14,51.42999824767266) (15,53.25179122796323) (16,55.0954124277313) (17,56.957265987245016) (18,58.83447252967211) (19,60.724706139204315) (20,62.62607212585374)};
\addlegendimage{sectorT,mark=*,only marks}\addlegendentry{Tensor}
\addlegendimage{black!65,mark=*,only marks}\addlegendentry{Numerical roots}
\addlegendimage{sectorV,mark=square*,only marks}\addlegendentry{Vector}
\addlegendimage{black!65,thick,dashed}\addlegendentry{$2N+a_\sigma+b_\sigma/N$}
\addlegendimage{sectorS,mark=triangle*,only marks}\addlegendentry{Scalar}
\addlegendimage{black!65,thin}\addlegendentry{Integrability boundary}
\end{axis}
\end{tikzpicture}
\caption{Crossing dimensions in the tensor, vector, and scalar sectors with $\ell=2$, for $N=1,\ldots,20$.
Markers give numerical determinant roots of the rational matrices with largest radial index $M=3$.
Dashed curves show $2N+a_\sigma+b_\sigma/N$ from Theorem \ref{thm:three-sector-crossing}; this approximation is asymptotic as $N\to\infty$.
Thin lines mark the integrability boundaries $2N+16$, $2N+18$, and $2N+20$ for the tensor, vector, and scalar sectors, respectively.
At small orders $M=3$ need not be the maximal admissible radial index; the theorem identifies these crossings with the first loss of positivity in each full $\ell=2$ sector for sufficiently large $N\in\mathbb N$.}
\label{fig:three-sector-crossings}
\end{figure}

The smallest integers above the three crossings are eventually $2N+17$, $2N+19$, and $2N+21$, respectively.
The tensor constants agree with those of the polynomial metric perturbations generated by algebraic Weyl tensors with homogeneous degrees $2,4,6,8$ in \cite[Theorem 3]{GKW26}, explaining the eventual bound in Theorem \ref{thm:ac-main}.
In \cite{GKW26}, we choose a metric perturbation profile for which the reduced energy has a nondegenerate critical point.
Placing rescaled copies in disjoint balls constructs one fixed smooth metric on $\Ss^n$, perturbed from $g_c$, for which \eqref{eq:constant-q-equation} has positive solutions with $L^\infty$ norms tending to infinity.

\begin{proof}[Proof of Theorem \ref{thm:three-sector-crossing}]
For the tensor-derived block with $\ell=2$, dividing $H^+$ in \eqref{eq:repair-Hplus} by $\langle\Phi_0^{(2)},\Phi_0^{(2)}\rangle_{w_2}$ gives $-S$ in \cite[(5.5)]{GKW26}.
Its local crossing and the constants $a_T,b_T$ therefore follow from \cite[Theorem 3]{GKW26}.
We prove the vector and scalar crossings and exclude earlier losses in all three sectors by controlling the smaller ranges of radial indices and their integrability endpoints.
Substitute $n=2N+t$ in the partial fractions \eqref{alg:eq:partialfractions}.
In these nine formulas, $c_j(\msv)=O(1)$ and $a_{j,i}(\msv)=O(\msv^3)$ as $\msv\to\infty$.
Since each denominator grows linearly in $N$ when $\msv=N+t/2$, the scalar pairings grow at most quadratically in $N$;
the norm polynomials $W_\sigma$ contribute the additional power $N^{s_\sigma}$.
The rational formulas for the radial pairings and the matrices representing $\Pi$ therefore prove the analytic extension
in \eqref{eq:three-sector-normalization}, locally uniformly with derivatives for $t>16+2s_\sigma$.
The growth bounds for $m_j$ and $W_\sigma$ also hold when the range of radial indices is reduced.
Only $m_0,\ldots,m_6$ are required in the tensor sector, $m_0,\ldots,m_7$ in the vector sector, and $m_0,\ldots,m_8$ in the scalar sector.

For the vector and scalar sectors, the limiting determinant numerators have degrees $18$ and $22$ after cancellation.
Their denominators are positive on $(18,19)$ and $(20,21)$, respectively.
Exact root counts certify one simple root in each interval and enclose it in a smaller interval with rational endpoints.
At these roots the $t$ derivatives of the determinants are negative and their $\vep$ derivatives are positive, while the first three leading principal minors are positive on these smaller intervals.
Appendix \ref{app:ancillary} specifies these root-count and sign tests.
The implicit-function theorem gives a unique local analytic solution in each of these two sectors,
\[
t=a_\sigma+b_\sigma\vep+O(\vep^2).
\]
Differentiating the determinant equation proves \eqref{eq:three-sector-b} and \eqref{eq:three-sector-expansion}.
The determinant decreases through zero while the first three leading minors remain positive, so an $LDL^{\msT}$ factorization proves the stated signs of the eigenvalues.

The largest radial index allowed by $2(2+s_\sigma+2j)<n-2N$ is $M$ exactly when
\[
4+2s_\sigma+4M<t\le8+2s_\sigma+4M.
\]
Exact root counts show that every leading minor of the limiting matrices for $M=0,1,2$ is positive throughout its interval.
For $M=3$, the leading principal minors of sizes one, two, and three are positive from the left endpoint through the smaller interval containing the root, and the determinant is positive up to that interval.
The behavior as $\vep\to0$ and $t$ approaches a left endpoint $e:=4+2s_\sigma+4M$ follows directly from the pairing formulas:
only $m_{2M+s_\sigma}$ has a pole at this endpoint of the integrability range, and only the last diagonal matrix entry contains that pairing.
Indeed, the relevant denominator factors as
\[
\left(N+\frac{t}{2}-i\right)^2-N^2=\frac{1}{4}(t-2i)(4N+t-2i).
\]
After the rescaling in \eqref{eq:three-sector-normalization}, the pole is fixed and simple, with all other factors analytic near $\vep=0,t=e$.
Writing $h_0:=t-e$, the rescaled matrix has the form
\[
\begin{pmatrix}
B(\vep,h_0)&q(\vep,h_0)\\
q(\vep,h_0)^{\msT}&\frac{r(\vep)}{h_0}+c(\vep,h_0)
\end{pmatrix} ,\quad\text{where }  B(0,0)>0 \quad\text{and}\quad  r(0)>0,
\]
where $B,q,c,r$ are analytic and bounded near $(0,0)$.
The positive residue and the leading minors of the leading principal matrix $B(0,0)$ are checked exactly for all twelve pairs $(\sigma,M)$, $0\le M\le3$; for $M=0$, that matrix is empty.
For small positive $\vep$ and $h_0$, $B$ remains positive definite and the last pivot in the $LDL^{\msT}$ factorization is
\[
\frac{r(\vep)}{h_0}+c-q^{\msT}B^{-1}q>0,
\]
with the term involving $B^{-1}$ omitted when $M=0$.
This proves positive definiteness uniformly for small $\vep>0$ as $t$ decreases to each integrability endpoint $e$.
On each remaining compact subinterval before the crossing, the positive leading minors of the limiting matrix are bounded away from zero.
Uniform convergence therefore preserves positive definiteness for sufficiently large $N\in\mathbb N$; together with the local crossing analysis, this excludes any earlier loss of positivity in the fixed sector.
Finally, $16+2s_\sigma<a_\sigma<17+2s_\sigma$ places the crossing inside the interval where $M=3$ is maximal.
\end{proof}

\section{Local corrections, Green functions, and mass}\label{sec:weighted-local-model}
Equation \eqref{eq:schur-york-identity-new} evaluates the second derivative at $t=0$ of $\mce_N(\exp(tH),U+t\Psi[H])$ through $\Pi(\Omega^2H)$, without solving the scalar equation separately.
To compare this quadratic form with a concentrating solution, we also need the decay and cutoff estimates for the scalar correction term in Corollary \ref{cor:repair-high-corrector}.

\subsection{Local scalar correction terms and critical-degree estimates}\label{sec:log-endpoint}
Fix $n>2N$ and use the standard bubble $U$ and its linearization $L_U$ from \eqref{eq:main-notation1} and \eqref{eq:bubble-linearization}.
We choose a real orthonormal basis $Y_a$, $1\le a\le n+1$, of first round spherical harmonics, pulled back by stereographic projection, and set $Z_a:=UY_a$.
Fix $0<\vth<1$ and let $X_\beta$ be the weighted $C^{2N,\vth}$ space defined by
\[
|\pp^j\phi(x)|\le C(1+|x|)^{-\beta-j} \quad\text{for }  0\le j\le2N,
\]
and by $[\pp^{2N}\phi]_{\vth,\{R<|x|<2R\}}\le CR^{-\beta-2N-\vth}$ for $R\ge1$, together with the ordinary $C^{2N,\vth}$ bound on $B_2$.
The norm $\|\phi\|_{X_\beta}$ is the least common bound in these estimates.
Let $\mcy_{\beta+2N}$ be the weighted $C^{0,\vth}$ source space defined in the same way, with decay exponent $\beta+2N$ and derivative order zero.

\begin{proposition}
\label{thm:repair-weighted-inverse}
For $0<\beta<n-2N$ and $f\in\mcy_{\beta+2N}$, the equation $L_U\phi=f$ has a solution in $X_\beta$ if and only if
\[
\int_{\R^n}f Z_a\,dx=0\quad(1\le a\le n+1).
\]
There is a unique such solution satisfying
\begin{equation}\label{eq:weighted-solution-orthogonality}
\int_{\Ss^n}\left(\frac{\phi}{U}\right)Y_a\,dv_{g_c}=0,
\end{equation}
and
\begin{equation}\label{eq:weighted-inverse-estimate}
\|\phi\|_{X_\beta}\le C\|f\|_{\mcy_{\beta+2N}}.
\end{equation}
\end{proposition}
\begin{proof}
Let $\mcl_c:=P_{2N,g_c}-(n+2N)\kappa_0/(n-2N)$.
The spectrum in \eqref{eq:round-scalar-spectrum} shows that its kernel is exactly the first harmonics.
Its inverse on their orthogonal complement has a kernel $K_c(\xi,\eta)$ with
\[
|K_c(\xi,\eta)|\le C\,d_{g_c}(\xi,\eta)^{2N-n}.
\]
The estimate follows from an elliptic parametrix, since the finite-rank projection has a smooth kernel.

Conformal covariance gives the inverse kernel in Euclidean coordinates:
\[
K_U(x,y):=U(x)U(y)K_c(\xi(x),\xi(y)).
\]
Here $\xi(x)$ denotes the point on the round sphere with stereographic coordinate $x$.
The chordal identity $|\xi(x)-\xi(y)|=\sqrt{\Omega(x)\Omega(y)}\,|x-y|$ yields $|K_U(x,y)|\le C|x-y|^{2N-n}$.
Splitting the convolution into $|y|\le|x|/2$, $|x-y|\le|x|/2$, and the remaining region proves
\[
\int_{\R^n}|x-y|^{2N-n}(1+|y|)^{-\beta-2N}\,dy
\le C_\beta(1+|x|)^{-\beta}
\]
precisely for $0<\beta<n-2N$.
The source obtained by conformal transformation to the sphere is locally integrable at the projection pole $p_+$, and its pairings against the first harmonics are $\int_{\R^n} fZ_a$.
Thus integration against $K_U$ constructs the solution satisfying \eqref{eq:weighted-solution-orthogonality} when these pairings vanish.
Scaled interior Schauder estimates on annuli give all derivative and H\"older bounds in $X_\beta$.

For uniqueness, if $L_U\phi=0$, then $\bph:=\phi/U$ solves $\mcl_c\bph=0$ off the pole and satisfies $|\nabla^j\bph|=O(\vrh_+^{\beta+2N-n-j})$.
Here $\vrh_+$ is the round distance to $p_+$, as in Section \ref{sec:scalar-complement-new}.
In a Green formula across a sphere of radius $\vrh_+$, every boundary term that could produce a distribution supported at the pole is $O(\vrh_+^\beta)$ and tends to zero as $\vrh_+\to0$.
Hence the equation holds in the weak sense on the closed sphere.
Elliptic regularity makes $\bph$ smooth, so it is a first harmonic and \eqref{eq:weighted-solution-orthogonality} eliminates it.
Finally, integration against $Z_a$, with the same vanishing boundary estimate on the coordinate spheres, proves necessity of the orthogonality conditions.
\end{proof}

\begin{corollary}
\label{cor:repair-high-corrector}
Let $H$ be a symmetric trace-free polynomial tensor of degree $d$ with $2\le d<n-2N$.
Extend $h=\Omega^2H$ as a locally integrable tensor on $\Ss^n$, and apply \eqref{eq:york-projection} in the weak sense.
Then
\[
\Psi[H]=U\frac{n-2N}{2n}\divg_{g_c}V[h]
\]
is the unique solution in $X_{n-2N-d}$ satisfying \eqref{eq:weighted-solution-orthogonality} and
\[
L_U\Psi[H]=-DP_{2N,\delta}[H]U.
\]
Let $\chi$ be smooth and radial, equal to one on $B_1$ and zero outside $B_2$, and set $H_R:=\chi(x/R)H$.
For every $0<\beta'<n-2N-d$,
\begin{equation}\label{eq:weighted-corrector-cutoff}
\|\Psi[H_R]-\Psi[H]\|_{X_{\beta'}}
\le C R^{d+2N-n+\beta'}.
\end{equation}
\end{corollary}
\begin{proof}
Since $d<n-2N<n$, the compactified tensor is locally integrable.
The map $h\mapsto V[h]$ defined by York decomposition has order $-1$ and local kernel bounded by $C\,d_{g_c}(x,y)^{1-n}$.
Splitting its integral near the projection pole, followed by scaled elliptic estimates, gives $|\nabla^jV[h]|=O(\vrh_+^{1-d-j})$ for $d>1$.
Consequently $q[h]=(n-2N)\divg_{g_c}V[h]/(2n)=O(\vrh_+^{-d})$ with differentiated bounds, and $Uq[h]\in X_{n-2N-d}$.

Conformal covariance of the trace-free conformal Killing operator and the identity $h+4q[h]g_c/(n-2N)=\Pi h+\mcl_{V[h]}g_c$ prove the asserted equation for $\Psi[H]$ first away from the pole.
The first variation of $g\mapsto P_{2N,g}1$ in the direction $\Pi h$ vanishes by \eqref{eq:round-tt-first-variation}, and its variation in the direction $\mcl_{V[h]}g_c$ vanishes because $P_{2N,g_c}1$ is constant.
Passing to the weak limit of the smooth cutoffs proves the equation on the whole sphere.
The normalization $V[h]\perp\ker D$ in \eqref{eq:york-projection} implies $\langle q[h],Y_a\rangle=-(n-2N)\langle V[h],\nabla Y_a\rangle/(2n)=0$, since $\nabla Y_a$ is a conformal Killing field.
Proposition \ref{thm:repair-weighted-inverse} proves uniqueness.

We define $f_R:=-DP_{2N,\delta}[H_R]U$ and $f:=-DP_{2N,\delta}[H]U$.
Their difference is supported where $|x|\ge R$ and obeys $|\pp^j(f_R-f)|\le C(1+|x|)^{d-n-j}$, including the cutoff annulus $R<|x|<2R$.
Thus its $\mcy_{\beta'+2N}$ norm is $O(R^{d+2N-n+\beta'})$.
Both sources are orthogonal to every $Z_a$ with $a\in\{1,\ldots,n+1\}$ by the equations for $\Psi[H_R]$ and $\Psi[H]$, respectively, and the necessity assertion in Proposition \ref{thm:repair-weighted-inverse}.
Estimate \eqref{eq:weighted-inverse-estimate} therefore proves \eqref{eq:weighted-corrector-cutoff}.
\end{proof}

At $n=2d+2N$, Proposition \ref{thm:repair-weighted-inverse} still applies with $\beta=d$, while $D^2\mscF(g_c)[\Pi h_R,\Pi h_R]$ has the leading logarithmic term computed in Proposition \ref{prop:log-energy}.

\begin{proposition}
\label{prop:log-energy}
Let $n=2N+2d$, $N\ge1$, $d>0$.
Let $h$ be a TT tensor on $\Ss^n\setminus\{p\}$, smooth away from $p$, with expansion in normal coordinates $y=r\theta$ at $p$
\[
h(y)=h_0(y)+O_{C^{N+1}}(r^{1-d}) ,\quad\text{where } 
h_0(y):=r^{-d}a(\theta) \quad\text{and}\quad  a\not\equiv0.
\]
Here the remainder notation includes the differentiated bounds of orders $0,\ldots,N+1$.
Fix a smooth function $\chi:[0,\infty)\to[0,1]$ with $\chi=0$ on $[0,1]$ and $\chi=1$ on $[2,\infty)$, and set $\chi_R(y):=\chi(Rr)$, extended by one outside the coordinate chart.
Set $t_R:=\Pi(\chi_Rh)$.
Write
\[
p_{N,n}(x)=-\hc_{N,n}x^N+O(x^{N-1}) ,\quad\text{where }  \hc_{N,n}>0.
\]
Then
\begin{equation}\label{eq:log-energy}
D^2\mscF(g_c)[t_R,t_R]
=-|\Ss^n|^{-\frac{n-2N}{n}}\hc_{N,n}\mcc(h_0)\log R+O(1),
\end{equation}
where, using ordered Cartesian derivative indices and Euclidean tensor components,
\[
\mcc(h_0):=\int_{\Ss^{n-1}}
|\pp^N h_0(\theta)|^2\,d\theta>0.
\]
The expansion holds as $R\to\infty$, with constants depending on $N,n,h$ and the cutoff profile.
\end{proposition}

\begin{proof}
For $0\le j<N$, the squared derivative of order $j$ is locally integrable: its radial exponent is $n-1-2d-2j=2(N-j)-1>-1$.
The top derivative of $h_0$ has squared radial integral
\[
\int_{\frac{2}{R}}^{r_0}r^{n-1-2d-2N}\,dr
\int_{\Ss^{n-1}}|\pp^Nh_0(\theta)|^2\,d\theta
=\mcc(h_0)\log R+O(1).
\]
Every cross term with the remainder is integrable, and the cutoff annulus $R^{-1}<r<2R^{-1}$ contributes $O(1)$ by scaling.
Replacing the Euclidean metric and derivatives by the round metric $g_c$ and its connection changes the energy by $O(1)$, since the normal-coordinate errors have positive powers of $r$.
If $\mcc(h_0)=0$, all derivatives of order $N$ vanish on the punctured ball, so every component of $h_0$ is a polynomial of degree at most $N-1$.
Its negative homogeneity forces that polynomial to vanish, a contradiction.

Extend $\mca=\nabla^*\nabla+2$ to trace-free symmetric tensors.
Integration by parts on the smooth tensor $\chi_Rh$ shows that the principal part of $\langle\chi_Rh,\mca^N\chi_Rh\rangle$ is $\|\nabla^N(\chi_Rh)\|_2^2$.
All remaining terms have differential order at most $2N-2$ and are bounded by the uniformly bounded $H^{N-1}$ norm.
Thus the energy of $\chi_Rh$ has the claimed logarithmic term before TT projection.

The weak divergence has no term supported at $p$: the boundary term on the coordinate sphere of radius $r$ is $O(r^{n-d-1})=o(1)$.
Consequently $\divg_{g_c}(\chi_Rh)$ is supported on the cutoff annulus, and
\[
\|\divg_{g_c}(\chi_Rh)\|_{H^{N-1}}\le C.
\]
This follows by applying up to $N-1$ derivatives to $(\nabla\chi_R)\mathbin{\lrcorner}h$; the largest squared norm scales as $R^{2(d+N)-n}=1$.
The elliptic estimate for $D^*D$ on $(\ker D)^\perp$ gives
\[
\|(I-\Pi)(\chi_Rh)\|_{H^N}\le C.
\]
On the round sphere $\mca$ is self-adjoint and preserves TT tensors, so it also preserves their $L^2$ orthogonal complement.
There is therefore no cross term between the two orthogonal components in the polynomial form $\langle\cdot,p_{N,n}(\mca)\cdot\rangle$.
Removing the component orthogonal to the TT tensors changes the energy by $O(1)$.
This proves \eqref{eq:log-energy}.
\end{proof}

\begin{lemma}
\label{lem:ac-endpoint}
If $n-2N$ is even and at least four, the coefficient of $\log R$ in the oriented Pohozaev form is positive definite on $\mfJ_{(n-2N)/2}$ for fixed $N,n$.
In particular, a nonzero homogeneous normal-coordinate term of degree $(n-2N)/2$ always has a nonzero leading term after applying $\Pi$.
\end{lemma}
\begin{proof}
Write $n=2N+2d$, so $d$ is the critical homogeneous degree.
In normal coordinates at the projection pole $p_+$,
the compactification of a homogeneous term in $\mfJ_d$ has the form $k=k_0+O_{C^{N+2}}(r^{2-d})$, with $k_0$ homogeneous of degree $-d$.
The remainder gains two powers of $r$ because inversion is combined with the even Taylor expansion of the round metric.
The principal part of $D^*D$ is, up to a nonzero scalar, $\mcl:=\Delta+(n-2)\nabla\divg/n$.
Substituting homogeneous vector fields reduces its principal equation to a finite linear system.
For a scalar angular harmonic $Y$ of degree $\ell\ge1$, its indicial determinant on vectors $r^\lambda(aY\theta+b\nabla_{\Ss^{n-1}}Y)$ is
\[
\frac{2(n-1)}{n}
(\lambda-\ell+1)(\lambda-\ell-1)
(\lambda+\ell+n-1)(\lambda+\ell+n-3).
\]
For tangential divergence-free vector fields, the corresponding indicial roots are $\ell,-\ell-n+2$; for radial vector fields with angular degree zero, they are $1,1-n$.
None of these indicial roots lies in $(2-n,0)$.
Since $1-d\in(2-n,0)$, the leading homogeneous term $V_0$ of $V[k]$ exists uniquely at degree $1-d$; hence $h_0:=k_0-\mcd_\delta V_0$ is a homogeneous Euclidean TT tensor.

The projected tensor also has the remainder required in Proposition \ref{prop:log-energy}.
Indeed, after cutting $V_0$ off away from the pole, the equation for $W:=V[k]-V_0$ has source $O_{C^N}(r^{1-d})$,
because the round vector operator differs from its principal flat operator by $O(r^2)\pp^2+O(r)\pp+O(1)$.
Its generalized inverse has kernel bounded by $C\,d_{g_c}(x,y)^{2-n}$.
Thus $W=O(1)$ for $d=2$, $W=O(1+|\log r|)$ for $d=3$, and $W=O(r^{3-d})$ for $d\ge4$, up to a smooth kernel term.
Scaled elliptic estimates give the differentiated bounds, so in every case
\[
\Pi k=h_0+O_{C^{N+1}}(r^{1-d}).
\]

If $h_0$ vanished, conformal covariance and inversion would express the original homogeneous normal-coordinate tensor
as $\mcd_\delta V$, with $V$ homogeneous of degree $d+1$ on the punctured space.
The field $V$ extends across zero: its growth and that of its first derivative exclude a source supported at zero, and elliptic regularity for $\mcd^*\mcd$ makes it smooth.
The field $V$ is therefore a homogeneous polynomial.
The proof of Proposition \ref{prop:repair-injectivity} applies: the normal-coordinate condition forces, for $b=\divg V$,
$r^2\Delta b=-(n-2)d(d-1)b$, incompatible with the nonnegative Fischer eigenvalues.
Hence the original tensor is zero.

The $O(1)$ remainder in the energy expansion \eqref{eq:log-energy} also holds when the original compactified tensor is cut off before applying $\Pi$.
Indeed, put $h:=\Pi k$, $V:=V[k]$, $t_R:=\Pi(\chi_Rh)$, and
\[
a_R:=\chi_R\mcd_{g_c}V-\mcd_{g_c}(\chi_RV) \quad\text{and}\quad 
\Pi(\chi_Rk)=t_R+\Pi a_R.
\]
The tensor $a_R$ is supported on the cutoff annulus $A_R:=\{R^{-1}<r<2R^{-1}\}$, and the bounds $V=O_{C^{N+1}}(r^{1-d})$ imply $\|a_R\|_{H^N}\le C$ by scaling.
The proof of Proposition \ref{prop:log-energy} gives $\|(I-\Pi)(\chi_Rh)\|_{H^N}\le C$; since $\|\chi_Rh\|_{H^N(A_R)}\le C$, it follows that $\|t_R\|_{H^N(A_R)}\le C$.
Because $p_{N,n}(\mca)$ preserves TT tensors,
\[
\langle t_R,p_{N,n}(\mca)\Pi a_R\rangle
=\langle t_R,p_{N,n}(\mca)a_R\rangle=O(1).
\]
For the last estimate, integrate by parts to use at most $N$ derivatives of each tensor; every resulting integral is supported on $A_R$ and is bounded by the two annular $H^N$ norms.
The self-pairing of $\Pi a_R$ is also bounded by continuity of the polynomial form and boundedness of $\Pi$ on $H^N$.
Consequently replacing $t_R$ by $\Pi(\chi_Rk)$ changes the energy by $O(1)$.
The same comparison holds after $R\pp_R$: the differentiated $a_R$ has the same support and $H^N$ bound,
while $R\pp_Rt_R=\Pi((R\pp_R\chi_R)h)$ has uniformly bounded global $H^N$ norm by annular scaling.

Proposition \ref{prop:log-energy} gives a negative coefficient of $\log R$ in $D^2\mscF(g_c)[t_R,t_R]$.
Dilation of a homogeneous degree-$d$ term multiplies this quadratic coefficient by $2d>0$, so the negative derivative represented by $H^+:=-S^{\mhom}$ has a positive coefficient.
Finite dimensionality and the injectivity just proved make the coefficient of $\log R$ positive definite on $\mfJ_d$.
\end{proof}

\begin{proposition}
\label{prop:repair-high-tail}
Let $0<\mu<r<r_0<1$, set $R:=r/\mu$, and let $H_d\in\mfJ_d$ range over a fixed finite set of degrees with uniformly bounded coefficients.
Choose a smooth radial function $\chi$ equal to one on $B_1$ and zero outside $B_2$, and write $|H_d|:=\|H_d\|_{L^2(\Ss^{n-1})}$.
Put
\[
h_d:=\Omega^2\chi\left(\frac{x}{R}\right)\mu^dH_d(x) \quad\text{and}\quad 
d_0:=\max\{2,\lfloor \frac{n-2N}{2}\rfloor+1\},
\]
and
\[
J_\mu(1)^2:=\sum_{2\le d<\frac{n-2N}{2}}\mu^{2d}|H_d|^2
+\mathbf{1}_{\{n\ {\rm even}\}}\mu^{n-2N}\log\left(\frac{1}{\mu}\right)|H_{\frac{n-2N}{2}}|^2.
\]
The term with degree $(n-2N)/2$ is omitted when $(n-2N)/2<2$.
The constants in this proposition may depend on the fixed finite range of Taylor degrees.
We write $h_{\mathrm{low}}:=\sum_{d=2}^{d_0-1}h_d$ and $h_{\mathrm{high}}:=\sum_{d\ge d_0}h_d$, with both sums restricted to the fixed set of degrees.
For $\mcb_N=D^2\mscF(g_c)$, the quadratic pairings after applying $\Pi$ that contain a Taylor term of degree greater than $(n-2N)/2$ satisfy, for every $\theta>0$,
\begin{equation}\label{eq:repair-small-r-tail}
\begin{split}
&|\mcb_N\left(\Pi h_{\mathrm{low}},\Pi h_{\mathrm{high}}\right)|
+|\mcb_N\left(\Pi h_{\mathrm{high}},\Pi h_{\mathrm{high}}\right)|
\\
&\qquad\le \theta J_\mu(1)^2+C_\theta\mu^{n-2N} r^{2d_0+2N-n}.
\end{split}
\end{equation}
Estimate \eqref{eq:repair-small-r-tail} also holds after one logarithmic scale derivative $\mu\pp_\mu$, with $r$ and the coefficients of $H_d$ fixed.
\end{proposition}
\begin{proof}
Integrating the derivative bounds for $h_d$ over its support, as in the proof of Proposition \ref{thm:strict-york-transmission-new}, gives
\[
\|h_d\|_{H^N}^2\le C\mu^{2d}|H_d|^2
\begin{cases}
1 &\text{for }2d<n-2N,\\
1+\log R &\text{for }2d=n-2N,\\
1+R^{2d+2N-n} &\text{for }2d>n-2N.
\end{cases}
\]
The lower-derivative terms satisfy the same upper bound, as do the product-rule terms with derivatives on $\chi(x/R)$.
For $2e>n-2N$ and $\mu<r$,
\[
\mu^{2e}(1+R^{2e+2N-n})\le C\mu^{n-2N} r^{2e+2N-n}.
\]
Since $1+\log R\le C_{r_0}\log(1/\mu)$, the sum of the terms with $d\le(n-2N)/2$ has $H^N$ norm at most $CJ_\mu(1)$, while the sum with $d>(n-2N)/2$ has norm at most $C\mu^{(n-2N)/2}r^{d_0-(n-2N)/2}$.
Lemma \ref{lem:sobolev-york-new}, estimate \eqref{eq:round-hessian-continuity}, and Young's inequality prove \eqref{eq:repair-small-r-tail}.
For the scale derivative, $\mu\pp_\mu$ multiplies the polynomial term by $d$ and replaces the cutoff by a fixed annular profile.
Its $H^N$ bounds are identical.
The product rule proves the last assertion.
\end{proof}

\subsection{The Pohozaev identity and the mass term}\label{sec:green-current}
We use $P_0=(-\Delta)^N$ and the dilation operator $\msZ$ from \eqref{eq:main-notation2}, with outward Euclidean unit normal $\nu$ on coordinate spheres.
All boundary integrals in this subsection use the Euclidean metric $\delta$.

For functions $f,g$ on a neighborhood of $\pp B_r$, the boundary expression in Green's formula is
\begin{equation}\label{eq:completion-Green-concomitant}
\mfB_N(f,g)
:=\sum_{j=0}^{N-1}
\left\{((-\Delta)^jg)\pp_\nu (-\Delta)^{N-1-j}f
-((-\Delta)^{N-1-j}f)\pp_\nu (-\Delta)^jg\right\}.
\end{equation}
We define the quadratic Pohozaev integral and its nonlinear counterpart by
\begin{align}
\mcp_N^0(r,u)
&:=\frac{1}{2}\int_{\pp B_r}
\left\{r u(-\Delta)^Nu-\mfB_N(u,\msZ u)\right\}\,d\sigma,
\label{eq:completion-flat-current}\\
\mcp_N^\kappa(r,u)
&:=\mcp_N^0(r,u)-\frac{(n-2N)\kappa r}{2n}
\int_{\pp B_r}u^{\frac{2n}{n-2N}}\,d\sigma.\label{eq:nonlinear-Pohozaev-integral}
\end{align}
\begin{lemma}
\label{lem:completion-flat-Pohozaev}
If $u$ is smooth on $\overline {B_r\setminus B_{r_1}}$, where $0<r_1<r$, and $u>0$, then
\begin{equation}\label{eq:completion-Pohozaev-identity}
\mcp_N^\kappa(r,u)-\mcp_N^\kappa(r_1,u)
=\int_{B_r\setminus B_{r_1}}(\msZ u)\left((-\Delta)^Nu-\kappa u^{\frac{n+2N}{n-2N}}\right)\,dx.
\end{equation}
If $u$ is smooth through the origin, the same formula holds with $r_1=0$ and $\mcp_N^\kappa(0,u)=0$.
\end{lemma}

\begin{proof}
For a smooth bounded domain $\mco\subset\R^n$, repeated Green's formula, $[(-\Delta)^N,\msZ]=2N(-\Delta)^N$, and the divergence theorem give
\[
2\int_{\mco}(\msZ u)(-\Delta)^Nu
=\int_{\pp\mco}\bigl((x\cdot\nu)u(-\Delta)^Nu-\mfB_N(u,\msZ u)\bigr).
\]
The critical nonlinear term satisfies
\[
\int_{\mco}(\msZ u)u^{\frac{n+2N}{n-2N}}
=\frac{n-2N}{2n}\int_{\pp\mco}(x\cdot\nu)u^{\frac{2n}{n-2N}}.
\]
Apply these identities to the annulus, with opposite normals on its two boundary components.
\end{proof}

Let $\omega_{n-1}:=|\Ss^{n-1}|$ and put
\begin{equation}\label{eq:completion-fundamental-constant}
\ga_{n,N}
=\frac{\Gamma\left(\frac{n}{2}-N\right)}{4^N\pi^{\frac{n}{2}}\Gamma(N)}.
\end{equation}
Then $(-\Delta)^N(\ga_{n,N}|x|^{2N-n})=\delta_0$ in the weak sense.

\begin{lemma}
\label{lem:completion-mass-current}
Suppose $A>0$ and, near the puncture,
\[
v(x)=A|x|^{2N-n}+B+\eta(x),\qquad
|\pp^j\eta(x)|\le C|x|^{1-j}\quad(0\le j\le2N).
\]
Then
\begin{equation}\label{eq:completion-mass-limit}
\lim_{r\downarrow0}\mcp_N^0(r,v)
=\frac{n-2N}{2\ga_{n,N}}AB.
\end{equation}
\end{lemma}

\begin{proof}
We write $\Phi:=r^{2N-n}$.
Direct differentiation gives
\[
(-\Delta)^j\Phi=C_jr^{2N-n-2j} ,\quad\text{where } 
C_j:=\prod_{t=0}^{j-1}(n-2N+2t)(2N-2-2t) \quad\text{for }  0\le j\le N-1,
\]
and $(-\Delta)^N\Phi=0$ away from the origin.
Thus
\[
C_{N-1}=4^{N-1}(N-1)!\left(\frac{n-2N}{2}\right)_{N-1} \quad\text{and}\quad 
\ga_{n,N}^{-1}=(n-2)\omega_{n-1}C_{N-1}.
\]
The last identity follows either from the gamma recurrence or by integrating $-\Delta(C_{N-1}r^{2-n})$ across a small sphere.

For $v_0=A\Phi+B$, one has $\msZ v_0=-(n-2N)A\Phi/2+(n-2N)B/2$.
The pure singular contribution to the sum $\mfB_N(v_0,\msZ v_0)$ cancels on pairing its $j$-th summand with its $(N-1-j)$-th summand.
The cross terms between $A\Phi$ and $B$ occur only at the two endpoints and give
\[
\int_{\pp B_r}\mfB_N(v_0,\msZ v_0)\,d\sigma
=-(n-2N)AB(n-2)\omega_{n-1}C_{N-1}.
\]
This remains the same formula for $N=1$, when both endpoint contributions occur in a single summand.
Since $(-\Delta)^Nv_0=0$ on the punctured ball, \eqref{eq:completion-flat-current} now gives the stated constant exactly for $v_0$.

Every cross term involving $A\Phi$ and $\eta$ in \eqref{eq:completion-Green-concomitant} is $O(r)$ after surface integration: its total derivative order is $2N-1$, and $r^{n-1}r^{2N-n}r^{1-(2N-1)}=r$.
Terms involving $B$ and $\eta$, or two remainders, tend to zero more rapidly.
Finally, $r\int_{\pp B_r} v(-\Delta)^N\eta\,d\sigma=O(r)+o(1)$ by the assumed derivative bound.
These estimates prove \eqref{eq:completion-mass-limit}.
\end{proof}

\subsection{The Green function expansion and mass identities} \label{sec:repair-Green}

\begin{theorem}
\label{thm:repair-Green}
Let $n\ge3$, $n>2N$, and suppose $P_{2N,g}$ is invertible on a smooth closed manifold.
Suppose that $g$ is the chosen conformal-normal representative at $p$ and that, in its normal coordinates,
\[
g=\exp H ,\quad\text{where } 
H=\sum_{e=d_0}^{n-2N}H_e+O_{C^K}(r^{n-2N+1}) ,\quad 
H_e\in\mfJ_e ,\quad\text{and}\quad  d_0>\frac{n-2N}{2},
\]
where $d_0\ge2$ is an integer and $K\ge4N+2$.
A sum with $d_0>n-2N$ is empty.
The trace-free and normal-coordinate conditions are required only through the displayed Taylor order; in normal coordinates, the exact identity $\det g=1$ near $p$ is a sufficient special case.
Then for every $0<\vep<\min\{1,2d_0+2N-n\}$,
\begin{equation}\label{eq:repair-Green-expansion}
G_g(p,x)=\ga_{n,N}r^{2N-n}
+\sum_{e=d_0}^{n-2N}r^{e+2N-n}a_e(\theta)+m_g(p)+O_{C^{2N}}(r^\vep).
\end{equation}
Each angular coefficient has zero spherical mean: $\int_{\Ss^{n-1}}a_e=0$.
The differentiated remainder, denoted by $\eta$, satisfies $|\pp^j\eta|\le C r^{\vep-j}$ for $0\le j\le2N$ as $r\to0$.
If $2N+1\le n\le2N+3$, every angular coefficient vanishes.
The constant $m_g(p)$ is the Green function mass in the chosen conformal normal coordinates,
namely the finite part remaining after the leading singular term and all displayed angular terms are subtracted.
Its identification with the geometric masses is proved in Proposition \ref{prop:low-order-mass}.
\end{theorem}
\begin{proof}
We write $\Phi:=\ga_{n,N}r^{2N-n}$.
Naturality, constant metric scaling and a Taylor expansion on the unit annulus imply that
\[
DP_{2N,\delta}[H_e]\Phi=r^{e-n}F_e(\theta).
\]
Every nonlinear term in $H$ acting on $\Phi$ contains metric factors whose Taylor degrees sum to at least $2d_0>n-2N$.
Likewise, a first metric variation of degree at least $d_0$, acting on a correction of degree $e+2N-n$, contributes $O(r^{d_0+e-n})$, with $d_0+e>n-2N$.
Thus each angular correction in \eqref{eq:repair-Green-expansion} is determined independently by its linear source $DP_{2N,\delta}[H_e]\Phi$.

The map $H_e\mapsto\int_{\Ss^{n-1}} F_e$ is an $O(n)$-invariant linear functional on $\mfJ_e$.
That space has no invariant vector: an invariant homogeneous tensor has the form $a r^e\delta_{ij}+b r^{e-2}x_ix_j$, and its trace and radial contraction force $a=b=0$.
Averaging over $O(n)$ therefore proves $\int_{\Ss^{n-1}} F_e=0$.
The source has a finite spherical harmonic expansion.
For $n\ge5$, Proposition \ref{thm:finite-SVT-exhaustion-new} shows that scalar summands occur only at degrees $e=\ell+2+2j\ge4$, with $\ell\ge2$ and $j\ge0$.
Equivariance of the scalar source therefore implies $F_2=F_3=0$ whenever these degrees are present.

For a degree-$j$ spherical harmonic, direct polar differentiation gives
\begin{equation}\label{eq:repair-indicial}
\begin{gathered}
(-\Delta)^N(r^\lambda Y_j)=I_{N,n}(\lambda,j)r^{\lambda-2N}Y_j,\\
I_{N,n}(\lambda,j):=\prod_{a=0}^{N-1}
\left[j(j+n-2)-(\lambda-2a)(\lambda-2a+n-2)\right].
\end{gathered}
\end{equation}
The indicial roots are $\lambda=j+2a\ge0$ and $\lambda=-j-n+2+2a\le2N-n$.
Thus $I_{N,n}(\lambda,j)$ is nonzero for $2N-n<\lambda<0$.
At $\lambda=0$ the only zero mode is $j=0$, which is absent from the source.
Taking $\lambda=e+2N-n$ therefore defines a unique mean-zero $a_e$ solving $(-\Delta)^N(r^{e+2N-n}a_e)=-DP_{2N,\delta}[H_e]\Phi$.

For $n=3,4$, necessarily $N=1$.
If $n=3$, then $n-2N=1<d_0$, so the angular sum is empty.
If $n=4$, the only possible retained term is $H_2$.
The degree-two normal-coordinate identity $\operatorname{Rm}'_\delta[H_2]_{ikj\ell}x^kx^\ell=-3(H_2)_{ij}$ and $\tr H_2=0$ imply $\Ric'_\delta[H_2]=0$,
so the linearized curvature determined by $H_2$ consists only of its Weyl component.
Moreover, $\det \exp(tH_2)=1$ and $\exp(tH_2)x=x$, so $\Delta_{\exp(tH_2)}f(r)=\Delta f(r)$ for every radial $f$.
Since $P_{2,g}=-\Delta_g+\operatorname{Scal}_g/6$ in dimension four, $DP_{2,\delta}[H_2]\Phi=0$ and $a_2=0$.

Let $v_0:=\Phi+\sum r^{e+2N-n}a_e$, with a smooth cutoff away from $p$.
Its error satisfies the weak equation
\[
P_{2N,g}v_0-\delta_p=:f,\qquad
|\pp^j f|\le C r^{\tau-2N-j}\quad(0\le j\le2N),\qquad
\tau:=\min\{1,2d_0+2N-n\}>0.
\]
This follows on each unit rescaled annulus: the metric Taylor remainder starts at degree $n-2N+1$, and every product of two retained metric terms has degree at least $2d_0$.
No additional distribution supported at $p$ occurs: for a correction $r^{e+2N-n}a_e$ with $e>0$, the largest boundary term on $\pp B_r$ in an order-$2N$ Green formula is $O(r^e)$;
all variable-coefficient corrections to the leading fundamental term also gain a positive power of $r$.

Choose $p_0>n/(2N)$ sufficiently close to $n/(2N-\tau)$ that $2N-n/p_0>\vep$, while $p_0<n/(2N-\tau)$.
Then $f\in L^{p_0}$.
Invertibility and the global elliptic estimate give $w=-P_{2N,g}^{-1}f\in W^{2N,p_0}\subset C^{0,\vep}$.
Hence $G_g(p,\cdot)=v_0+w$ and $w(p)=m_g(p)$.
Finally, we apply the scaled interior Schauder estimate to $w-w(p)$ on annuli.
Its supremum is $O(r^\vep)$; its rescaled source is $O(r^\tau)$, and the extra source $w(p)P_{2N,g}1$ is bounded before rescaling.
This gives $|\pp^j(w-w(p))|\le C r^{\vep-j}$ through order $2N$.
\end{proof}

\begin{proposition}
\label{prop:repair-Green-example}
For $N=1$ and $n=8$, there are metrics in normal coordinates with determinant one and $H_2=H_3=0$ for which the parametrix for the Green function contains a nonzero $r^{-2}$ angular term.
\end{proposition}
\begin{proof}
Take $Y:=x_1^2-x_2^2$ and
\[
H_{ij}:=r^4\pp_i\pp_jY
-r^2(x_i\pp_jY+x_j\pp_iY)
+\frac{2}{7}r^2Y\delta_{ij}+\frac{12}{7}Yx_ix_j.
\]
It belongs to $\mfJ_4$.
Set $g_\ep:=\exp(\ep H)$ locally, with $\ep\ne0$.
Then $\det g_\ep=1$ and $(g_\ep)_{ij}x^j=x_i$ exactly.
A direct Cartesian calculation gives
\[
\operatorname{Scal}'_\delta[H]=\pp_i\pp_jH_{ij}=\frac{864}{7}Y.
\]
The Laplace--Beltrami operator on any radial function is exactly the Euclidean one in this determinant-one normal metric.
Therefore, with $P_{2,g}=-\Delta_g+3\operatorname{Scal}_g/14$,
\[
DP_{2,\delta}[H](\ga_{8,1}r^{-6})
=\frac{3}{14}\frac{864}{7}\ga_{8,1}r^{-4}Y_2(\theta) ,\quad\text{where }  Y=r^2Y_2.
\]
Since $-\Delta(r^{-2}Y_2)=24r^{-4}Y_2$, the degree-four correction is
\[
-\frac{54}{49}\ep\ga_{8,1}r^{-2}Y_2(\theta).
\]
Products of two metric perturbations first have Taylor degree eight, so their contributions cannot cancel the correction arising from the degree-four term $H$.
\end{proof}

\begin{lemma}
\label{lem:repair-angular-mass}
Suppose
\[
v=A r^{2N-n}+B+\psi+\eta ,\quad\text{where }  A>0 \quad\text{and}\quad 
\int_{\Ss^{n-1}}\psi(r,\theta)\,d\theta=0,
\]
where $|\pp^j\psi|\le C r^{d_0+2N-n-j}$, $d_0>(n-2N)/2$, and $|\pp^j\eta|\le C r^{\vep-j}$, $\vep>0$, through order $2N$.
For the Pohozaev integral \eqref{eq:completion-flat-current} with the outward normal,
\begin{equation}\label{eq:repair-angular-mass}
\lim_{r\downarrow0}\mcp_N^0(r,v)
=\frac{n-2N}{2\ga_{n,N}}AB.
\end{equation}
\end{lemma}
\begin{proof}
Spherical averaging commutes with $-\Delta$, $\msZ$, and radial differentiation.
Consequently every term in the polarized Pohozaev integral pairing a radial function with $\psi$ integrates to zero.
The $\psi$--$\psi$ terms are $O(r^{2d_0+2N-n})$: the boundary expression \eqref{eq:completion-Green-concomitant} has total derivative order $2N-1$, and the $r v(-\Delta)^Nv$ term has the same scaling.
The $\psi$--$\eta$ terms are $O(r^{d_0+\vep})$; the singular--$\eta$ terms are $O(r^\vep)$.
Applying Lemma \ref{lem:completion-mass-current} to $Ar^{2N-n}+B$ gives \eqref{eq:repair-angular-mass}.
\end{proof}

For a finite sum $v:=\sum b_jG_g(p_j,\cdot)$ with distinct poles and $p_i\in\mcz$, take its conformal-normal representative $\Lambda_{p_i}^{-1}v$ at $p_i$.
The transformation law $G_{g_{p_i}}(p_j,x)=\Lambda_{p_i}(p_j)^{-1}\Lambda_{p_i}(x)^{-1}G_g(p_j,x)$ and $\Lambda_{p_i}(p_i)=1$ identify its finite part as
\begin{equation}\label{eq:repair-green-sum-finite-part}
B=b_i m_g(p_i)
+\sum_{j\ne i}b_j\Lambda_{p_i}(p_j)G_{g_{p_i}}(p_j,p_i)
=b_i m_g(p_i)+\sum_{j\ne i}b_jG_g(p_j,p_i).
\end{equation}

\begin{proposition}
\label{prop:low-order-mass}
Let $N\ge1$, $n>2N$, and suppose that $(M^n,g)$ is closed, connected and smooth and $\ker P_{2N,g}=0$.
Fix the conformal-normal family of Section \ref{sec:main-mass}, and use $d_0$ from \eqref{eq:main-mass-expansion}.
At $\xi\in\mcz$, write $G:=G_{g_\xi}(\xi,\cdot)$.
Let $F$ be any smooth positive function on $M\setminus\{\xi\}$ agreeing with $G$ near $\xi$, and define
\[
\hg_\xi:=
\begin{cases}
\left(\frac{F}{\ga_{n,1}}\right)^{\frac{4}{n-2}}g_\xi &\text{for }N=1,\\
F^{\frac{4}{n-2N}}g_\xi &\text{for }N\ge2.
\end{cases}
\]
The manifold $(M\setminus\{\xi\},\hg_\xi)$ is complete, has one asymptotically flat end, and satisfies $Q_{2N,\hg_\xi}=0$ near infinity.
In an asymptotic chart, with $\rho$ the Euclidean radius, its decay is
\begin{equation}\label{eq:low-order-AF}
\hg_{\xi,ij}-\delta_{ij}
=O_{C^{2N}}(\rho^{-\tau}) ,\quad\text{where } 
\tau:=\min\{d_0,n-2N\}>\frac{n-2N}{2}.
\end{equation}
For $N\ge2$, we define
\begin{equation}\label{eq:all-order-mass-integral}
\mathscr E_{2N}(\hg_\xi)
:=\int_{M\setminus\{\xi\}}
(-\Delta_{\hg_\xi})^{N-1}\operatorname{Scal}_{\hg_\xi}\,dv_{\hg_\xi}.
\end{equation}
This integral is absolutely convergent.
With $\hg_\xi=(F/\ga_{n,1})^{4/(n-2)}g_\xi$ for $N=1$, the mass identities are
\begin{equation}\label{eq:low-order-masses}
m_{\mADM}(\hg_\xi)=\frac{2}{\ga_{n,1}}m_g(\xi)\quad(N=1),\qquad
\mathscr E_{2N}(\hg_\xi)=\frac{4(n-1)}{n-2N}m_g(\xi)\quad(N\ge2).
\end{equation}
The higher-order identity extends \cite[Theorem 5.1]{Michel} to every $n>2N$ at the poles in $\mcz$.
The corresponding mass integrals are
\begin{equation}\label{eq:low-order-mass-definitions}
\begin{aligned}
m_{\mADM}(\hg)
&:=\frac{1}{2(n-1)\omega_{n-1}}
\lim_{\rho\to\infty}\int_{|y|=\rho}
(\pp_j\hg_{ij}-\pp_i\hg_{jj})\nu_i\,d\sigma,\\
\mathscr E_{2N}(\hg)
&=\lim_{\rho\to\infty}\int_{|y|=\rho}
\pp_k(-\Delta)^{N-2}
\bigl(\pp_j\pp_j\hg_{ii}
-\pp_i\pp_j\hg_{ij}\bigr)\nu_k\,d\sigma \quad\text{for }  N\ge2.
\end{aligned}
\end{equation}
Here the Laplacian, normal and surface measure are Euclidean.
If $G_g(x,y)>0$ for $x\ne y$, one may take $F=G$, and then $Q_{2N,\hg_\xi}=0$ everywhere.
\end{proposition}
\begin{proof}
Fix $\xi\in\mcz$, suppress the pole from the notation, and write $\ga:=\ga_{n,N}$ and $m:=m_g(\xi)$.
Theorem \ref{thm:repair-Green} gives $G\sim\ga r^{2N-n}>0$ near $\xi$.
Choose a cutoff $\zeta$ supported in a neighborhood where $G>0$, equal to one on a smaller neighborhood.
Then $F:=\zeta G+(1-\zeta)c$, for any constant $c>0$, is an admissible positive extension.
Conformal covariance gives $Q_{2N,\hg}=0$ near infinity, and everywhere if $F=G>0$ globally.
Since $F\sim\ga r^{2N-n}$, the metric is comparable to $r^{-4}g_\xi$ near the puncture, so it is complete and has one end.

To evaluate the mass integrals, we use the metric $\bg:=(F/\ga)^{4/(n-2N)}g_\xi$ and the inverted coordinates $y=x/r^2$, $\rho:=|y|$.
Thus $\hg=\bg$ when $N=1$, whereas $\hg=\ga^{4/(n-2N)}\bg$ when $N\ge2$.
Choose $0<\vep<\min\{1,2d_0+2N-n\}$.
Set $b:=m/\ga$.
Theorem \ref{thm:repair-Green} and the Taylor expansion of $g_\xi=\exp H$ yield
\begin{equation}\label{eq:low-order-inverted-expansion}
\bg_{ij}
=\delta_{ij}+K_{ij}
+\frac{4}{n-2N}\bigl(b\rho^{2N-n}+\psi\bigr)\delta_{ij}
+O_{C^{2N}}(\rho^{2N-n-\vep}).
\end{equation}
Here $K$ is the inversion of the finite sum of the homogeneous terms $H_d$ in $\log g_\xi$, $d_0\le d\le n-2N$, and
$\psi:=\ga^{-1}\sum_{d=d_0}^{n-2N}\rho^{-d}a_{\xi,d}(\theta)$.
The normal-coordinate and trace-free identities for these jets, and the mean-zero identities in \eqref{eq:repair-Green-expansion}, imply
\[
Ky=0 ,\quad  \tr K=0 ,\quad 
\int_{\Ss^{n-1}}\psi(\rho\theta)\,d\theta=0 ,\quad\text{and}\quad  K,\psi=O_{C^{2N}}(\rho^{-d_0}).
\]
Indeed, the differential of inversion is a scalar multiple of $I-2\theta\otimes\theta$, which fixes each retained tangential tensor under conjugation.
Finite-order determinant normalization makes each retained $H_d$ trace-free.
The quadratic terms of $\exp H$ have degree at least $2d_0>n-2N$, while the metric Taylor remainder and all omitted products involving the relative correction to the Green function have degree at least $n-2N+1$.
These bounds give the differentiated remainder in \eqref{eq:low-order-inverted-expansion}.
They prove \eqref{eq:low-order-AF} for $\bg$; the further coordinate dilation $y\mapsto\ga^{2/(n-2N)}y$ gives that estimate for $\hg$ when $N\ge2$.

\medskip
We next prove absolute convergence of the curvature integral \eqref{eq:all-order-mass-integral} for $N\ge2$.
In the Euclidean end, we write
\begin{equation}\label{eq:scalar-curvature-linearization}
\mcs(e):=\operatorname{Scal}'_\delta[e]
=\pp_i\pp_j e_{ij}-\Delta\tr e.
\end{equation}
The flat linearization of the $Q_{2N}$-curvature satisfies
\begin{equation}\label{eq:all-order-Q-linearization}
Q'_{2N,\delta}[e]=\frac{1}{2(n-1)}(-\Delta)^{N-1}\mcs(e).
\end{equation}
To see this, naturality and constant metric scaling make this linearization an $O(n)$-invariant constant-coefficient scalar operator of differential order $2N$.
Its Fourier symbol is a linear combination of $|\zeta|^{2N}\tr e$ and $|\zeta|^{2N-2}\zeta_i\zeta_j e_{ij}$.
Vanishing on infinitesimal diffeomorphisms forces it to be a multiple of $(-\Delta)^{N-1}\mcs$.
For the conformal variation $e=2\vph\delta$, covariance gives $Q'_{2N,\delta}[2\vph\delta]=(-\Delta)^N\vph$, while
$\mcs(2\vph\delta)=-2(n-1)\Delta\vph$.
This determines the coefficient and sign in \eqref{eq:all-order-Q-linearization}.

For $\bg=\delta+e$, Taylor expansion on unit rescaled annuli now gives
\begin{equation}\label{eq:nonlinear-curvature-remainder}
Q_{2N,\bg}-\frac{1}{2(n-1)}(-\Delta_{\bg})^{N-1}\operatorname{Scal}_{\bg}
=O(\rho^{-2\tau-2N}).
\end{equation}
Indeed, the difference has scaling weight $-2N$ and vanishes to second order at the flat metric.
Since $Q_{2N,\bg}=0$ near infinity and $2\tau>n-2N$, the absolute integral on the end is bounded by a constant times
$\int_{1}^\infty\rho^{n-2N-1-2\tau}\,d\rho<\infty$.
The integrand is smooth on the remaining compact region.
Constant metric scaling gives absolute integrability for $\hg$ as well.

The divergence theorem on the compact region bounded by $\{|y|=\rho\}$ gives
\begin{equation}\label{eq:all-order-linear-flux}
\mathscr E_{2N}(\bg)
=-\lim_{\rho\to\infty}\int_{|y|=\rho}
\pp_{\nu_{\bg}}\bigl((-\Delta_{\bg})^{N-2}\operatorname{Scal}_{\bg}\bigr)\,d\sigma_{\bg}
=-\lim_{\rho\to\infty}\int_{|y|=\rho}
\pp_\rho(-\Delta)^{N-2}\mcs(e)\,d\sigma.
\end{equation}
The second equality follows by linearizing the curvature, iterated Laplacian, normal and surface measure.
The pointwise error is $O(\rho^{-2\tau-2N+1})$, hence its integrated contribution is $O(\rho^{n-2N-2\tau})=o(1)$.
Since $-\mcs(e)=\pp_j\pp_j\bg_{ii}-\pp_i\pp_j\bg_{ij}$, the last expression is the second boundary integral in \eqref{eq:low-order-mass-definitions} for every $N\ge2$.

\medskip
The angular terms in \eqref{eq:low-order-inverted-expansion} cancel in this linear boundary integral and in the ADM mass integral.
For $K$, differentiating $Ky=0$ gives
\[
y_i\pp_jK_{ij}=-\tr K=0.
\]
Thus $K$ has zero contribution to the ADM mass integral, and $\divg K$ is tangent to each coordinate sphere.
The scalar $\mcs(K)=\divg_\delta\divg_\delta K$ in \eqref{eq:scalar-curvature-linearization} therefore has zero spherical mean.
Spherical averaging commutes with the Euclidean Laplacian and radial differentiation, so $K$ also contributes zero to \eqref{eq:all-order-linear-flux}.
The mean-zero scalar term $\psi\delta$ contributes zero to each mass integral for the same reason, using $\mcs(\psi\delta)=-(n-1)\Delta\psi$.
The remainder contributes $O(\rho^{-\vep})$, since the integrand has differential order $2N-1$ and is integrated over a sphere of dimension $n-1$.

Only the radial term $(4b/(n-2N))\rho^{2N-n}\delta$ remains.
For $N=1$ its ADM mass, defined in \eqref{eq:low-order-mass-definitions}, is $2b$, giving the first identity in \eqref{eq:low-order-masses}.
For $N\ge2$, successive applications of $-\Delta$ give
\[
-\int_{|y|=\rho}\pp_\rho(-\Delta)^{N-1}(\rho^{2N-n})\,d\sigma
=2^{N-1}(N-1)!\,\omega_{n-1}\prod_{j=1}^{N}(n-2j)
=\ga^{-1}.
\]
Therefore $\mathscr E_{2N}(\bg)=4(n-1)m/((n-2N)\ga^2)$.
For any constant $\lambda>0$, the mass integral \eqref{eq:all-order-mass-integral} satisfies
$\mathscr E_{2N}(\lambda^2\bg)=\lambda^{n-2N}\mathscr E_{2N}(\bg)$.
Since $\hg=\ga^{4/(n-2N)}\bg$, this proves the second identity in \eqref{eq:low-order-masses}.
At $N=2$ this also proves the fourth-order identity of \cite[Proposition 8.3]{GKW} under the finite-order determinant normalization used here.
\end{proof}

\begin{corollary}
\label{cor:low-order-compactness}
Let $N\in\{1,2\}$, $n>2N$, and let $(M^n,g)$ be a closed connected smooth manifold not conformally diffeomorphic to the unit round sphere $(\Ss^n,g_c)$.
Assume $\ker P_{2N,g}=0$ and $G_g(x,y)>0$ for $x\ne y$; when $N=2$, assume in addition that $Y_g>0$.
Then
\[
m_g(\xi)>0 \quad \text{for every }\xi\in\mcz,
\]
so the mass hypothesis \eqref{eq:mass-hypothesis} holds.
Consequently, the conclusion of Theorem \ref{thm:ac-main} holds for $3\le n\le24$ when $N=1$, and for $5\le n\le24$ when $N=2$.
\end{corollary}
\begin{proof}
Fix $\xi\in\mcz$ and take $F=G_{g_\xi}(\xi,\cdot)$ in Proposition \ref{prop:low-order-mass}.
Write $\hg:=\hg_\xi$ for the resulting complete asymptotically flat Green metric and $G:=G_{g_\xi}(\xi,\cdot)$.
For $N=1$, the Green metric is scalar flat and \eqref{eq:low-order-AF} has rate greater than $(n-2)/2$.
The smooth case of \cite[Theorem 1.1]{KWW} therefore gives nonnegative ADM mass and Euclidean rigidity; the decay satisfies \cite[Definition 2.2]{KWW}.
For $N=2$, the Yamabe quotient of a compactly supported function $v$ on the Green metric equals that of $G^{(n-2)/(n-4)}v$ on $g_\xi$, and is consequently at least $Y_g>0$.
The Green metric also has $Q_{4,\hg}=0$ and the $C^4$ decay required by \cite[Theorem A]{ALL}.
That theorem gives $\mathscr E_4(\hg)\ge0$, with equality only for Euclidean space.
In either order, Euclidean rigidity implies that $(M,g)$ is conformally diffeomorphic to the round sphere:
the conformal map on the punctured manifold extends across the pole by Liouville's theorem in a local conformally flat chart.
Such a chart exists because local conformal flatness holds away from the pole and extends to the pole by smoothness.
The exclusion of the round conformal class therefore makes the corresponding geometric mass strictly positive.
The identities \eqref{eq:low-order-masses} give $m_g(\xi)>0$ at every $\xi\in\mcz$, proving \eqref{eq:mass-hypothesis}.
Theorem \ref{thm:ac-main} then gives the stated compactness conclusions.
\end{proof}

\section{Isolated blow-up points and local estimates}\label{sec:actual-comparison}
We derive a local Pohozaev estimate, prove that blow-up points are isolated simple, and obtain quantitative vanishing of the metric coefficients.
The estimate retains the annular norm $A_i(\rho)$ of the approximation error, defined in \eqref{eq:ac-display-2-5}, for the refinement in Section \ref{sec:green-matching}.

\subsection{Uniform local tools and operator estimates}\label{sec:ac-tools}
The local results used below are point selection \cite[Proposition 4.1]{MP}, convergence to a bubble \cite[Lemma 4.2]{MP}, and annular Harnack and derivative estimates \cite[Lemma 4.3]{MP}.
These statements apply in every dimension $n>2N$ to positive global solutions when $P_{2N,g}$ is invertible and its Green function is positive.
We use them for the critical equation with coefficient $\kappa_0$, as normalized in Section \ref{sec:main-mass}.
At a selected center $\xi_i$, we write $g_i$ for the coordinate expression of $g_{\xi_i}$ in its conformal-normal chart.
In all local expressions, $u_i$ denotes the original solution divided by $\Lambda_{\xi_i}$ and expressed in that chart.
In these coordinates, we define the rescaled solution and the centered bubble by
\begin{equation}\label{eq:ac-local-bubble}
v_i(y):=\mu_i^{\frac{n-2N}{2}}u_i(\mu_i y) \quad\text{and}\quad 
B_i(x):=\mu_i^{-\frac{n-2N}{2}}U\left(\frac{x}{\mu_i}\right).
\end{equation}
In statements on $M$, $u_i$ denotes the original global solution and $B_i$ denotes the corresponding bubble multiplied by $\Lambda_{\xi_i}$.

We use the notions of \emph{isolated blow-up point} and \emph{isolated simple blow-up point} introduced by Schoen; see \cite[Section 2]{KMS} and \cite[Definition 2.2]{GKW}.
\begin{definition}[Isolated blow-up point]\label{def:ac-isolated}
In the conformal-normal coordinates centered at $\xi_i$, a limit $\xi_i\to\xi$ of local maximum points with $u_i(0)\to\infty$ is an \emph{isolated blow-up point} for $\{u_i\}$ if
\begin{equation}\label{eq:ac-isolation-bound}
|x|^{\frac{n-2N}{2}}u_i(x)\le C \quad (0<|x|<r_0)
\end{equation}
for constants $C,r_0>0$ independent of $i$.
\end{definition}

\begin{definition}[Isolated simple blow-up point]\label{def:ac-isolated-simple}
For an isolated blow-up point as in Definition \ref{def:ac-isolated}, let
\begin{equation}\label{eq:ac-spherical-average}
\bar u_i(r):=\frac{1}{|\pp B_r|_{g_i}}\int_{\pp B_r}u_i\,d\sigma_{g_i}.
\end{equation}
The point is an \emph{isolated simple blow-up point} if, after decreasing $r_0$, the function $r^{(n-2N)/2}\bar u_i(r)$ has exactly one critical point in $(0,r_0)$ for all sufficiently large $i\in\mathbb N$.
\end{definition}

Before \eqref{eq:ac-isolation-bound} has been established with a radius independent of $i$, we use the same growth bound on balls whose radii may tend to zero:
\begin{equation}\label{eq:ac-local-growth}
|x|^{\frac{n-2N}{2}}u_i(x)\le C\quad(0<|x|<8\rho_i) \quad\text{and}\quad 
\frac{\mu_i}{\rho_i}\to0\quad\text{as }i\to\infty.
\end{equation}
Here $8\rho_i$ is smaller than the fixed conformal-normal coordinate radius.
Point selection supplies the growth bound when $\rho_i$ is at most one sixteenth of the distance to the other selected points, up to a fixed change of coordinate radius.
The estimates below are first proved on balls where
\begin{equation}\label{eq:ac-relative-bubble}
\sup_{|x|<4\rho_i}\left|\frac{u_i(x)}{B_i(x)}-1\right|\le\vep \quad\text{and}\quad  \frac{\mu_i}{\rho_i}\to0,
\end{equation}
for a sufficiently small fixed $\vep>0$.
Proposition \ref{prop:ac-saturation} will establish this comparison on fixed balls at isolated blow-up points and will prove that they are simple.
We work on $B_{4R}$ in the rescaled variables and choose cutoffs supported strictly inside this ball, decreasing the physical radius by fixed factors when necessary.
For the local estimates, we temporarily omit the sequence index $i$, writing $\mu=\mu_i$, $\rho=\rho_i$, $v=v_i$, and $H_d=H_{\xi_i,d}$ from \eqref{eq:conformal-normal-metric-expansion}.
In the $y$-variables, we define the rescaled metric polynomial $K$, its truncation $K_{\mathrm{tr}}$, and the scalar correction $\psi$ by
\begin{equation}
\begin{gathered}
R:=\frac{\rho}{\mu} \quad\text{and}\quad  K(y):=\sum_{d=2}^{T-1}\mu^dH_d(y),\\
K_{\mathrm{tr}} :=\sum_{2\le d<n-2N}\mu^dH_d \quad\text{and}\quad  \psi:=\sum_{2\le d<n-2N}\mu^d\Psi[H_d].
\end{gathered}             \label{eq:ac-display-2-1}
\end{equation}
Corollary \ref{cor:repair-high-corrector} yields $L_U\Psi[H_d]=-DP_0[H_d]U$ and $\Psi[H_d]=O_{C^{2N}}((1+|y|)^{d+2N-n})$.
The subscript on $DP_0$ denotes metric variation at the flat metric.
On $|y|\le4R$,
\begin{equation}
\sum_{j\le2N+2}(1+|y|)^j|\pp^jK(y)|
\le C\mu^2(1+|y|)^2\le C\rho^2.                 \label{eq:ac-display-2-2}
\end{equation}
The coefficients of each homogeneous Taylor term $H_d$ are bounded uniformly as the center varies.

For a scalar function, we write
\begin{equation}\label{eq:scalar-energy-norm}
\|v\|_{\mce}^2:=\int_{\R^n} vP_0v ,\quad\text{where }  \mce:=\dot H^N(\R^n),
\end{equation}
let $\mce^*$ denote the continuous dual of this space, and define the following norm for a compactly supported metric perturbation $K$:
\begin{equation}
\|K\|_{\mch}^2
:=\sum_{j=0}^{N}\int_{\R^n} U^2(1+|y|)^{2j-2N}|\pp^jK|^2\,dy.
\label{eq:ac-display-2-3}
\end{equation}
For the local polynomial $K$ in \eqref{eq:ac-display-2-1}, we apply this norm after multiplication by a smooth radial cutoff function supported in $B_{4R}$.
Each derivative acts on the entire product.
The norm $\|v\|_{\mce}$ is equivalent, under $v=U\bar v$, to the conformal energy norm \eqref{eq:conformal-energy-space}.
The norm $\|K\|_{\mch}$ is equivalent to $\|\Omega^2K\|_{H^N(\Ss^n,g_c)}$;
this follows from the conformal change of connection and the Hardy inequalities at the projection pole $p_+$.
The constants in these norm equivalences may depend on $N,n$, but not on $\mu,R$, or the center.

To measure the homogeneous Taylor terms of degree at most $(n-2N)/2$, we define
\begin{equation}
J_\mu(\rho)^2:=
\sum_{2\le d<\frac{n-2N}{2}}\mu^{2d}|H_d|^2
+\mathbf{1}_{\{n\ {\rm even}\}}\mu^{n-2N}\log R\,|H_{\frac{n-2N}{2}}|^2,
\label{eq:ac-display-2-4}
\end{equation}
where the term of degree $(n-2N)/2$ is omitted when $(n-2N)/2<2$.
We assume $R\to\infty$ whenever we use this norm.

For the \emph{approximation error} $w:=v-U-\psi$, with $\psi$ from \eqref{eq:ac-display-2-1}, let
\begin{equation}
A:=\sum_{j=0}^{2N}R^j\|\pp^jw\|_{L^\infty(A_R^*)}
+R^{2N+\vth}[\pp^{2N}w]_{\vth,A_R^*},
\label{eq:ac-display-2-5}
\end{equation}
where $A_R^*:=\{y:R<|y|<3R\}$ and $0<\vth<1$.
The annulus $A_R^*$ contains the support of each derivative of the function $\chi_0$ used to form $\chi_0w$ in Lemma \ref{lem:ac-decomposition}.
Before matching, the differentiated bubble bounds and Corollary \ref{cor:repair-high-corrector} imply $A=O(R^{2N-n})$ as $R\to\infty$.

\begin{lemma}
\label{lem:ac-tame}
Let $K$ be a compactly supported symmetric trace-free matrix field specifying the metric $\exp K$, assume \eqref{eq:ac-display-2-2}, and put $\delta:=C\rho^2$.
For sufficiently small $\delta$, the weak form of $P(K):=P_{2N,\exp K}$ satisfies
\begin{align}
\|DP_0[K]U\|_{\mce^*}&\le C\|K\|_{\mch},\label{eq:ac-linear-source}\\
|\langle(P(K)-P_0)f,h\rangle|&\le C\delta\|f\|_{\mce}\|h\|_{\mce},\label{eq:ac-form}\\
\|(P(K)-P_0-DP_0[K])U\|_{\mce^*}
&\le C\delta\|K\|_{\mch}.                 \label{eq:ac-source}
\end{align}
The analogous cubic Taylor error with two $U$ factors is bounded by $C\delta\|K\|_{\mch}^2$.
These estimates also hold after applying $y\cdot\pp$ once to the metric coefficients, and with smooth radial cutoffs at scale $R$.
\end{lemma}
\begin{proof}
For fixed $N,n$, naturality and homogeneity under constant metric scaling express $P(K)$ as a finite sum of smooth contractions of inverse metrics and metric derivatives,
with total differential order $2N$ counting derivatives of both the metric and the scalar argument.
We apply Taylor's formula to these smooth coefficients in $K$.
Repeatedly integrate by parts until each scalar test factor has at most $N$ derivatives.
For a monomial with two or more metric factors, continue the same redistribution until any factor to be estimated in $L^2$ has at most $N$ derivatives.
Since the total differential order is $2N$, at most one factor can carry more than $N$ derivatives at any step.
Remaining metric factors and their weighted derivatives are bounded using \eqref{eq:ac-display-2-2}.
Derivatives of $U$ obey $|\pp^jU|\le C U(1+|y|)^{-j}$.
The resulting $L^2$ metric expressions are precisely the summands of \eqref{eq:ac-display-2-3}.
The scalar expressions are bounded by
\[
\|(1+|y|)^{-a}\nabla^{N-a}f\|_2\le C\|f\|_{\mce} \quad\text{for }  0\le a\le N,
\]
using the iterated Hardy inequality, valid because $n>2N$.
For \eqref{eq:ac-linear-source}, estimate one metric factor in the weighted $L^2$ norm
and one scalar test factor in the $L^2$ norm controlled by the Hardy inequality, and take the supremum over scalar test functions of unit energy norm.
Retaining a second small metric factor proves \eqref{eq:ac-source}, and the same redistribution proves \eqref{eq:ac-form}.
At cubic order, retaining two metric factors in $L^2$ and one in the weighted supremum norm yields the stated cubic error.
Derivatives of a radial cutoff contribute $R^{-j}$ on its cutoff annulus and satisfy the same weighted bounds.
Finally $y\cdot\pp$ preserves every derivative weight and the bound \eqref{eq:ac-display-2-2}.
\end{proof}

\begin{lemma}[Uniform cutoff estimates]\label{lem:ac-uniform-cutoffs}
Let $K,K_{\mathrm{tr}}$ be as in \eqref{eq:ac-display-2-1}, with a fixed finite range of degrees.
Fix smooth radial cutoffs $\chi=1$ on $B_1$, supported in $B_2$, and $\eta=1$ on $B_2$, supported in $B_4$.
We define $\chi_R(y):=\chi(y/R)$ and $\eta_R(y):=\eta(y/R)$ for $R\ge2$.
Then the metric norm satisfies
\begin{equation}\label{eq:ac-uniform-jet-norm}
\|\eta_RK\|_{\mch}^2
\asymp
\sum_d\mu^{2d}|H_d|^2
\left(1+\int_1^R t^{2d+2N-n-1}\,dt\right).
\end{equation}
The scalar correction term satisfies
\begin{equation}\label{eq:ac-uniform-corrector-cutoff}
\|\chi_R\Psi[K_{\mathrm{tr}}]\|_{\mce}
\le C\|\eta_RK_{\mathrm{tr}}\|_{\mch}.
\end{equation}
The constants are independent of $R,\mu$, and the center.
The same upper bounds hold when either cutoff is differentiated once by $R\pp_R$, or when one coordinate-dilation derivative is applied to the quantity on the left.

\end{lemma}
\begin{proof}
To prove \eqref{eq:ac-uniform-jet-norm}, first rescale a dyadic annulus of radius $t\ge1$ to a fixed unit annulus.
The norm of a polynomial on that annulus controls the sum of squares of its homogeneous coefficients, with a constant depending only on the finite degree range.
Since $U\asymp t^{2N-n}$ there, the contribution of degree $d$ is comparable to $\mu^{2d}|H_d|^2t^{2d+2N-n}$.
Summing these annuli and adding the unit ball proves the two-sided estimate with constants independent of the outer radius.
The radial integral is bounded for $2d<n-2N$, logarithmic at equality, and of order $R^{2d+2N-n}$ for $2d>n-2N$.

For each integer $d$ with $2\le d<n-2N$, Corollary \ref{cor:repair-high-corrector} implies
\begin{equation}\label{eq:ac-corrector-decay}
|\pp^j\Psi[\mu^dH_d](y)|
\le C\mu^d|H_d|(1+|y|)^{d+2N-n-j} \quad\text{for }  0\le j\le2N.
\end{equation}
Its squared order-$N$ derivative norm on the same annulus is bounded by $C\mu^{2d}|H_d|^2t^{2d+2N-n}$.
Terms in which derivatives fall on $\chi_R$ have the same power at $t\asymp R$.
Sum over the finite set of degrees and use the lower bound in \eqref{eq:ac-uniform-jet-norm} to obtain \eqref{eq:ac-uniform-corrector-cutoff}.
A logarithmic cutoff derivative or a coordinate-dilation derivative preserves these annular powers, proving the derivative assertions.
\end{proof}

\subsection{Approximation error and the coarse Pohozaev estimate}\label{sec:ac-shrinking}

Let $\mck:=\Span\{UY_a\}$ be the space spanned by the bubble Jacobi fields,
where $Y_1,\ldots,Y_{n+1}$ form an $L^2(\Ss^n,g_c)$-orthonormal basis of the first spherical harmonics and are pulled back by stereographic projection.
We take the orthogonal complement $\mck^\perp$ in the energy space $\mce$ of \eqref{eq:scalar-energy-norm}.
Orthogonality to $\mck$ in the energy inner product is equivalently expressed by
\[
\mck^\perp
=\left\{\phi\in\mce:
\int_{\R^n}\phi U^{\frac{n+2N}{n-2N}}Y_a\,dy=0
\text{ for }a\in\{1,\ldots,n+1\}\right\},
\]
because $P_0(UY_a)=(n+2N)\kappa_0U^{(n+2N)/(n-2N)}Y_a/(n-2N)$.
Let $\mcl_c^\dagger$ be the round pseudoinverse of $P_{2N,g_c}-(n+2N)\kappa_0/(n-2N)$: it vanishes on the first harmonics and inverts the operator on their orthogonal complement.
The self-adjoint reduced inverse of $L_U$ is
\[
G_Uf:=U\mcl_c^\dagger
\left(\frac{f}{U^{\frac{n+2N}{n-2N}}}\right).
\]
Since $dv_{g_c}=U^{2n/(n-2N)}\,dy$, conformal covariance gives distinct projections on the source and solution spaces:
\begin{equation}\label{eq:ac-source-projection}
L_UG_Uf=(I-\Pi_{\msrc})f,\qquad
\Pi_{\msrc}f:=U^{\frac{n+2N}{n-2N}}\sum_{a=1}^{n+1}Y_a\int_{\R^n}fUY_a\,dy.
\end{equation}
\begin{equation}\label{eq:ac-solution-projection}
G_UL_U\phi=(I-\Pi_{\msol})\phi,\qquad
\Pi_{\msol}\phi:=\sum_{a=1}^{n+1}UY_a\int_{\R^n}\phi U^{\frac{n+2N}{n-2N}}Y_a\,dy.
\end{equation}
These identities first hold for smooth data and extend to the energy spaces by duality.
The spectrum \eqref{eq:round-scalar-spectrum} and the kernel construction in Proposition \ref{thm:repair-weighted-inverse} give
\begin{equation}
\|G_Uf\|_{\mce}\le C\|f\|_{\mce^*} \quad\text{and}\quad 
|G_U(y,z)|\le C|y-z|^{2N-n}.                             \label{eq:ac-display-3-1}
\end{equation}

\begin{lemma}
\label{lem:ac-cap}
Suppose $\mcb_R$ is an operator of order at most $2N$, supported in $B_{4R}$, whose metric part has coefficient estimates
\[
|\pp^j b_a(y)|\le C\mu^2(1+|y|)^{2+a-2N-j} \quad\text{for }  0\le a\le2N.
\]
Here $b_a$ denotes any coefficient multiplying a derivative of order $a$, and $j$ is a nonnegative derivative order.
Assume its additional zeroth-order part is bounded, in scaled $C^{0,\vth}$ norm, by $C\vep U^{4N/(n-2N)}$.
Assume $\rho=\mu R$ and $\vep$ are sufficiently small.
If $e$ is supported in $A_R^*$ and has scaled $C^{0,\vth}$ norm at most $AR^{-2N}$, the solution on $\mck^\perp$ of
\[
b=-G_U\mcb_Rb+G_Ue
\]
satisfies
\begin{equation}
\|b\|_{\mce}\le C A R^{\frac{n-2N}{2}} \quad\text{and}\quad 
|\pp^j b(y)|\le C A(1+|y|)^{-j}\quad(|y|\le2R,\ j\le2N).
\label{eq:ac-display-3-2}
\end{equation}
The corresponding scaled H\"older estimate also holds.
Outside $B_{2R}$, the pointwise bound includes the additional decay factor $(R/|y|)^{n-2N}$.
\end{lemma}
\begin{proof}
We set $W_R(y):=\min\{1,(R/(1+|y|))^{n-2N}\}$.
We use the Banach norm formed from the supremum norms of $W_R^{-1}(1+|y|)^j\pp^jb$ and the corresponding scaled H\"older seminorms.
The kernel estimate \eqref{eq:ac-display-3-1}, followed by interior Schauder estimates on annuli rescaled to unit size, bounds $G_Ue$ by $CA$ in this weighted norm.
The source norm in $\mce^*$ is at most $CAR^{(n-2N)/2}$ by dilation and Sobolev duality.

For the metric part, the pointwise source is bounded by $C\mu^2(1+|y|)^{2-2N}W_R(y)\|b\|$ on $B_{4R}$.
Splitting the Riesz integral into $|z|\ll|y|$, $|z-y|\ll|y|$, and their complement gives
\[
\int_{|z|<4R}|y-z|^{2N-n}\mu^2(1+|z|)^{2-2N}W_R(z)\,dz
\le C\rho^2 W_R(y).
\]
Here the radial integral is $O(\mu^2\int_1^{4R}t\,dt)=O(\rho^2)$.
For the potential part the corresponding bound is
\[
\int_{\R^n} |y-z|^{2N-n}U(z)^{\frac{4N}{n-2N}}W_R(z)\,dz\le C W_R(y);
\]
this follows from $U^{4N/(n-2N)}=O((1+|z|)^{-4N})$ and $n-2N>0$.
Annular Schauder estimates give the identical operator bounds in the full derivative/H\"older norm.
The finite-rank source projection in \eqref{eq:ac-source-projection} produces smooth terms with the same bounds; their coefficients are estimated by the same convergent integrals against $UY_a$.
Consequently $G_U\mcb_R$ has norm at most $C(\rho^2+\vep)$, uniformly in $R$, in this Banach space.

In $\mce$, the form estimate \eqref{eq:ac-form} and the Hardy bound for $U^{4N/(n-2N)}$ give the same small operator norm.
The Neumann series converges in both spaces to the same solution and proves \eqref{eq:ac-display-3-2}.
\end{proof}

\begin{lemma}
\label{lem:ac-decomposition}
Let $v$ be the rescaling \eqref{eq:ac-local-bubble} of a positive solution of \eqref{eq:constant-q-equation} with $\kappa=\kappa_0$.
Suppose that $|v/U-1|\le\vep$ on $B_{4R}$ and that $v$ satisfies the derivative bounds furnished by \cite[Lemma 4.3]{MP} on interior annuli.
Let $\psi$ be defined by \eqref{eq:ac-display-2-1}, and let $A$ be the scaled annular norm of $v-U-\psi$ in \eqref{eq:ac-display-2-5}.
Choose a smooth radial cutoff function $0\le\chi_0\le1$, equal to one on $B_{2R}$ and supported in $B_{(5/2)R}$.
After adjusting the bubble center and scale so that $\chi_0(v-U-\psi)$ is orthogonal to the Jacobi fields in the energy inner product, as specified by \eqref{eq:ac-solution-projection}, we have $\chi_0(v-U-\psi)=a+b$.
Here $a$ solves the reduced equation with the forcing from the metric perturbation and its associated nonlinear terms, and $b$ solves it with the forcing from the cutoff commutators, with
\begin{align}
\|a\|_{\mce}
&\le C\rho^2\|K\|_{\mch}
+C\|K-K_{\mathrm{tr}}\|_{\mch}
+C\mu^{\frac{n-2N}{2}}\rho^{T-\frac{n-2N}{2}},                \label{eq:ac-metric-rem}\\
\|b\|_{\mce}&\le C A R^{\frac{n-2N}{2}} \quad\text{and}\quad 
|\pp^jb|\le CA(1+|y|)^{-j}\quad(|y|\le2R).       \label{eq:ac-cap-rem}
\end{align}
The last term in \eqref{eq:ac-metric-rem} bounds the part of $a$ generated by the metric Taylor remainder.
Replacing the operator of the rescaled metric by $P(K)$, the operator of the metric with truncated logarithm $K$, in a localized Pohozaev integral contributes an additional $O(\mu^{n-2N}\rho^{T+2N-n})$.
After the center and scale have been fixed, multiplying $v-U-\psi$ by a smooth radial function supported in a smaller ball need not preserve orthogonality to $\mck$; Proposition \ref{prop:ac-fine} estimates the resulting projection onto $\mck$ and its contribution to the reduced equation.
\end{lemma}
\begin{proof}
The cutoffs are supported strictly inside $B_{4R}$, so the required derivatives are controlled by interior estimates on $A_R^*$.
Equation \eqref{eq:ac-corrector-decay} implies $|\pp^j\psi|\le C\rho^2U(1+|y|)^{-j}$ on the ball.
After multiplication by any of the radial cutoff functions used below, $\psi$ has energy norm bounded by $C\|K_{\mathrm{tr}}\|_{\mch}$ by \eqref{eq:ac-uniform-corrector-cutoff}, with the outer metric cutoff included in the norm.

We write the logarithm of the actual rescaled metric as $K+K_T$ on $B_{4R}$ and apply the same outer cutoff to $K_T$.
Taylor's theorem, including the trace and cutoff derivatives, yields
\begin{align*}
|\pp^jK_T(y)|&\le C\mu^T(1+|y|)^{T-j} \quad (0\le j\le2N+2),\\
\|K_T\|_{\mch}^2
&\le C\mu^{2T}\int_0^{4R}(1+r)^{2T+2N-n-1}\,dr
\le C\mu^{n-2N}\rho^{2T+2N-n}.
\end{align*}
The integration-by-parts argument of Lemma \ref{lem:ac-tame}, which also applies to the trace of $K_T$, and the differentiated bounds for $U+\psi$ imply
\begin{align*}
\|[P(K+K_T)-P(K)](U+\psi)\|_{\mce^*}
\le C\|K_T\|_{\mch}
\le C\mu^{\frac{n-2N}{2}}\rho^{T-\frac{n-2N}{2}}.
\end{align*}
Since $\mu(1+|y|)\le C\rho$ on $B_{4R}$, the Taylor remainder also satisfies the coefficient bounds required in Lemma \ref{lem:ac-cap}.
For the approximation error $w=v-U-\psi$, subtract the equation of $U$ and $L_U\psi=-DP_0[K_{\mathrm{tr}}]U$.
The equation for $w$ consists of: (i) the linear source $-DP_0[K-K_{\mathrm{tr}}]U$ from the uncorrected Taylor degrees; (ii) terms at least quadratic in $(K,\psi)$;
(iii) $(P_0-P(K+K_T))w$; (iv) $-[P(K+K_T)-P(K)](U+\psi)$; and (v) a multiplication operator times $w$ from the nonlinear Taylor remainder.
We define
\[
\mcr(t):=(U+t)^{\frac{n+2N}{n-2N}}-U^{\frac{n+2N}{n-2N}}
-\frac{n+2N}{n-2N}U^{\frac{4N}{n-2N}}t.
\]
For the multiplication operator in (v), use the exact divided difference of $\mcr$ between $t=\psi$ and $t=v-U$.
Its magnitude is at most $C\vep U^{4N/(n-2N)}$ once $\rho^2\le\vep$, since both $U+\psi$ and $v$ are positive and differ from $U$ by at most $C\vep U$.
Scaled elliptic estimates for $v-U$, whose supremum is at most $\vep U$ and whose metric coefficients differ from the Euclidean coefficients by $O(\rho^2)$,
give relative differentiated bounds of size $C(\vep+\rho^2)$ on $B_{3R}$, including $A_R^*$.
Thus the divided difference has scaled H\"older norm at most $C\vep U^{4N/(n-2N)}$ when $\rho^2\le\vep$.

Multiply the equation for $w$ by the function $\chi_0$ chosen in the statement.
Commutators of $L_U$ and of $P(K+K_T)-P_0$ with $\chi_0$ are supported in $A_R^*$, involve at most $2N-1$ derivatives of $w$, and have scaled source size at most $CAR^{-2N}$.
Choose a smooth radial cutoff $\chi_1$ equal to one on the support of $\chi_0$ and supported in $B_{3R}$.
Use $\chi_1$ to extend the nonlinear coefficient defined by the divided difference by zero, and extend the metric coefficients with the metric cutoff supported in $B_{4R}$.
We collect the terms linear in $w$ into $-\mcb_Rw$, with the cutoff extensions just specified.
The product $\chi_0w$ satisfies the exact equation
\[
L_U(\chi_0w)=-\mcb_R(\chi_0w)+F+e,
\]
where $F$ collects the forcing terms independent of $w$ and $e$ consists of the cutoff commutators.
Applying \eqref{eq:ac-solution-projection} gives
\[
(I-\Pi_{\msol})(\chi_0w)
=-G_U\mcb_R(\chi_0w)+G_U(F+e).
\]
The center and scale adjustment at the end of this proof imposes $\Pi_{\msol}(\chi_0w)=0$, so the projected equation becomes
\[
\chi_0w=-G_U\mcb_R(\chi_0w)+G_U(F+e).
\]
Lemma \ref{lem:ac-tame}, \eqref{eq:ac-linear-source}, and the relative nonlinear Taylor inequality give
\[
\|F\|_{\mce^*}
\le C\rho^2\|K\|_{\mch}+C\|K-K_{\mathrm{tr}}\|_{\mch}
+C\mu^{\frac{n-2N}{2}}\rho^{T-\frac{n-2N}{2}}.
\]
The coefficient $C\vep U^{4N/(n-2N)}$ is included in $\mcb_R$.
On $B_{4R}$, $|\psi|\le C\rho^2U$, so the relative Taylor formula gives
\[
|\mcr(\psi)|\le C\rho^2U^{\frac{4N}{n-2N}}|\psi|.
\]
The bubble potential form bound and \eqref{eq:ac-uniform-corrector-cutoff} imply
\[
\|\chi_0\mcr(\psi)\|_{\mce^*}
\le C\rho^2\|\chi_1\psi\|_{\mce}
\le C\rho^2\|K_{\mathrm{tr}}\|_{\mch}.
\]
This applies Taylor's formula to $s\mapsto s^{(n+2N)/(n-2N)}$ for $s=(U+\psi)/U$ in a fixed compact subinterval of $(0,\infty)$, and holds in every dimension $n>2N$.
Solve $a=G_U(F-\mcb_Ra)$ and $b=G_U(e-\mcb_Rb)$ in $\mck^\perp$, keeping the same coefficients in the small operator $\mcb_R$ for both equations.
Their differential equations are
\[
L_Ua=(I-\Pi_{\msrc})(F-\mcb_Ra) \quad\text{and}\quad 
L_Ub=(I-\Pi_{\msrc})(e-\mcb_Rb).
\]
Uniqueness of the contraction gives $a+b=\chi_0w$; the sum of the two unprojected right-hand sides is the actual source $L_U(\chi_0w)$ and hence has zero pairing with every Jacobi field.
Applying \eqref{eq:ac-display-3-1} to the equation for $a$, and absorbing $G_U\mcb_Ra$ using its small operator norm, proves \eqref{eq:ac-metric-rem}.
Lemma \ref{lem:ac-cap}, applied to the equation for $b$, proves \eqref{eq:ac-cap-rem}.

To impose $\Pi_{\msol}(\chi_0w)=0$, adjust the bubble scale and center and recenter the conformal-normal chart.
Consider the $n+1$ pairings of $\chi_0(v-U)$ against $U^{(n+2N)/(n-2N)}Y_a$, for $a\in\{1,\ldots,n+1\}$.
At the bubble, their derivative with respect to the scale divided by $\mu$ and the center displacement divided by $\mu$ is the Gram matrix of the Jacobi fields with weight $\chi_0U^{4N/(n-2N)}$, up to fixed nonzero column factors, and is therefore invertible.
Local bubble convergence and the integrable tail of $U^{2n/(n-2N)}$ allow the implicit-function theorem.
Along the blow-up sequence, the adjusted scale divided by the initial scale tends to one, as in \eqref{eq:bubble-adjusted-scale}, and the center displacement divided by the initial scale tends to zero.
By Proposition \ref{thm:finite-SVT-exhaustion-new}, the polynomial scalar correction terms have scalar angular degrees at least two,
so their integrals against $U^{(n+2N)/(n-2N)}Y_a$ vanish after multiplication by $\chi_0$.
Therefore the center and scale chosen to make $\chi_0(v-U)$ orthogonal to the Jacobi fields also satisfy $\Pi_{\msol}(\chi_0w)=0$.

For a fixed smooth radial cutoff $\chi$ supported in $B_1$, the coefficient bounds for $K_T$ and the differentiated bubble bounds for $v$ imply
\begin{align*}
\left|\int_{B_R}\chi\left(\frac{y}{R}\right)(\msZ v)[P(K+K_T)-P(K)]v\,dy\right|
&\le C\mu^T\int_0^{4R}(1+r)^{T+2N-n-1}\,dr\\
&\le C\mu^{n-2N}\rho^{T+2N-n}.
\end{align*}
The assumptions $T>n+4$ and $R\to\infty$ justify the differentiations and both Taylor-remainder estimates.
\end{proof}

\medskip
We next relate the quadratic part of the localized Pohozaev integral to the dilation derivative of the energy form in \eqref{eq:schur-york-identity-new}.
For trace-free compactly supported $K$, coordinate-dilation covariance and self-adjointness with respect to $dy$ give
\begin{equation}
\int_{\R^n} (\msZ f)(P_0-P(K))f
=\frac{1}{2}\langle f,DP(K)[y\cdot\pp K]f\rangle. \label{eq:ac-display-4-1}
\end{equation}
Indeed $\msZ ^*=-\msZ-2N$ and $[P(K),\msZ]-2NP(K)=-DP(K)[y\cdot\pp K]$.
For $A_K:=P(K)-P_0$, the localized version is the exact identity
\begin{equation}\label{eq:ac-localized-dilation}
\begin{aligned}
2\int_{\R^n}\chi(\msZ f)(P_0-P(K))f
&=\langle f,\chi\dot A_K f\rangle
-\langle f,[A_K,\chi]\msZ f\rangle\\
&\quad+\langle f,(y\cdot\pp\chi)A_Kf\rangle,
\end{aligned}
\end{equation}
where $\dot A_K:=DP(K)[y\cdot\pp K]$.
The last two operators have total order at most $2N$ and are supported in the cutoff annulus; the commutator loses one derivative before $\msZ$ adds one.
Integration by parts bounds their weak forms using only $N$ derivatives of each energy argument.

At second order, when $L_U\psi=-DP_0[K]U$, \eqref{eq:ac-display-4-1} becomes
\begin{equation}
\frac{1}{2}\langle U,D^2P_0[K,y\cdot\pp K]U\rangle
+\langle\psi,DP_0[y\cdot\pp K]U\rangle.       \label{eq:ac-display-4-2}
\end{equation}
For $K_\tau(y):=K(e^\tau y)$, expression \eqref{eq:ac-display-4-2} equals $|\Ss^n|^{(n-2N)/n}/4$ times the derivative at $\tau=0$ of
\[
\left.\pp_t^2\right|_{t=0}\mce_N(\exp(tK_\tau),U+t\Psi[K_\tau]),
\]
whose value is given by \eqref{eq:schur-york-identity-new}.
To identify the mixed term, differentiate $L_U\psi=-DP_0[K]U$ and use self-adjointness of $L_U$.
Thus, in each angular block, \eqref{eq:ac-display-4-2} is a positive multiple of the quadratic form with matrix $S^{\mhom}=-H^+$, by Propositions \ref{thm:schur-york-identity-new} and \ref{prop:repair-dilation}.
For the truncated scalar correction term in \eqref{eq:ac-display-2-1}, the representation by the matrix $S^{\mhom}$ applies to $K_{\mathrm{tr}}$; the contributions of the remaining Taylor degrees are estimated in Lemma \ref{lem:ac-quadratic}.
The linear term vanishes for every $H_d\in\mfJ_d$: its average over a centered sphere is an $O(n)$-invariant linear functional of $H_d$, whereas the decomposition in Proposition \ref{thm:finite-SVT-exhaustion-new} contains no angular degree-zero summand.
The same cancellation holds after radial integration.

\medskip
Fix a smooth nonincreasing $\chi$ with $\chi=1$ near zero and $\chi=0$ beyond one.
We define the localized Pohozaev integral
\begin{equation}
\mscD_i(\rho):=\int_{B_\rho}\chi\left(\frac{x}{\rho}\right)(\msZ u_i)(P_0-P_{2N,g_i})u_i\,dx.
\label{eq:ac-display-4-3}
\end{equation}
Lemma \ref{lem:completion-flat-Pohozaev}, applied to the local representative specified in Section \ref{sec:ac-tools}, gives
\begin{equation}
\mscD_i(\rho)=
\int_0^\infty-\chi'(t)\mcp_N^{\kappa_0}(t\rho,u_i)\,dt.
\label{eq:ac-display-4-4}
\end{equation}

\begin{lemma}
\label{lem:ac-quadratic}
Assume \eqref{eq:established-A} and consider the model $U+\psi$ on $B_{4R}$, where $R=\rho/\mu\to\infty$.
If the Taylor coefficients are uniformly bounded, then, for every $\theta>0$, the quadratic Taylor term of \eqref{eq:ac-display-4-3} at $(\delta,U)$, evaluated on the metric and scalar perturbations $(K,\psi)$, is bounded above by
\[
-cJ_\mu(\rho)^2+\theta J_\mu(\rho)^2
+C_\theta R^{2N-n}\rho^4
+C_\theta\mu^{n-2N}\rho^{2d_0+2N-n}.
\]
The last term bounds the contributions from Taylor degrees greater than $(n-2N)/2$ after their mixed pairings are absorbed into $\theta J_\mu(\rho)^2$.

For a sequence of scales and Taylor coefficients, write $J_i(\rho):=J_{\mu_i}(\rho)$ and omit $\rho$ when it is fixed.
If, for fixed $\rho$,
\[
J_i^2=O(\mu_i^{n-2N}) \quad\text{and}\quad 
H_{i,d}\to0\quad\left(2\le d\le\frac{n-2N}{2}\right)
\quad\text{as }i\to\infty,
\]
then multiplying the metric and scalar perturbations of degrees below $(n-2N)/2$ by the cutoff functions changes the quadratic term by $o_i(\mu_i^{n-2N})$.
At the critical degree, when present, the term remaining after subtraction of the leading logarithmic term has the same bound.
\end{lemma}
\begin{proof}
For a degree $d<(n-2N)/2$, the difference $\eta_R\mu^dH_d-\mu^dH_d$, with $\eta_R$ as in Lemma \ref{lem:ac-uniform-cutoffs}, has form norm at most $C\mu^d|H_d|R^{d-(n-2N)/2}$.
The scalar correction term has the same energy tail, by \eqref{eq:ac-corrector-decay} and the reduced-inverse bounds \eqref{eq:ac-display-3-1}.
Put
\begin{equation}\label{eq:subcritical-cutoff-tail}
t_i(\rho):=\mu_i^{\frac{n-2N}{2}}\sum_{2\le d<\frac{n-2N}{2}}
|H_{i,d}|\rho^{d-\frac{n-2N}{2}}.
\end{equation}
Continuity of the localized polarized form, including terms with derivatives of the cutoff functions, bounds the change in the quadratic form caused by these multiplications by $C(J_it_i+t_i^2)$.
Young's inequality gives the first claimed estimate, because $t_i^2\le C\mu_i^{n-2N}\rho^{4+2N-n}=CR_i^{2N-n}\rho^4$ when the Taylor coefficients are uniformly bounded.
If $J_i=O(\mu_i^{(n-2N)/2})$ and $H_{i,d}\to0$ as $i\to\infty$, the expression $C(J_it_i+t_i^2)$ divided by $\mu_i^{n-2N}$
tends to zero at each fixed $\rho$ directly, before applying Young's inequality.

At degree $d=(n-2N)/2$, the quadratic model term has a negative definite coefficient of $\log R$ by Lemma \ref{lem:ac-endpoint}.
Pairings of such a term with one of degree below $(n-2N)/2$ converge as $R\to\infty$.
The bounded term after subtracting the logarithmic term is $O(\mu^{n-2N}|H_{(n-2N)/2}|^2)$.
For the cutoff comparison at the critical degree, write $n=2N+2d$ and fix the term $\mu^dH_d$.
Use the cutoffs $\chi_R,\eta_R$ of Lemma \ref{lem:ac-uniform-cutoffs}.
Compare $\Psi[\eta_R\mu^dH_d]$, obtained by multiplying the metric variation by $\eta_R$ before solving its scalar equation, with $\chi_R\Psi[\mu^dH_d]$, obtained by multiplying the solution of the original scalar equation by $\chi_R$.
The source of $\Psi[\eta_R\mu^dH_d]-\Psi[\mu^dH_d]$ is supported outside $B_{2R}$ and is bounded by $C\mu^d|H_d|(1+|y|)^{d-n}$, including where derivatives fall on $\eta_R$.
The kernel bound in \eqref{eq:ac-display-3-1}, integrated over $|z|>2R$, and interior estimates give
\[
|\pp^j(\Psi[\eta_R\mu^dH_d]-\Psi[\mu^dH_d])(y)|
\le C\mu^d|H_d|R^{-d}(1+|y|)^{-j} \quad\text{for }  |y|<R \quad\text{and}\quad  0\le j\le2N.
\]
Pairing this difference with a degree-$d$ metric source on the interior is bounded by
\[
C\mu^{2d}|H_d|^2R^{-d}
\left(1+\int_1^R t^{d-1}\,dt\right)
\le C\mu^{2d}|H_d|^2.
\]
After integration by parts, differentiated pairings have the same total radial power.
On annuli with $|y|\asymp R$, both scalar correction terms have derivative bounds $C\mu^d|H_d|R^{-d-j}$ and hence give the same bound by dilation.
For $r=|y|>8R$, the derivatives of $\Psi[\eta_R\mu^dH_d]$ satisfy the bounds $C\mu^d|H_d|R^d r^{-2d-j}$, whereas $\chi_R\Psi[\mu^dH_d]=0$.
Its exterior energy is bounded by
\[
C\mu^{2d}|H_d|^2R^{2d}
\int_{8R}^{\infty}r^{-2d-1}\,dr
\le C\mu^{2d}|H_d|^2.
\]
The same estimates hold after $R\pp_R$, since the differentiated cutoff source is supported where $|y|\asymp R$ with the same scale.
Replacing $\Psi[\eta_R\mu^dH_d]$ by $\chi_R\Psi[\mu^dH_d]$ therefore changes the bounded term by at most $C\mu^{n-2N}|H_{(n-2N)/2}|^2$ and leaves the logarithmic coefficient unchanged.
Young's inequality absorbs the bounded mixed pairings into the definite diagonal terms below and at the critical degree, whose negative sign in the sign convention for the local Pohozaev integral controls $-J_\mu(\rho)^2$.
The bounded error divided by $\mu_i^{n-2N}$ vanishes when $H_{i,(n-2N)/2}\to0$ as $i\to\infty$.

For a Taylor term of degree $e>(n-2N)/2$, estimate \eqref{eq:ac-uniform-jet-norm} gives $\|\eta_R\mu^eH_e\|_{\mch}^2\le C\mu^{n-2N}\allowbreak\rho^{2e+2N-n}|H_e|^2$.
Applying the form bounds and Young's inequality to each pairing with a term of degree at most $(n-2N)/2$ yields $\theta J_i^2+C_\theta\mu_i^{n-2N}\rho^{2d_0+2N-n}$.
The same argument applies to the scalar correction terms of degrees less than $n-2N$, since the bounds for $\chi_R\psi$ in \eqref{eq:ac-uniform-corrector-cutoff} have the same power.
\end{proof}

\begin{proposition}
\label{prop:ac-coarse}
Assume \eqref{eq:established-A} and the relative bubble estimate \eqref{eq:ac-relative-bubble}, and set $R_i:=\rho_i/\mu_i$.
For sufficiently small $\rho_i$, uniformly in the moving center,
\begin{equation}
\mscD_i(\rho_i)
\le -cJ_i(\rho_i)^2
+C\rho_i^2\left(\frac{\mu_i}{\rho_i}\right)^{n-2N}
+o_i(J_i(\rho_i)^2).                            \label{eq:ac-display-4-5}
\end{equation}
Changing the cutoff radii by fixed numerical factors only changes $c,C$.
Consequently, if $\rho_i\to0$ as $i\to\infty$,
\begin{equation}
\limsup_{i\to\infty}\left(\frac{\rho_i}{\mu_i}\right)^{n-2N}\mscD_i(\rho_i)\le0.     \label{eq:ac-display-4-6}
\end{equation}
\end{proposition}
\begin{proof}
Apply Lemma \ref{lem:ac-decomposition} with $A\le CR_i^{2N-n}$, retaining $A$ in the estimates needed later in Proposition \ref{prop:ac-fine}.
The part $b$ generated by the cutoff commutators has energy norm at most $CR_i^{-(n-2N)/2}$; to estimate the term linear in $b$, use the sharper pointwise bound \eqref{eq:ac-display-3-2}.
That term has on an annulus $|y|\simeq t$ the bound
\[
C A\,\mu_i^2(1+t)^2;
\]
this follows from $U=O(t^{2N-n})$, total differential order $2N$, and the $n$-dimensional annular volume element.
The radial integral is bounded by
\begin{equation}\label{eq:ac-cutoff-linear-error}
CA\mu_i^2\left(1+\int_1^{4R_i}t\,dt\right)
\le CA\rho_i^2.
\end{equation}
By \eqref{eq:ac-form} and \eqref{eq:ac-cap-rem}, the terms quadratic in $b$ are bounded by
\begin{equation}\label{eq:ac-cutoff-quadratic-error}
C\rho_i^2\|b\|_{\mce}^2\le C\rho_i^2 A^2R_i^{n-2N}.
\end{equation}
With $A\le CR_i^{2N-n}$, both \eqref{eq:ac-cutoff-linear-error} and \eqref{eq:ac-cutoff-quadratic-error} are bounded by $C\rho_i^2R_i^{2N-n}$.
Bound the mixed term involving $\psi$ and $b$ by \eqref{eq:ac-uniform-corrector-cutoff} and Young's inequality:
\[
C\rho_i^2\|K\|_{\mch}\|b\|_{\mce}
\le C\rho_i^2\|K\|_{\mch}^2
+C\rho_i^2A^2R_i^{n-2N}.
\]
The same estimates include the cutoff commutators in \eqref{eq:ac-localized-dilation}.

The part of the approximation error generated by the polynomial metric perturbation contributes at most $C\rho_i^2\|K\|_{\mch}^2$,
together with pairings containing Taylor terms of degree greater than $(n-2N)/2$, which are controlled by Young's inequality.
By \eqref{eq:ac-metric-rem} and linearity of the inverse with $\mcb_R$ fixed, the correction generated by the metric Taylor remainder $K_T$ contributes at most
\begin{align*}
C\mu_i^{\frac{n-2N}{2}}\rho_i^{T-\frac{n-2N}{2}}\|K\|_{\mch}
\le \theta\|K\|_{\mch}^2
+C_\theta\mu_i^{n-2N}\rho_i^{2T+2N-n}.
\end{align*}
Quadratic pairings of this correction with the approximation error are bounded by the same expression and by the bound $C\rho_i^2A^2R_i^{n-2N}$ for terms quadratic in $b$.
The direct operator error is $C\mu_i^{n-2N}\rho_i^{T+2N-n}$ by Lemma \ref{lem:ac-decomposition}.
These estimates use the localized dilation identity \eqref{eq:ac-localized-dilation}.
Cubic metric terms satisfy the same bounds.
A smooth radial cutoff gives identical bounds on its support and cutoff annulus, since its derivatives preserve the total differential order in the estimates.
The first-order term has zero radial average because the decomposition in Proposition \ref{thm:finite-SVT-exhaustion-new} has no angular degree-zero summand.
Lemma \ref{lem:ac-quadratic} identifies the remaining quadratic part and its sign.
Finally, \eqref{eq:ac-uniform-jet-norm} and \eqref{eq:ac-display-2-4} imply
\[
\|K\|_{\mch}^2\le C J_i^2+C\mu_i^{n-2N}\rho_i^{2d_0+2N-n} \quad\text{and}\quad 
\mu_i^{n-2N}\rho_i^{2d_0+2N-n}=R_i^{2N-n}\rho_i^{2d_0}.
\]
Choose the Young parameter and then $\rho_0$ small so that the positive $J_i^2$ error terms can be absorbed into the negative quadratic contribution from Lemma \ref{lem:ac-quadratic}.
The additional Taylor errors satisfy
\begin{align*}
\mu_i^{n-2N}\rho_i^{2T+2N-n}+\mu_i^{n-2N}\rho_i^{T+2N-n}
=R_i^{2N-n}(\rho_i^{2T}+\rho_i^T)
\le C\rho_i^2R_i^{2N-n}.
\end{align*}
The other errors have the same bound because $d_0\ge2$.
Convergence of the finite matrices representing the cutoff forms supplies the term $o_i(J_i^2)$ in \eqref{eq:ac-display-4-5}, including pairings between terms of degree below and equal to $(n-2N)/2$.
This proves \eqref{eq:ac-display-4-5}.
For sufficiently large $i\in\mathbb N$, its $o_i(J_i^2)$ term is bounded above by $cJ_i^2/2$.
Dropping the remaining negative term, multiplying by $(\rho_i/\mu_i)^{n-2N}$, and letting $i\to\infty$ with $\rho_i\to0$ proves \eqref{eq:ac-display-4-6}.
\end{proof}

\subsection{Isolated simple blow-up and quantitative vanishing} \label{sec:ac-separation}
The argument that an isolated blow-up point is simple is due to Schoen; see \cite[Lemma 8.2]{KMS}.
For orders four and six, the corresponding conclusions are \cite[Corollary 7.3]{GKW} and \cite[Corollary 10.33]{GKW}.

\begin{proposition}
\label{prop:ac-saturation}
Assume \eqref{eq:established-A}, $\ker P_{2N,g}=0$, and $G_g>0$.
For local representatives satisfying \eqref{eq:ac-local-growth}, if $\rho_i\to0$, then
\begin{equation}\label{eq:ac-shrinking-comparison}
\sup_{|x|<\rho_i}\left|\frac{u_i(x)}{B_i(x)}-1\right|\to0 \quad \text{as }i\to\infty.
\end{equation}
At an isolated blow-up point, for every sufficiently small $\vep>0$, there is a fixed $\rho_0>0$ for which \eqref{eq:ac-relative-bubble} holds with $\rho_i=\rho_0$ for all sufficiently large $i\in\mathbb N$.
Every isolated blow-up point is isolated simple.
\end{proposition}
\begin{proof}
Fix a small $\vep>0$ and let $r_i\le\rho_i$ be the supremum of radii on which $|u_i/B_i-1|\le\vep$.
The bubble convergence in \cite[Lemma 4.2]{MP} implies $\mu_i/r_i\to0$.
If $r_i<\rho_i$ and $r_i\to0$, the rescaling limit in \cite[Proposition 4.4]{MP} is
\[
r_i^{n-2N}\mu_i^{-\frac{n-2N}{2}}u_i(r_i z)
\to A_0|z|^{2N-n}+H(z) ,\quad\text{where } 
P_0H=0 \quad\text{and}\quad  H(0)>0,
\]
in $C^{2N}_{\mathrm{loc}}(B_2\setminus\{0\})$, where $A_0:=2^{(n-2N)/2}$.
Apply \eqref{eq:ac-display-4-6} at $\tau r_i$ with $\tau>0$ fixed and $4\tau<1$.
After multiplication by $(r_i/\mu_i)^{n-2N}$, the nonlinear boundary term in \eqref{eq:nonlinear-Pohozaev-integral} tends to zero because it contains the additional factor $(\mu_i/r_i)^{2N}$.
Equations \eqref{eq:ac-display-4-4} and \eqref{eq:completion-Pohozaev-identity} therefore imply
\[
\mcp_N^0(\tau,A_0|z|^{2N-n}+H)\le0.
\]
The left-hand side is independent of $\tau$ by \eqref{eq:completion-Pohozaev-identity}, whereas \eqref{eq:completion-mass-limit} identifies its value as $(n-2N)A_0H(0)/(2\ga_{n,N})>0$.
Thus $r_i<\rho_i$ cannot hold along a subsequence with $r_i\to0$.
For shrinking $\rho_i$, applying this conclusion for each fixed sufficiently small $\vep>0$ proves \eqref{eq:ac-shrinking-comparison}; for fixed $\rho_i$, it supplies the asserted comparison on a smaller fixed ball.

To verify simplicity in Definition \ref{def:ac-isolated-simple}, use the spherical average \eqref{eq:ac-spherical-average}.
The rescaled bubble has one nondegenerate critical radius for its weighted spherical average, and local $C^2$ convergence gives exactly one critical point on $0<r<L\mu_i$ for fixed sufficiently large $L>1$.
On the remaining annuli, scaled interior estimates for the equation of $u_i-B_i$ give
\begin{equation}\label{eq:ac-simple-derivative}
r|\nabla(u_i-B_i)|\le C(\vep+r^2)B_i \quad (L\mu_i\le r=|x|\le\rho_0).
\end{equation}
Choosing $\vep$ and then $\rho_0$ small in \eqref{eq:ac-simple-derivative} shows that $(r^{(n-2N)/2}\bar u_i(r))'<0$ there; the conformal-normal volume error is absorbed by the $r^2$ term.
\end{proof}

\begin{proposition}
\label{prop:ac-separation}
Assume \eqref{eq:established-A}, $\ker P_{2N,g}=0$, and $G_g>0$.
For each $i\in\mathbb N$, let $\xi_{j,i}$, $1\le j\le m_i$, be the critical points selected for $u_i$ by \cite[Proposition 4.1]{MP}, where $m_i$ is their number.
They satisfy $d_g(\xi_{j,i},\xi_{k,i})\ge c>0$ for $j\ne k$, after passage to a subsequence.
Their number is uniformly bounded, and every blow-up point is isolated simple.
\end{proposition}
\begin{proof}
We apply the proof of \cite[Proposition 7.1]{MP}, replacing its dimension-dependent bubble comparison \cite[Proposition 6.4]{MP} by \eqref{eq:ac-shrinking-comparison}.
All remaining inputs---point selection, annular Harnack estimates, positivity of the Green function, and classification of entire positive solutions---hold for every $n>2N$ under the present hypotheses.
The argument for the pair of selected points with the smallest separation therefore gives the stated lower bound on pairwise distances.
This is the higher-order version of \cite[Proposition 8.3]{KMS}; the order-four and order-six conclusions also appear in \cite[Corollary 7.4]{GKW} and \cite[Corollary 10.34]{GKW}.
The uniform separation makes sufficiently small balls about the selected points disjoint; compactness of $M$ therefore bounds their number.
The point-selection estimate then gives \eqref{eq:ac-isolation-bound} on a ball of fixed radius about each blow-up center.
Thus every blow-up point is isolated in the sense of Definition \ref{def:ac-isolated}, and Proposition \ref{prop:ac-saturation} proves simplicity.
\end{proof}

\begin{corollary}
\label{cor:ac-jet-vanishing}
Under the assumptions of Proposition \ref{prop:ac-separation}, at each blow-up point one can choose a fixed radius $\rho_0>0$
for which \eqref{eq:ac-relative-bubble} holds on the ball of radius $4\rho_0$.
For each $i\in\mathbb N$, adjust the bubble scale and center so that $\Pi_{\msol}(\chi_0(v_i-U-\psi_i))=0$, using the function $\chi_0$ from Lemma \ref{lem:ac-decomposition} with $R=\rho_0/\mu_i$, and retain these parameters for every smaller radius $0<\rho\le\rho_0$.
In those moving conformal-normal charts,
\begin{equation}
J_i(\rho_0)^2=O(\mu_i^{n-2N}) \quad\text{and}\quad 
|H_{i,d}|=O\left(\mu_i^{\frac{n-2N}{2}-d}\right)\quad\left(2\le d<\frac{n-2N}{2}\right),
\label{eq:ac-display-5-2}
\end{equation}
and, when $(n-2N)/2$ is an integer at least two,
\begin{equation}
|H_{i,\frac{n-2N}{2}}|=O\!\left(\left(\log\frac{1}{\mu_i}\right)^{-\frac{1}{2}}\right).
\label{eq:ac-display-5-3}
\end{equation}
Consequently each limiting blow-up point belongs to $\mcz$.
\end{corollary}
\begin{proof}
Proposition \ref{prop:ac-separation} supplies disjoint fixed neighborhoods of the blow-up points.
Proposition \ref{prop:ac-saturation} then gives \eqref{eq:ac-relative-bubble} on a smaller fixed ball for a sufficiently small tolerance.
The differentiated bubble bounds give $|\mcp_N^{\kappa_0}(t\rho_0,u_i)|\le C_{\rho_0}\mu_i^{n-2N}$
for $t\in\supp\chi'$, the radii appearing in the average \eqref{eq:ac-display-4-4}.
Combining the resulting lower bound for the averaged Pohozaev integral with \eqref{eq:ac-display-4-5}, and absorbing its final $o(J_i^2)$, gives \eqref{eq:ac-display-5-2}.
The individual summands in \eqref{eq:ac-display-2-4} yield the coefficient estimates in \eqref{eq:ac-display-5-2}--\eqref{eq:ac-display-5-3}.
Smooth dependence of the conformal-normal family on its center shows that every homogeneous metric term of degree $2\le d\le(n-2N)/2$ vanishes at the limiting point.
Equation \eqref{eq:metric-weyl-jets} therefore places that point in $\mcz$ when $n\ge2N+4$; in lower dimensions $\mcz=M$ by definition.
\end{proof}

\begin{remark}[Local equations and varying metrics]\label{rem:local-scope}
Under \eqref{eq:established-A}, the coefficient estimates \eqref{eq:ac-display-5-2}--\eqref{eq:ac-display-5-3} also hold for positive solutions of the critical equation on a fixed ball with metrics $g_i\to g$ smoothly.
It suffices to assume, in the moving conformal-normal charts, the relative bubble estimate \eqref{eq:ac-relative-bubble} at a fixed sufficiently small radius, local $C^{2N}$ convergence of the rescaled solutions to $U$, and the uniform differentiated annular bounds used in Lemma \ref{lem:ac-decomposition}.
Indeed, these bounds imply $|\mcp_N^{\kappa_0}(t\rho_0,u_i)|\le C_{\rho_0}\mu_i^{n-2N}$ on the fixed annulus.
Combining this bound with \eqref{eq:ac-display-4-5} gives $J_i(\rho_0)^2=O(\mu_i^{n-2N})$, and each nonnegative summand in \eqref{eq:ac-display-2-4} gives the corresponding coefficient estimate.
The global Green function hypotheses are used to establish \eqref{eq:ac-isolation-bound} on fixed balls and obtain these local bounds; they are not needed once the bounds are assumed.

\medskip
We restrict our statements to the critical equation.
The method extends to subcritical powers approaching the critical exponent by the standard changes to the blow-up scaling and Pohozaev identity described in \cite[Sections 2--3]{KMS} and \cite[Section 4]{MP}.
\end{remark}

\section{Proof of the compactness theorem}\label{sec:green-matching}
By Theorem \ref{thm:repair-Green} and \eqref{eq:ac-display-7-1}, the bubble and its scalar correction reproduce the leading radial singularity and all singular angular terms of the limit expressed as a sum of Green functions.
Subtracting them leaves a bounded term, which gives the annular approximation estimate \eqref{eq:ac-display-7-3}.
Here $\rho$ is the radius in the original conformal-normal coordinates about the chosen blow-up center.
We first let $i\to\infty$ at fixed $\rho$ and then let $\rho\downarrow0$; this order of limits is retained throughout the refined Pohozaev estimate and the mass contradiction.

\subsection{The limiting Green functions and matching at fixed poles}\label{sec:ac-matching} 
\begin{proposition}\label{prop:ac-green-limit}
Assume \eqref{eq:established-A}, that $P_{2N,g}$ is invertible with positive Green function, and that $u_i>0$ solves \eqref{eq:constant-q-equation} with $\kappa=\kappa_0$ and $\|u_i\|_{L^\infty(M)}\to\infty$.
After applying Proposition \ref{prop:ac-separation} and discarding selected balls with uniformly bounded local supremum,
denote the remaining blow-up scales by $\mu_{j,i}$, $1\le j\le m$, so that $\mu_{j,i}\to0$ as $i\to\infty$.
Pass to a subsequence so that $\xi_{j,i}\to p_j$ as $i\to\infty$ and the scale ratios below converge.
Relabel so that
\[
\mu_i:=\max_{1\le j\le m}\mu_{j,i}=\mu_{1,i} \quad\text{and}\quad 
\tau_j:=\lim_{i\to\infty}\left(\frac{\mu_{j,i}}{\mu_i}\right)^{\frac{n-2N}{2}}\in[0,1] ,\quad\text{where } \tau_1=1.
\]
Then, on fixed closed sets away from the limiting blow-up points,
\begin{equation}
\mu_i^{-\frac{n-2N}{2}}u_i\to
v:=\frac{A_0}{\ga_{n,N}}\sum_{j=1}^m\tau_jG_g(p_j,\cdot)
\label{eq:ac-display-6-1}
\end{equation}
in every $C^k$ norm with $k\in\mathbb N\cup\{0\}$ as $i\to\infty$.
\end{proposition}
\begin{proof}
On a fixed annulus about a blow-up point, Corollary \ref{cor:ac-jet-vanishing} gives $u_i\le C\mu_{j,i}^{(n-2N)/2}\le C\mu_i^{(n-2N)/2}$.
The positive Green kernel has a positive global lower bound: it diverges positively at the diagonal, and is positive and continuous off it.
Evaluating the global Green representation at a point in that annulus yields
\[
\int_Mu_i^{\frac{n+2N}{n-2N}}\,dv_g\le C\mu_i^{\frac{n-2N}{2}}.
\]
Since $\sup_y\int_MG_g(x,y)\,dv_g(x)<\infty$, the same representation gives $\|u_i\|_{L^1}\le C\mu_i^{(n-2N)/2}$.
Bubble convergence at the first center gives the complementary bound $\int_Mu_i^{(n+2N)/(n-2N)}\,dv_g\ge c\mu_i^{(n-2N)/2}$;
the positive lower bound for $G_g$ therefore yields $u_i\ge c\mu_i^{(n-2N)/2}$ on $M$.
On compact sets away from selected points, the selection estimate bounds $u_i$ uniformly, so $P_{2N,g}-\kappa_0u_i^{4N/(n-2N)}$ is a uniformly elliptic operator with bounded zeroth-order coefficient.
Interior $L^1$ estimates give $u_i=O(\mu_i^{(n-2N)/2})$ there.
The normalized function $\wtv_i:=\mu_i^{-(n-2N)/2}u_i$ consequently lies in a fixed positive interval on each such compact set and satisfies
\[
P_{2N,g}\wtv_i=\kappa_0\mu_i^{2N}\wtv_i^{\frac{n+2N}{n-2N}}.
\]
Elliptic bootstrapping gives uniform local $C^k$ bounds.
A selected ball with uniformly bounded local supremum is treated by the same local estimate and therefore contributes no source term in the limit.

On a fixed ball around $\xi_{j,i}$ supplied by Proposition \ref{prop:ac-saturation},
\[
u_i(x)\le C\mu_{j,i}^{-\frac{n-2N}{2}}
\left(1+\frac{|x-\xi_{j,i}|}{\mu_{j,i}}\right)^{2N-n}.
\]
Consequently the integral of the nonlinear source within that fixed ball but outside the ball of radius $L\mu_{j,i}$, divided by $\mu_i^{(n-2N)/2}$, is at most $C L^{-2N}$ uniformly in $i$.
On the smaller ball use local bubble convergence.
It follows, first with $i\to\infty$ and then $L\to\infty$, that
\[
\kappa_0\mu_i^{-\frac{n-2N}{2}}u_i^{\frac{n+2N}{n-2N}}\,dv_g
\rightharpoonup \kappa_0\left(\int_{\R^n} U^{\frac{n+2N}{n-2N}}\right)
\sum_j\tau_j\delta_{p_j}.
\]
The source outside these fixed balls tends to zero, since there $u_i=O(\mu_i^{(n-2N)/2})$ and its normalized nonlinear source therefore contains the factor $\mu_i^{2N}$.
The representation of $U$ by the fundamental solution \eqref{eq:completion-fundamental-constant} gives
\begin{equation}\label{eq:ac-bubble-source-mass}
\kappa_0\int_{\R^n}U^{\frac{n+2N}{n-2N}}\,dy
=\frac{A_0}{\ga_{n,N}} ,\quad\text{where }  A_0=2^{\frac{n-2N}{2}}.
\end{equation}
Invertibility of $P_{2N,g}$ identifies the limit as \eqref{eq:ac-display-6-1}, and local elliptic estimates yield smooth convergence away from $p_1,\ldots,p_m$.
\end{proof}

\begin{lemma}
\label{lem:ac-corrector-asymptotic}
Let $H_d\in\mfJ_d$, $2\le d<n-2N$.
The scalar correction term $\Psi[H_d]$ specified in Corollary \ref{cor:repair-high-corrector} has the expansion
\begin{equation}
\Psi[H_d](y)=|y|^{d+2N-n}a_d\left(H_d,\frac{y}{|y|}\right)
+o_{C^{2N}}(|y|^{d+2N-n}),                    \label{eq:ac-display-7-1}
\end{equation}
as $|y|\to\infty$, where $a_d$ has zero spherical mean.
The function $a_d$ equals the angular coefficient in the expansion of the Green function associated with $H_d$,
when the leading singular term is normalized as $A_0r^{2N-n}$ and the same conformal normal coordinates are used.
\end{lemma}
\begin{proof}
The forcing is a finite sum of scalar spherical harmonics with radial coefficients.
Naturality and homogeneity give
\[
-DP_0[H_d]U=r^{d-n}f_d(H_d,\theta)+O_{C^{2N}}(r^{d-n-2}).
\]
Only scalar spherical harmonic components present in $\mfJ_d$ can occur.
Their degrees are at least two by Proposition \ref{thm:finite-SVT-exhaustion-new}, so $f_d$ has zero mean and no first-harmonic component.
For each angular degree $\ell$, the indicial polynomial $I_{N,n}(\lambda,\ell)$ in \eqref{eq:repair-indicial} is nonzero for $2N-n<\lambda<0$.
We define $a_d$ by $I_{N,n}(d+2N-n,\ell)a_{d,\ell}=f_{d,\ell}$, where the subscript $\ell$ denotes projection onto the degree-$\ell$ spherical harmonics.
We subtract $r^{d+2N-n}a_d$ with a cutoff near zero from $\Psi[H_d]$.
At infinity the remaining source decays faster than $r^{d-n}$ by a positive power; the potential $(n+2N)\kappa_0U^{4N/(n-2N)}/(n-2N)=O(r^{-4N})$ also contributes only terms of smaller order.
For each angular component, the radial equation in $t=\log r$ has constant-coefficient indicial part and coefficients decaying exponentially as $t\to\infty$.
The decay estimate \eqref{eq:ac-corrector-decay} excludes nondecaying indicial solutions; every decaying indicial term has homogeneity at most $2N-n$, so it is $o(r^{d+2N-n})$, including possible logarithms at lower powers.
Variation of constants for the associated finite first-order system, followed by scaled interior Schauder estimates, gives \eqref{eq:ac-display-7-1} with derivatives through order $2N$ and convergence in the scaled $C^{2N,\vth}$ norm on annuli.
The constants are uniform for $H_d$ on the unit sphere of $\mfJ_d$.

Since $U=A_0r^{2N-n}+O(r^{2N-n-2})$, the source $f_d$ is the negative metric variation applied to $A_0r^{2N-n}$.
The construction in Theorem \ref{thm:repair-Green} solves the same indicial equation with the same sign and singular normalization.
\end{proof}

\begin{proposition}
\label{prop:ac-matched-cap}
Fix a pole $p_j$ of \eqref{eq:ac-display-6-1} with $\tau_j>0$ and use its scale $\mu_i:=\mu_{j,i}$ in this proposition.
Retain the center and scale chosen at the fixed radius of Corollary \ref{cor:ac-jet-vanishing}.
Let $A_i(\rho)$ be the scaled annular norm in \eqref{eq:ac-display-2-5} of the approximation error $v_i-U-\sum_{2\le d<n-2N}\mu_i^d\Psi[H_{i,d}]$.
Then
\begin{equation}
\sup_{0<\rho<\rho_0}\limsup_{i\to\infty}
\mu_i^{2N-n}A_i(\rho)<\infty.      \label{eq:ac-display-7-3}
\end{equation}
\end{proposition}
\begin{proof}
We first take $j=1$, so $\mu_i=\mu_{1,i}$ and $\tau_1=1$.
We use the smooth conformal-normal chart family at the modulated centers.
On fixed annuli in the $x$-variables, the chart, conformal factor, and their finite derivatives converge to those at $p_1$.
Since $\tau_1=\Lambda_{p_1}(p_1)=1$, the conformal-normal representative $\Lambda_{p_1}^{-1}v$ of the limit in \eqref{eq:ac-display-6-1} has singular coefficient $A_0$ at $p_1$.
For $d_0$ as in \eqref{eq:main-mass-expansion} and $0<\vep<\min\{1,2d_0+2N-n\}$, Corollary \ref{cor:ac-jet-vanishing} and Theorem \ref{thm:repair-Green} imply
\begin{equation}
(\Lambda_{p_1}^{-1}v)(x)=A_0r^{2N-n}+\sum_{d=d_0}^{n-2N}r^{d+2N-n}a_d(\theta)
+B+O_{C^{2N}}(r^\vep).
\label{eq:ac-display-7-4}
\end{equation}
The Green functions with poles at the other blow-up points contribute smooth terms.
Every displayed angular term is linear in the homogeneous metric terms of degrees $d_0\le d\le n-2N$, because $2d_0>n-2N$.
For each integer $d$ with $2\le d<n-2N$, Lemma \ref{lem:ac-corrector-asymptotic} implies, on a fixed annulus in the $x$-variables as $i\to\infty$,
\[
\mu_i^{d+2N-n}\Psi[H_{i,d}]\left(\frac{x}{\mu_i}\right)
\to r^{d+2N-n}a_d(H_{p_1,d},\theta).
\]
Subtracting this convergence and $\mu_i^{2N-n}U(x/\mu_i)\to A_0r^{2N-n}$ from \eqref{eq:ac-display-7-4} removes every singular angular term and leaves
\[
B+a_{n-2N}(\theta)+O_{C^{2N}}(r^\vep).
\]
Here $a_{n-2N}=0$ when the angular sum in \eqref{eq:ac-display-7-4} is empty.
The scaled interior Schauder estimates used in the proof of Theorem \ref{thm:repair-Green} bound the scaled $C^{2N,\vth}$ norms of this expression uniformly as $\rho\downarrow0$.
This proves \eqref{eq:ac-display-7-3} at $p_1$.
For $\tau_j>0$, the scale ratio $(\mu_{j,i}/\mu_{1,i})^{(n-2N)/2}\to\tau_j$, \eqref{eq:ac-display-6-1}, and convergence of the charts and conformal factors give
\[
\mu_{j,i}^{-\frac{n-2N}{2}}u_i\to\frac{\Lambda_{p_j}^{-1}v}{\tau_j} \quad (i\to\infty)
\]
on fixed annuli, using the local representative at $\xi_{j,i}$ specified in Section \ref{sec:ac-tools}.
This limit again has singular coefficient $A_0$, so the same subtraction, using Corollary \ref{cor:ac-jet-vanishing} and Lemma \ref{lem:ac-corrector-asymptotic} at $p_j$, proves \eqref{eq:ac-display-7-3} with $\mu_i=\mu_{j,i}$.
\end{proof}

\subsection{The refined Pohozaev estimate and compactness}\label{sec:fine-completion}
\begin{proposition}\label{prop:ac-fine}
At a pole $p_j$ with $\tau_j>0$, use the scale $\mu_i:=\mu_{j,i}$ and the hypotheses of Proposition \ref{prop:ac-matched-cap}.
Set $\eta:=\min\{2,2d_0+2N-n\}>0$.
Under \eqref{eq:established-A}, there exist $c>0$ and $C<\infty$ such that, for each fixed $\rho\in(0,\rho_0)$ sufficiently small,
\begin{equation}
\mscD_i(\rho)
\le -cJ_i(\rho)^2+C\rho^\eta\mu_i^{n-2N}+o_i(\mu_i^{n-2N}).                    \label{eq:ac-display-8-1}
\end{equation}
The constant $C$ is independent of small $\rho$.
The last term divided by $\mu_i^{n-2N}$ tends to zero as $i\to\infty$ for each fixed $\rho\in(0,\rho_0)$; this convergence need not be uniform in $\rho$.
In particular
\begin{equation}
\limsup_{\rho\downarrow0}\limsup_{i\to\infty}
\mu_i^{2N-n}\mscD_i(\rho)\le0.\label{eq:ac-display-8-2}
\end{equation}
\end{proposition}
\begin{proof}
For each $i\in\mathbb N$, we retain the scale and center chosen at radius $\rho_0$ in Corollary \ref{cor:ac-jet-vanishing} and change only the cutoff radius $\rho$, with $R_i=\rho/\mu_i$.
Choose $\chi_0$ as in Lemma \ref{lem:ac-decomposition} with $R=R_i$, and set $\zeta_i:=\Pi_{\msol}(\chi_0(v_i-U-\psi_i))\in\mck$.
The coefficients of $\zeta_i$ are $O(R_i^{-n})$: the product formed with the function at radius $\rho_0/\mu_i$ has zero projection, the two functions agree on $B_{2R_i}$, and $v_i-U-\psi_i=O(U)$ on their supports.
Formula \eqref{eq:ac-solution-projection} and $\int_{|y|>R_i}U^{2n/(n-2N)}\,dy=O(R_i^{-n})$ give the bound.
Subtracting $\zeta_i$ adds the source $-\mcb_{R_i}\zeta_i$, whose correction is
\[
-(I+G_U\mcb_{R_i})^{-1}G_U\mcb_{R_i}\zeta_i.
\]
The weighted operator bound for $G_U\mcb_{R_i}$ proved in Lemma \ref{lem:ac-cap}, together with its energy bound and the Neumann series, bounds the induced correction in both norms by $C(\rho^2+\vep)R_i^{-n}$.
Since
\[
\mu_i^{2N-n}R_i^{-n}=\rho^{-n}\mu_i^{2N}\to0 \quad (i\to\infty,\ \rho\text{ fixed}),
\]
both $\zeta_i$ and its induced correction contribute $o_i(\mu_i^{n-2N})$ to the localized Pohozaev integral.

We apply Lemma \ref{lem:ac-decomposition} with the improved annular estimate \eqref{eq:ac-display-7-3}.
For fixed $\rho\in(0,\rho_0)$ and all sufficiently large $i\in\mathbb N$, $A_i\le (C+1)\mu_i^{n-2N}$, with $C$ independent of $\rho$.
Equations \eqref{eq:ac-cutoff-linear-error} and \eqref{eq:ac-cutoff-quadratic-error} bound the contribution linear in $b$ by $C\rho^2A_i$ and the contribution quadratic in $b$ by $C\rho^2A_i^2R_i^{n-2N}$.
These two terms satisfy
\begin{equation}
\mu_i^{2N-n}\{\rho^2A_i+\rho^2A_i^2R_i^{n-2N}\}
\le C(\rho^2+\rho^{n-2N+2}).       \label{eq:ac-display-8-3}
\end{equation}

To estimate the change in the quadratic term for $U+\psi$ caused by the cutoff functions, we apply the conclusion of Lemma \ref{lem:ac-quadratic} for fixed $\rho$.
Equations \eqref{eq:ac-display-5-2} and \eqref{eq:ac-display-5-3} imply $J_i=O(\mu_i^{(n-2N)/2})$ and $H_{i,d}\to0$ as $i\to\infty$ for $2\le d\le (n-2N)/2$.
Thus, with $t_i$ as in \eqref{eq:subcritical-cutoff-tail}, $C(J_it_i+t_i^2)=o_i(\mu_i^{n-2N})$ at fixed $\rho$.
At degree $(n-2N)/2$, when this degree is an integer at least two, the bounded term $O(\mu_i^{n-2N}|H_{i,(n-2N)/2}|^2)$ after subtraction of the logarithmic term is also $o_i(\mu_i^{n-2N})$.
Lemma \ref{lem:ac-quadratic} bounds the sum of the pairings containing a Taylor term of degree greater than $(n-2N)/2$ by $\theta J_i^2+C_\theta\mu_i^{n-2N}\rho^{2d_0+2N-n}$.
By \eqref{eq:ac-metric-rem}, the part of the approximation error generated by the metric perturbation contributes at most $C\rho^2J_i^2$,
the same bound for pairings containing terms of degree greater than $(n-2N)/2$, and
\begin{align*}
\theta J_i^2+C_\theta\mu_i^{n-2N}\rho^{2T+2N-n}.
\end{align*}
The cubic metric terms satisfy the same bounds, and the direct Taylor operator error is $C\mu_i^{n-2N}\rho^{T+2N-n}$.
Take $\theta$ small and then reduce the fixed upper bound on $\rho$ so that the positive $J_i^2$ error terms can be absorbed into the negative quadratic term from Lemma \ref{lem:ac-quadratic}.
Since $2T+2N-n>T+2N-n>2\ge\eta$, the two additional Taylor errors are bounded by $C\mu_i^{n-2N}\rho^\eta$, uniformly for small $\rho$.
Together with \eqref{eq:ac-display-8-3}, this proves \eqref{eq:ac-display-8-1}.
First divide by $\mu_i^{n-2N}$ and let $i\to\infty$ at fixed $\rho$; then let $\rho\downarrow0$ to obtain \eqref{eq:ac-display-8-2}.
\end{proof}

\begin{proof}[Proof of Theorem \ref{thm:ac-main}]
Suppose that $u_i>0$ solves \eqref{eq:constant-q-equation} with $\kappa=\kappa_0$ and $\|u_i\|_{L^\infty(M)}\to\infty$.
At the pole $p_1$ selected in Proposition \ref{prop:ac-green-limit}, we use $\mu_i=\mu_{1,i}$ and its fixed limiting conformal-normal gauge.
The finite-part identity \eqref{eq:repair-green-sum-finite-part} identifies the constant term in \eqref{eq:ac-display-7-4} as
\begin{equation}
B=\frac{A_0}{\ga_{n,N}}\left(
m_g(p_1)+\sum_{j\ne1}\tau_jG_g(p_j,p_1)\right)>0.
\label{eq:ac-display-8-4}
\end{equation}
Formula \eqref{eq:ac-display-8-4} uses \eqref{eq:mass-hypothesis} at $p_1\in\mcz$ and positivity of $G_g$.
Theorem \ref{thm:repair-Green} and Lemma \ref{lem:repair-angular-mass} imply
\begin{equation}
\lim_{\rho\downarrow0}\mcp_N^0(\rho,\Lambda_{p_1}^{-1}v)
=\frac{n-2N}{2\ga_{n,N}}A_0B>0.                \label{eq:ac-display-8-5}
\end{equation}
Here $u_i$ in the Pohozaev integral is the local representative in the moving conformal-normal chart, as specified in Section \ref{sec:ac-tools}.
For every fixed $\rho\in(0,\rho_0)$, the annular convergence \eqref{eq:ac-display-6-1} implies $\mu_i^{2N-n}\mcp_N^{\kappa_0}(t\rho,u_i)\to\mcp_N^0(t\rho,\Lambda_{p_1}^{-1}v)$ as $i\to\infty$, uniformly for $t\in\supp\chi'$.
The nonlinear integral over the coordinate sphere vanishes with a factor $\mu_i^{2N}$.
The averaging weight in \eqref{eq:ac-display-4-4} has integral one, so the Pohozaev integral and its radial average have the same double limit \eqref{eq:ac-display-8-5}.
The positive limit in \eqref{eq:ac-display-8-5} contradicts the upper bound \eqref{eq:ac-display-8-2}.

The resulting $L^\infty$ bound and elliptic estimates make the positive solution set precompact in $C^0$.
The bounded inverse gives $\|u\|_\infty\le C\kappa_0\|u\|_\infty^{(n+2N)/(n-2N)}$, so a uniform limit cannot vanish identically.
Positivity of the Green representation makes every such limit strictly positive, yielding a uniform lower bound.
Elliptic bootstrapping then proves compactness in every $C^m$ topology with $m\in\mathbb N\cup\{0\}$.
\end{proof}

\appendix
\section{Reconstruction of TT tensors and normalization}\label{app:weighted}
\subsection{Reconstruction and pairings with a first spherical harmonic} \label{app:weighted-reconstruction}
We supply the reconstruction and normalization identities used in the angular reduction of Section \ref{sec:ambient-york}.
For each fixed angular harmonic, the trace and divergence equations determine the corresponding TT tensor from one scalar radial function.
We compute the action of $\mca=\nabla^*\nabla+2$ from \eqref{eq:round-tensor-operator}, together with scalar formulas for the tensor norm and for the tensor pairing after multiplication by the signed height $z$.
Let $p_-$ be the image of the Euclidean origin under stereographic compactification, and let $\vrh$ denote round distance from $p_-$.
We write
\[
g_c=d\vrh^2+S^2\ga ,\quad\text{where }  S:=\sin\vrh ,\quad  z=\cos\vrh ,\quad\text{and}\quad  t:=1-z^2,
\]
where $\ga$ is the unit metric on $\Ss^{n-1}$.
In this appendix, we use $\ga$ for derivatives and norms of angular fields on $\Ss^{n-1}$.
For an angular degree $\ell\ge2$, we write
\[
\kappa_\ell=\ell(\ell+n-2) ,\quad 
\mu_\ell=(\ell-1)(\ell+n-1) ,\quad\text{and}\quad 
\beta_\ell=\frac{n-2}{n-1}(\kappa_\ell-n+1).
\]
In each nonzero angular space, we choose a harmonic of unit $L^2(\ga)$ norm: $\msM$ is a TT tensor with eigenvalue $\kappa_\ell-2$ for the nonnegative rough Laplacian,
$v$ is a divergence-free vector with eigenvalue $\kappa_\ell-1$ for that Laplacian, and $-\Delta_\ga Y=\kappa_\ell Y$.
They are the angular restrictions of the harmonic tensors in Proposition \ref{thm:finite-SVT-exhaustion-new}.
We define
\[
\msV(v)_{ab}:=\nabla_av_b+\nabla_bv_a \quad\text{and}\quad 
\msS(Y)_{ab}:=\nabla_a\nabla_bY+\frac{\kappa_\ell}{n-1}Y\ga_{ab}.
\]
With the divergence convention $(\divg_\ga h)_b:=\nabla^ah_{ab}$, their divergence identities and norms are
\[
\divg_\ga\msV(v)=-\mu_\ell v ,\quad 
\|\msV(v)\|_2^2=2\mu_\ell ,\quad 
\divg_\ga\msS(Y)=-\beta_\ell\,dY ,\quad\text{and}\quad 
\|\msS(Y)\|_2^2=\kappa_\ell\beta_\ell.
\]
For the scalar reconstruction functions, we use the weight \eqref{alg:eq:cyclic} and operator \eqref{eq:radial-operator}:
\[
w_\ell=t^{\frac{n}{2}+\ell-1} \quad\text{and}\quad 
A_\ell=-t\pp_z^2+(n+2\ell)z\pp_z+\ell(\ell+n-1).
\]
We first prove the identities for scalar reconstruction functions in $C_c^\infty(-1,1)$.
They extend by approximation whenever the pairings converge and the endpoint terms at $z=\pm1$ tend to zero.

\begin{proposition}
\label{prop:weighted-reconstruction}
We define
\[
P_V:=-t\pp_z+(n+\ell-1)z ,\quad 
P:=-t\pp_z+(n+\ell-2)z ,\quad\text{and}\quad 
P_S:=\frac{P^2}{\kappa_\ell}-\frac{1}{n-1}.
\]
The nonzero covariant components of the three reconstructions are
\begin{equation}\label{eq:reconstruction-ansatz}
\begin{aligned}
(\mcr_Tu)_{ab}&=S^{\ell+2}u\msM_{ab},\\
(\mcr_Vu)_{\vrh a}&=S^\ell uv_a\quad\text{and}\quad
(\mcr_Vu)_{ab}=\frac{S^{\ell+1}}{\mu_\ell}(P_Vu)\msV(v)_{ab},\\
(\mcr_Su)_{\vrh\vrh}&=S^{\ell-2}uY\quad\text{and}\quad
(\mcr_Su)_{\vrh a}=\frac{S^{\ell-1}}{\kappa_\ell}(Pu)\nabla_aY,\\
(\mcr_Su)_{ab}&=S^\ell\left(-\frac{u}{n-1}Y\ga_{ab}
+\frac{P_Su}{\beta_\ell}\msS(Y)_{ab}\right).
\end{aligned}
\end{equation}
These tensors are TT.
Taking the adjoint of $\mcr_\sigma$ with respect to the spherical $L^2$ tensor pairing and the scalar pairing in $L^2(w_\ell\,dz)$, we obtain
\begin{equation}\label{eq:reconstruction-intertwining-norm}
\mca\mcr_\sigma=\mcr_\sigma A_\ell,\qquad
\mcr_\sigma^*\mcr_\sigma
=c_{\sigma,\ell}W_\sigma(A_\ell),\qquad
c_{T,\ell}:=1,\quad c_{V,\ell}:=\frac{2}{\mu_\ell},\quad
c_{S,\ell}:=\frac{1}{\kappa_\ell\beta_\ell},
\end{equation}
where $W_T(x)=1$, $W_V(x)=x$, and $W_S(x)=x(x+n-2)$.
\end{proposition}
\begin{proof}
The tensor-derived reconstruction is TT because $\msM$ has both vanishing trace and vanishing divergence on $\Ss^{n-1}$.
In the polar-coordinate calculations below, primes denote derivatives with respect to $\vrh$.
For the vector ansatz $h_{\vrh a}=bv_a$, $h_{ab}=S^2d\msV(v)_{ab}$, the remaining divergence equation is $\mu_\ell d=b'+(n-1)\cot\vrh\,b$.
Substituting $b=S^\ell u$ yields the coefficient of $\msV(v)$ in \eqref{eq:reconstruction-ansatz}.
For the scalar ansatz
\[
h_{\vrh\vrh}=aY ,\quad  h_{\vrh a}=b\nabla_aY ,\quad\text{and}\quad 
h_{ab}=S^2(cY\ga_{ab}+d\msS(Y)_{ab}),
\]
the trace and divergence equations are
\[
c=-\frac{a}{n-1} ,\quad 
\frac{\kappa_\ell b}{S^2}=a'+n\cot\vrh\,a ,\quad\text{and}\quad 
\beta_\ell d=b'+(n-1)\cot\vrh\,b+c.
\]
With $a=S^{\ell-2}u$, these equations yield the scalar reconstruction in \eqref{eq:reconstruction-ansatz}.

Computing $\mca=\nabla^*\nabla+2$ in polar coordinates and using the trace and divergence equations yields
\begin{align*}
(\mca h)_{\vrh a}
&=\left(-b''-(n-1)\cot\vrh\,b'+\frac{\kappa_\ell}{S^2}b\right)v_a
&&\text{in the vector branch},\\
(\mca h)_{\vrh\vrh}
&=\left(-a''-(n+3)\cot\vrh\,a'
+\left(\frac{\kappa_\ell}{S^2}-2n\cot^2\vrh+2\right)a\right)Y
&&\text{in the scalar branch}.
\end{align*}
Substituting $b=S^\ell u$ and $a=S^{\ell-2}u$ yields $S^\ell(A_\ell u)v_a$ and $S^{\ell-2}(A_\ell u)Y$, respectively.
For the tensor reconstruction, writing $h_{ab}=S^2f\msM_{ab}$ yields
\[
(\mca h)_{ab}
=S^2\left(-f''-(n-1)\cot\vrh\,f'+\frac{\kappa_\ell}{S^2}f\right)\msM_{ab} ,\quad\text{where }  f:=S^\ell u.
\]
This equals $S^{\ell+2}(A_\ell u)\msM_{ab}$.
The trace and divergence commutators on the round sphere show that $\mca$ preserves TT tensors; its polar-coordinate formula also preserves each angular harmonic type.
For the vector and scalar reconstructions, the displayed component determines the other coefficients uniquely through the trace and divergence equations.
This proves the intertwining identity for every component of the reconstructed tensor.

Integrating the squared angular components gives
\begin{align*}
\|\mcr_Tu\|_2^2&=\int_{-1}^{1} w_\ell u^2,\\
\|\mcr_Vu\|_2^2&=\frac{2}{\mu_\ell}
\int_{-1}^{1} t^{\frac{n}{2}+\ell-2}\bigl((P_Vu)^2+\mu_\ell u^2\bigr),\\
\|\mcr_Su\|_2^2&=
\int_{-1}^{1} t^{\frac{n}{2}+\ell-3}
\left(\frac{n}{n-1}u^2+\frac{2}{\kappa_\ell}(Pu)^2
+\frac{\kappa_\ell}{\beta_\ell}(P_Su)^2\right).
\end{align*}
The adjoint of $P_V$ in $L^2(t^{n/2+\ell-2}\,dz)$ is $P_V^\dagger=t\pp_z-(\ell-1)z$, and
\[
t^{-1}(P_V^\dagger P_V+\mu_\ell)=A_\ell.
\]
The adjoint of $P$ in $L^2(t^{n/2+\ell-3}\,dz)$ is
\begin{equation}\label{eq:reconstruction-adjoint-operators}
P^\dagger:=t\pp_z-(\ell-2)z.
\end{equation}
For a multiplication operator $F$, we define
\[
\mcu_\ell(F):=(P^\dagger)^2FP^2-\frac{\kappa_\ell}{n-1}((P^\dagger)^2F+FP^2)
+2\beta_\ell P^\dagger FP
+\left(\frac{\kappa_\ell^2}{(n-1)^2}
+\frac{n\kappa_\ell\beta_\ell}{n-1}\right)F.
\]
Polarizing the scalar norm and integrating by parts gives
\[
c_{S,\ell}^{-1}\mcr_S^*\mcr_S
=t^{-2}\mcu_\ell(I)
=A_\ell(A_\ell+n-2).
\]
Substituting $P^\dagger$ and $P$ and collecting the coefficients of each derivative proves the last identity and hence \eqref{eq:reconstruction-intertwining-norm}.
\end{proof}

We next compute $\mcr_\sigma^*z\mcr_\sigma$, corresponding to multiplication by $z$ in the tensor inner product.

\begin{proposition}
\label{prop:weighted-height}
For an operator $T$ on scalar reconstruction functions, we define $\mcj_\ell(T):=A_\ell T+TA_\ell$.
In \eqref{eq:reconstruction-height}, $z$ denotes multiplication by the height function, and a polynomial in $\mcj_\ell$ acts on this multiplication operator.
For every angular degree $\ell\ge2$, multiplication by $z$ in the tensor pairing satisfies
\begin{equation}\label{eq:reconstruction-height}
\begin{gathered}
\mcr_\sigma^*z\mcr_\sigma
=c_{\sigma,\ell}\wtbeta_\sigma(\mcj_\ell)z,\\
\wtbeta_T(u):=1 ,\quad 
\wtbeta_V(u):=\frac{u+n-2}{2} ,\quad\text{and}\quad 
\wtbeta_S(u):=\frac{(u+n-2)(u+3n-8)}{4}.
\end{gathered}
\end{equation}
\end{proposition}
\begin{proof}
Insert multiplication by $z$ in the polarized norm formulas.
The tensor identity is immediate, and the vector identity follows from
\[
t^{-1}(P_V^\dagger zP_V+\mu_\ell z)
=\frac{1}{2}(A_\ell z+zA_\ell+(n-2)z).
\]
For the scalar identity the same calculation gives
\[
c_{S,\ell}^{-1}\mcr_S^*z\mcr_S
=t^{-2}\mcu_\ell(z)
=\frac{1}{4}(\mcj_\ell+n-2)(\mcj_\ell+3n-8)z.
\]
Both sides are differential operators of order at most four; equality on $1,z,\ldots,z^4$ therefore proves the identity, with $n,\ell$ left symbolic.
\end{proof}

\medskip
We now apply the reconstruction formulas to the normalized first moment in \eqref{eq:height-curvature-functional}.
In the coordinates used here, the fixed first spherical harmonic is $Y_{p_+}=-z$.
Conformal covariance and self-adjointness give
\begin{equation}\label{eq:completion-moment-variation}
D\mscF(g)[2Y_{p_+}g]=(n-2N)\mcm_{p_+}(g).
\end{equation}
For $\dot g=2\phi g$, the total curvature $\mcq(g)=\int_{\Ss^n}Q_{2N,g}\,dv_g$ has derivative $(n-2N)\allowbreak\int_{\Ss^n}\phi Q_{2N,g}\,dv_g$,
whereas $D\Vol[2\phi g]=n\int_{\Ss^n}\phi\,dv_g$.
Differentiating the volume factor subtracts the average curvature in the definition of $\mcm_{p_+}$.

For operators on TT tensors, multiplication by $z$ means $h\mapsto\Pi(zh)$.
We set
\[
\mcj_{\mfull}(T):=\mca T+T\mca ,\quad 
M_{\mfull}:=\frac{1}{2}(\mca z+z\mca-nz) ,\quad\text{and}\quad 
M_\ell:=\frac{1}{2}(\mcj_\ell-n)z.
\]
Let $p_{N,n}$ be the Hessian polynomial in \eqref{eq:ambient-hessian-coefficient}.
The Hessian of the normalized first moment \eqref{eq:height-curvature-functional} is represented on TT tensors by
\[
D^2\mcm_{p_+}(g_c)[h_1,h_2]
=|\Ss^n|^{-\frac{n-2N}{n}}
\langle h_1,\mcn_{N,n}h_2\rangle_{L^2(g_c)},
\]
where
\[
\mcn_{N,n}:=\frac{G^*p_{N,n}(\mca)+p_{N,n}(\mca)G}{n-2N} ,\quad\text{where } 
G:=\frac{1}{2}[\mca,z]-\frac{n}{2}z.
\]
Indeed, diffeomorphism invariance replaces the conformal direction $-2zg$ in \eqref{eq:completion-moment-variation} by $(\mcl_X-2z)g$, where $X=-\nabla z$ as in Proposition \ref{prop:repair-dilation}.
This direction and $D\mscF$ vanish at $g_c$, so $\mcm_{p_+}(g_c)=0$ and $D\mcm_{p_+}(g_c)=0$.
The derivative of $g\mapsto(\mcl_X-2z)g$, restricted to TT tensors and followed by $\Pi$, is $G$.
Two derivatives of \eqref{eq:completion-moment-variation} therefore leave exactly the two Hessian terms in the numerator; stationarity makes this second derivative independent of the metric chart.
We define the polynomial $\mathfrak r_{N,n}$ by
\begin{align*}
(n-2N)(u-n)\mathfrak r_{N,n}(u)
&=(y-x)(p_{N,n}(y)-p_{N,n}(x))
-n(p_{N,n}(x)+p_{N,n}(y)),\\
u&=x+y \quad\text{and}\quad  (y-x)^2=2u+n(n-2).
\end{align*}
The right-hand side is symmetric in $x,y$ and hence is a polynomial in $u$ under the displayed relation.
At $u=n$, that relation implies that $x,y$ are $0,n$ in either order; since $p_{N,n}(0)=0$, the right-hand side vanishes there.
Thus $\mathfrak r_{N,n}$ has degree at most $N-1$.
The $A_\ell$ eigenfunctions are Gegenbauer polynomials in $z$, indexed by spherical degree $k=\ell,\ell+1,\ldots$.
Multiplication by $z$ has nonzero matrix entries only between adjacent degrees; \eqref{eq:reconstruction-height} proves the same property after TT projection.
Between eigenfunctions with eigenvalues $x=k(k+n-1)$ and $y=(k+1)(k+n)$, the matrix entry of $\mcn_{N,n}$ is the corresponding matrix entry of multiplication by $z$ times
\[
\frac{(y-x)(p_{N,n}(y)-p_{N,n}(x))-n(p_{N,n}(x)+p_{N,n}(y))}{2(n-2N)}
=\frac{u-n}{2}\mathfrak r_{N,n}(u).
\]
Consequently
\[
\mcn_{N,n}=\mathfrak r_{N,n}(\mcj_{\mfull})M_{\mfull}.
\]
Equations \eqref{eq:reconstruction-intertwining-norm} and \eqref{eq:reconstruction-height} therefore express this operator on scalar reconstruction functions as
\begin{equation}\label{eq:reconstruction-first-moment}
\mcr_\sigma^*\mcn_{N,n}\mcr_\sigma
=c_{\sigma,\ell}\mathfrak r_{N,n}(\mcj_\ell)
\wtbeta_\sigma(\mcj_\ell)M_\ell.
\end{equation}

\subsection{Gauge terms and normalization of the Pohozaev integral}\label{app:low-order}
The original polynomial tensors satisfy the normal-coordinate conditions; their projections under $\Pi$ satisfy the TT conditions.
The \emph{first jet} of a tensor at a point consists of its value and covariant derivative there.
Although the lowest vector-derived and scalar-derived homogeneous tensors vanish to orders three and four, their compactified TT projections can have nonzero first jets at $p_-$.

\begin{proposition}
At angular degree $\ell=2$, consider the vector and scalar TT tensors reconstructed from $u=\Phi_0=(1+z)^{-2}$.
The vector reconstruction has a nonzero first derivative at $p_-$, and the scalar reconstruction has a nonzero value there.
\end{proposition}
\begin{proof}
Use polar distance $\vrh$ from $p_-$, so $z=\cos\vrh$ and $\Phi_0=1/4+O(\vrh^2)$.
The reconstruction formulas \eqref{eq:reconstruction-ansatz} extend smoothly across $p_-$: in normal coordinates their leading terms are multiples of $\pp_iV_j+\pp_jV_i$ and $\pp_i\pp_jY$, respectively, and the remaining radial coefficients are smooth functions of $\vrh^2$.
For the vector reconstruction the covariant mixed coefficient is $b=\sin^2\vrh\,\Phi_0$.
In an orthonormal angular frame the mixed component is $(\vrh/4+O(\vrh^3))v_a$, where $v$ is the nonzero degree-two vector harmonic.
For the scalar reconstruction the meridional component is $h_{\vrh\vrh}=\Phi_0Y=(1/4+O(\vrh^2))Y$.
\end{proof}

We determine the coefficient relating the global Pohozaev integral to the Hessian of $\mcm_{p_+}$.

\begin{proposition}
\label{prop:completion-global-TT}
Let $H_1,H_2$ be smooth trace-free Euclidean metric directions whose compactifications $h_i:=\Omega^2H_i$ are smooth TT tensors with respect to $g_c$, using \eqref{eq:main-notation1}.
In this proposition and its proof, we assume that the Euclidean integrals converge, or define them by smooth spherical compactification.
Set $g_{\ep,t}:=\exp(\ep H_1+tH_2)$.
For a twice differentiable positive family $u_{\ep,t}=U+O(\ep^2+|\ep t|+t^2)$, the following quotient of the Pohozaev integral has mixed derivative
\begin{equation}\label{eq:completion-global-coefficient}
\left.\pp_\ep\pp_t\right|_0
\frac{\int_{\R^n} (\msZ u_{\ep,t})
\left(P_{2N,g_{\ep,t}}u_{\ep,t}-\kappa_{\ep,t}u_{\ep,t}^{\frac{n+2N}{n-2N}}\right)\,dx}
{\left(\int_{\R^n} u_{\ep,t}^{\frac{2n}{n-2N}}\,dx\right)^{\frac{n-2N}{n}}}
=-\frac{(n-2N)^2}{4}\,D^2\mcm_{p_+}(g_c)[h_1,h_2].
\end{equation}
The coefficient is independent of the mixed second scalar variation and of the spatially constant multiplier $\kappa_{\ep,t}$.
In these TT directions, the solution of the projected linearized equation has zero first scalar variation when its projection onto the span of the bubble Jacobi fields is fixed to zero.
\end{proposition}

\begin{proof}
We first compute the coefficient with the bubble fixed, then check that a second-order scalar change does not affect it.
For every sufficiently decaying positive $u$, integration by parts and the definition of $\msZ$ yield
\[
\int_{\R^n}(\msZ u)u^{\frac{n+2N}{n-2N}}\,dx=0.
\]
Thus the multiplier contributes zero to the Pohozaev integral.
For a fixed determinant-one metric $g$, the energy \eqref{eq:normalized-energy-quotient} uses $dv_g=dx$.
Differentiating in the scalar argument with $g$ fixed, and using conformal covariance, gives
\[
\frac{1}{2} D_u\mce_N(g,u)[\msZ u]
=\frac{\int_{\R^n} (\msZ u)P_{2N,g}u\,dx}{\left(\int_{\R^n} u^{\frac{2n}{n-2N}}\,dx\right)^{\frac{n-2N}{n}}} \quad\text{and}\quad 
\mce_N(g,u)=\frac{n-2N}{2}\mscF\left(u^{\frac{4}{n-2N}}g\right).
\]
At $u=U$, the metric variation induced by $\msZ u$ is $4(\msZ U/U)\Omega^2g/(n-2N)=-2Y_{p_+}\Omega^2g$.
Equation \eqref{eq:completion-moment-variation} proves the exact identity
\[
\frac{\int_{\R^n} (\msZ U)P_{2N,g}U\,dx}{\left(\int_{\R^n} U^{\frac{2n}{n-2N}}\,dx\right)^{\frac{n-2N}{n}}}
=-\frac{(n-2N)^2}{4}\mcm_{p_+}(\Omega^2g).
\]
Taking the mixed derivative in the exponential family $g_{\ep,t}$ proves \eqref{eq:completion-global-coefficient} when the scalar function is fixed at $U$.

It remains to check that an $O(\ep^2+|\ep t|+t^2)$ scalar change has no effect.
The standard bubble $U$ satisfies $P_0U-\kappa_0U^{(n+2N)/(n-2N)}=0$, whose linearization in the scalar variable is $L_U$ in \eqref{eq:bubble-linearization}.
The scalar derivative of the quotient in \eqref{eq:completion-global-coefficient} in a direction $\phi$ equals $(\int_{\R^n} U^{2n/(n-2N)})^{-(n-2N)/n}\allowbreak\int_{\R^n} (\msZ U)L_U\phi\,dx=0$, by self-adjointness and $L_U(\msZ U)=0$.
This proves independence of the mixed scalar variation.
The vanishing of the first scalar correction term is \cite[Proposition 4.12]{GKW26}, with its projection onto the span of the Jacobi fields fixed to zero.
\end{proof}

\begin{remark}
For polynomial directions of degree less than $(n-2N)/2$, Proposition \ref{thm:strict-york-transmission-new}
and Corollary \ref{cor:strict-hessian-transmission-new} prove convergence of the Hessian pairings under cutoff, including their derivatives under $H(x)\mapsto H(\lambda x)$ at $\lambda=1$.
Propositions \ref{prop:scalar-corrector-unique-new} and \ref{thm:schur-york-identity-new} identify the limits of the scalar correction term and the corresponding second variation of $\mce_N$.
Proposition \ref{prop:repair-dilation} then identifies the quadratic coefficient of the Pohozaev integral with that in Proposition \ref{prop:completion-global-TT}, even when a Euclidean radial cutoff does not preserve the TT conditions.
\end{remark}

To compare with the energy and scale-derivative matrices in \cite{GKW26}, we finally translate this coefficient into the fixed-volume exponential convention.
The notation $D^2_{\exp}\int_{\Ss^n} Q_{2N}$ in the next proposition means the mixed second derivative in the determinant-one exponential chart.
The functional $\mcq$ is not stationary under arbitrary metric variations at the round metric $g_c$, so its second derivative depends on the chosen chart.
In the fixed-volume chart used here, this derivative is $D^2\mscF(g_c)$ multiplied by $|\Ss^n|^{(n-2N)/n}$.

\begin{proposition}
\label{prop:completion-Weyl-calibration}
Let $H_i$ be radial powers times the quadratic tensors associated with algebraic Weyl tensors in Proposition \ref{thm:finite-SVT-exhaustion-new}, centered at the Euclidean origin.
Assume that their homogeneous degrees satisfy $d_i<(n-2N)/2$, and set
\[
E_{ij}:=D^2_{\exp}\!\left(\int_{\Ss^n} Q_{2N}\,dv\right)_{g_c}
[\Omega^2H_i,\Omega^2H_j] \quad\text{and}\quad 
S_{ij}:=(d_i+d_j)E_{ij}.
\]
The matrix of the mixed second derivative of the normalized Pohozaev integral in \eqref{eq:completion-global-coefficient} equals
\begin{equation}\label{eq:completion-S-calibration}
-\frac{n-2N}{4|\Ss^n|^{\frac{n-2N}{n}}}S.
\end{equation}
Thus a negative definite $S$ gives a positive definite global Pohozaev quadratic form.
The local Pohozaev integral in \eqref{eq:ac-display-4-3} uses $P_0-P_{2N,g}$ and has the opposite leading sign.
\end{proposition}
\begin{proof}
The reduced energy in \cite[Proposition 4.12]{GKW26} is the value, after solving the projected scalar equation, of
\[
\mathcal J(g,u):=\frac12\int_{\R^n}uP_{2N,g}u\,dv_g
-\frac{(n-2N)\kappa_0}{2n}\int_{\R^n}u^{\frac{2n}{n-2N}}\,dv_g.
\]
There is no division by a volume factor in $\mathcal J$.
Its quadratic coefficient is $(n-2N)E_{ij}/8$, so its mixed second derivative is $(n-2N)E_{ij}/4$.
Since $H_i(\lambda x)=\lambda^{d_i}H_i(x)$, differentiation in the bubble scale at $\lambda=1$ multiplies this entry by $d_i+d_j$; see \cite[Proposition 5.3]{GKW26}.
For $U_\lambda(x):=\lambda^{-(n-2N)/2}U(x/\lambda)$, one has $\pp_{\log\lambda}U_\lambda|_{\lambda=1}=-\msZ U$.
The scale derivative of this reduced energy is therefore the negative of the numerator in \eqref{eq:completion-global-coefficient}.
Division by $|\Ss^n|^{(n-2N)/n}$ proves \eqref{eq:completion-S-calibration}.
\end{proof}

\section{Rational scalar pairings and exact verification}\label{app:verification}
\subsection{Rational formulas for the scalar pairings}\label{app:rational-moments} 
We use the scalar extraction in \cite[Lemma 6.2]{GKW26} and the finite-summation method in \cite[Theorem 6.5]{GKW26} to prove Proposition \ref{alg:prop:moments}.
The cited calculation treats the polynomial metric perturbations generated by algebraic Weyl tensors of homogeneous degrees $2,4,6,8$.
The nine scalar pairings $m_0,\ldots,m_8$ determine all $218$ rational compactness matrices and their $430$ pivots; each block uses only $m_0,\ldots,m_K$, with $K$ specified in \eqref{alg:eq:momentbound}.
They also determine the matrices used to prove Corollary \ref{cor:first-pohozaev-failure} and Theorem \ref{thm:three-sector-crossing}; all the linear systems below have size at most ten.
Appendix \ref{app:ancillary} describes the common matrix and pivot construction and the verification program.

Throughout this subsection, we use $t=r^2$, the radial operator $A$ in \eqref{alg:eq:triA}, and the measure, normalization, and scalar pairings in \eqref{alg:eq:normalized-pairing}--\eqref{alg:eq:scalar-pairings}.
We first treat $\msv$ as an indeterminate in the exact rational calculations and then evaluate it at $n/2$.
Fix $0\le j\le8$ and $n>2N+2j+4$.
We define the polynomial
\begin{equation}\label{alg:eq:ambientpoly}
F_j(w,t):=w^2\sum_{a=0}^{j}
\frac{(-1)^a(-j)_a(j+2)_a}{(1-\msv)_a a!}
[1+(t-1)w]^{j-a}[(w-1)(1+tR)]^a.
\end{equation}
We set
\[
w:=\frac{1+\frac{\rho}{2}}{1-\frac{\rho t}{2}} ,\quad 
\upsilon:=\frac{2\rho}{\left(1+\frac{\rho}{2}\right)^2} ,\quad\text{and}\quad  c:=\sqrt{1-\upsilon}.
\]
Here $\rho$ is the ambient coordinate.
The polynomial identity to verify is
\begin{align}\label{alg:eq:ambientidentity}
0={}&w^2(w-1)F_{ww}-2Rt(w-1)F_{wt}
-\frac{1}{2}w(2w+2\msv-4)F_w\notag\\
&-t(1+t)F_{tt}+\frac{1}{2}((2\msv-8)t-2\msv-4)F_t+(2\msv-2)F.
\end{align}
Substitution of $F=F_j$ verifies this identity and the initial condition $F_j(1,t)=t^j$.
This is \cite[(6.5)]{GKW26}, obtained from the ambient equation \cite[(3.2)]{GKW26} with $b=-A$.
The cancellation using successive summand coefficients in \cite[Proposition 6.1]{GKW26} remains valid for $j=4,\ldots,8$; clearing $(1-\msv)_j$ verifies the resulting polynomial identity for each of these values.

Since $F_0=w^2$, the integral arising in the mixed second variation of $\log\det g_\rho$ simplifies to
\[
\frac{1+t}{2}\int_1^w u^2(F_0)_u(F_j)_u\,du
=\sum_{L=5}^{2j+5}\kappa_{j,L}(t)(w^L-1) ,\quad\text{where } 
\kappa_{j,L}:=\frac{(L-3)(1+t)}{L}[w^{L-3}]F_j.
\]
Here $[w^{L-3}]$ extracts the indicated polynomial coefficient.
We write $\langle f\rangle_\mu:=\int_0^\infty f\,d\mu$ for the corresponding scalar average.
The polarized determinant equation \cite[(3.6)]{GKW26} therefore yields
\begin{equation}\label{alg:eq:integrated-determinant-factor}
\mscH_j(\upsilon)
:=-\frac{\langle1,At^j\rangle_\mu}{2(n-1)}
-\sum_{L=5}^{2j+5}\langle\kappa_{j,L}(w^L-1)\rangle_\mu.
\end{equation}
The quantity $\mscH_j$ is $(1+\rho/2)^2$ times the $\rho$ derivative of the mixed second variation of $\log\det g_\rho$, integrated against $d\mu$ and written as a function of $\upsilon$ by inverting $\upsilon=2\rho/(1+\rho/2)^2$ near zero.

Take $\rho<0$ near zero, so $c>1$, and use the symmetric measure and Cauchy transform of \cite[Lemma A.1]{GKW26}, with $n=2\msv$ and the integration variable replaced by its negative:
\[
d\nu_0(\zeta):=\frac{(1-\zeta^2)^{\msv-3}}{B\left(\frac{1}{2},\msv-2\right)}\,d\zeta \quad\text{and}\quad  \mfF(c):=\int_{-1}^1\frac{d\nu_0(\zeta)}{c+\zeta} \quad\text{for }  c>1.
\]
Polynomial division and differentiation as in \cite[(A.5)]{GKW26} give
\begin{equation}\label{alg:eq:Hdecomp}
\mscH_j(\upsilon)
=A_j(\upsilon)\mfF(\sqrt{1-\upsilon})
+\sqrt{1-\upsilon}\,B_j(\upsilon)+c_j(\msv).
\end{equation}
For $0\le j\le8$, the Laurent powers of $A_j$ and $B_j$ lie in $[-j,2]$ and $[-j,1]$, respectively, with coefficients in $\Q(\msv)$.
Appendix \ref{app:ancillary} specifies the finite checks of these bounds.
For an indeterminate $k$, we define $r_j(k,\msv)$ by \cite[(A.11)]{GKW26}, using the coefficients of $A_j,B_j$ in the expression for $R_{ij}$ in \cite[(A.8)]{GKW26}.
These formulas use only the Laurent coefficients and are unchanged when the profile index exceeds three.

The inequality $n>2N+2j+4$ implies $N+j-1<\msv-3$.
Thus the Cauchy expansion in \cite[(A.3)]{GKW26} supplies all coefficients through order $N+j-1$; its omitted part, multiplied by $A_j$, contributes nothing through order $N-1$.
The scalar coefficient formula \cite[(6.11)]{GKW26} is independent of the radial profile and applies because \eqref{alg:eq:ambientidentity} has been verified.
Together with the coefficient ratio \cite[(A.10)]{GKW26}, it yields
\begin{equation}\label{alg:eq:sum}
m_j=c_j(\msv)+\sum_{k=0}^{N-1}
\frac{(1-N)_k(1+N)_k}{(2-\msv)_k(1+\msv)_k}r_j(k,\msv).
\end{equation}
The finite symbolic calculation verifies that
\begin{equation}\label{alg:eq:rhs}
q_j(k,\msv):=(k+2-\msv)_{j+2}r_j(k,\msv)
\quad\text{is a polynomial in $k$ of degree at most $j+1$}.
\end{equation}
The computation first constructs a polynomial candidate of degree at most $j+1$ from the $j+2$ values at $k=\msv-i-1$, $1\le i\le j+2$,
evaluating the right-hand side of \eqref{alg:eq:rhs} after cancellation of its removable singularities.
The difference between the candidate and that rational expression has a common denominator dividing
\[
(2k+1)\prod_{m=1}^j(k+m)
\prod_{a=-j}^{2}(k-a-\frac{1}{2})
\prod_{a=-j}^{2}(\msv-3-k+a).
\]
This product has degree $3j+7$ in $k$.
Since $r_j(k,\msv)=O(k^{-1})$, multiplying the difference by this common denominator produces a polynomial of degree at most $4j+8$.
Its exact vanishing at $k=0,\ldots,4j+8$, as identities over $\Q(\msv)$, therefore proves \eqref{alg:eq:rhs}.

\begin{proof}[Proof of Proposition \ref{alg:prop:moments}]
For $0\le j\le8$, we define the difference operator from \cite[(6.18)]{GKW26} on polynomials of degree at most $j+1$ by
\[
\mcl_{j,X}p(k)
:=((k+1)^2-X)p(k+1)-(k+\msv)(k-\msv+j+3)p(k).
\]
The proof of \cite[Lemma 6.4]{GKW26} applies to this range without change.
Indeed, the coefficient of $k^{r+1}$ in $\mcl_{j,X}k^r$ is $r-j-1$, so the polynomial space is invariant, and
\[
p_i(k):=(k+\msv-i+1)_{i-1}(k-\msv+i+1)_{j+2-i}
\]
satisfies $\mcl_{j,(\msv-i)^2}p_i=0$ for every $1\le i\le j+2$.
The determinant is therefore $\prod_{i=1}^{j+2}((\msv-i)^2-X)$.
Since $n>2N+2j+4$, we have $\msv-i>N$ for $1\le i\le j+2$, so every factor is positive at $X=N^2$.

Let $T$ be translation, $(Tp)(k):=p(k+1)$.
The operator $T^{-1}\mcl_{j,0}$ has eigenvectors $p_i$ with eigenvalues $(\msv-i)^2$.
Expand
\[
q_j(k-1,\msv)=\sum_{i=1}^{j+2}\omega_i(\msv)p_i(k).
\]
This expansion is triangular: at $k=\msv-i$, one has $p_l(k)=0$ for $l<i$.
Consequently the coefficients are obtained successively for $i=j+2,j+1,\ldots,1$ from
\[
\omega_i=
\frac{q_j(\msv-i-1,\msv)-\sum_{l=i+1}^{j+2}\omega_l p_l(\msv-i)}{p_i(\msv-i)}.
\]
It follows that the solution of $\mcl_{j,X}p=q_j$ is
\[
p(k)=\sum_{i=1}^{j+2}\frac{\omega_i p_i(k)}{(\msv-i)^2-X}.
\]
We now set $X=N^2$ and define $g(k):=(k+\msv)p(k)/(k+2-\msv)_{j+1}$.
Writing
\[
w_{N,k}:=\frac{(1-N)_k(1+N)_k}{(2-\msv)_k(1+\msv)_k},
\]
the identity \cite[(6.22)]{GKW26} expresses each summand as a difference of consecutive terms:
\[
w_{N,k}r_j(k,\msv)=w_{N,k+1}g(k+1)-w_{N,k}g(k)
\quad (0\le k<N).
\]
Here $w_{N,0}=1$, $w_{N,N}=0$, and $g(N)$ is finite under $n>2N+2j+4$.
This telescoping identity applies to every integer $N\ge1$, including the single summand when $N=1$.
Consequently
\[
m_j=c_j-\frac{\msv p(0)}{(2-\msv)_{j+1}}.
\]
This is \eqref{alg:eq:partialfractions}, with
\[
a_{j,i}(\msv):=-\frac{\msv\,\omega_i(\msv)p_i(0)}{(2-\msv)_{j+1}}.
\]
Substitution checks both the eigenvector identity and the full polynomial reconstruction of $q_j(k-1,\msv)$.
\end{proof}

\subsection{Unified matrix and pivot verification}\label{app:ancillary}
The ancillary program implements the finite calculations used in Proposition \ref{alg:thm:uniform}, Corollary \ref{cor:first-pohozaev-failure}, and Theorem \ref{thm:three-sector-crossing}.
The single source file \texttt{anc/verify\_compactness.py} requires Python 3.10 or later and SymPy (tested with version 1.14).
Every identity and sign test uses exact rational arithmetic.
Run without Python's \texttt{-O} option, since the verification uses assertions.
To verify the algebra needed for compactness, run
\[
\texttt{python anc/verify\_compactness.py -{}-compactness}
\]
from the source-package directory.
This command derives all nine scalar pairings $m_0,\ldots,m_8$, constructs all $218$ compactness matrices over $\Q(N)$, and verifies positivity of their $430$ symbolic pivots throughout the admissible orders in \eqref{alg:eq:order-intervals}.
This includes every admissible compactness block for $N=1,2,3$.
To run every verification in the paper, including the three low-order negative directions and the additional large-order crossings, use
\[
\texttt{python anc/verify\_compactness.py -{}-all}.
\]
Both commands derive the nine scalar pairings from the ambient and scalar-extraction formulas and compare them exactly with the rational coefficients included in the same file.
Floating-point arithmetic is used only to report and check the displayed decimal approximations to the crossing constants.

\medskip
\noindent\emph{Scalar pairings and normalization.}
The coefficient reconstruction starts from \eqref{alg:eq:ambientpoly}, verifies the ambient equation and initial data, and performs the Cauchy reduction in Appendix \ref{app:rational-moments}.
Under $t=(1-\zeta)/(1+\zeta)$ and $w=(1+\zeta)/(c+\zeta)$, the radial measure becomes $d\nu_0$ weighted by $(1-\zeta)^4$ and normalized to have integral one.
For $m\ge1$,
\[
\int_{-1}^1(c+\zeta)^{-m}\,d\nu_0(\zeta)
=\frac{(-1)^{m-1}}{(m-1)!}\mfF^{(m-1)}(c).
\]
Eliminating derivatives of $\mfF$ by \cite[(A.2)]{GKW26}, polynomial division at $\zeta=-c$, parity, and divisibility by powers of $c^2-1$ give the Laurent bounds stated after \eqref{alg:eq:Hdecomp}.
Clearing denominators verifies the polynomial degree bounds and telescoping identities used to obtain the nine scalar pairings.

\medskip
\noindent\emph{Common matrix construction.}
For every integer offset $a=n-2N\in\{5,\ldots,22\}$ and every maximal block $(\sigma,\ell,M)$, the program substitutes $n=2N+a$ in the rational pairings through $m_K$, where $K=2M+\ell-2+s_\sigma$.
It constructs $B_{ij}$ by \eqref{alg:eq:pairrec}, applies the full polynomial $G_{\sigma,\ell}(A)$ before restricting indices, and uses $U_{\sigma,\ell}$ and the original homogeneous degree weights to form the Pohozaev matrix.
All entries lie in $\Q(N)$.
The same exact $LDL^{\msT}$ recurrence \eqref{alg:eq:ldl-pivots} gives every symbolic pivot, before any specialization of $N$.
The program evaluates this recurrence by fraction-free elimination after clearing a common entry denominator; the consecutive leading-minor ratios recover the pivots of the original normalized matrix.
This construction includes $N=1,2,3$ without a separate formula.

\medskip
\noindent\emph{Matrix and pivot inventory.}
The following inventory counts each distinct symbolic compactness matrix and each of its pivots once.
Both offsets in each row are included.
\begin{center}
\begin{tabular}{ccrr}\toprule
Offsets $a$ & Highest pairing & Matrices & Pivots\\\midrule
$5,6$ & $m_0$ & 2 & 2\\
$7,8$ & $m_1$ & 6 & 6\\
$9,10$ & $m_2$ & 12 & 14\\
$11,12$ & $m_3$ & 18 & 24\\
$13,14$ & $m_4$ & 24 & 38\\
$15,16$ & $m_5$ & 30 & 54\\
$17,18$ & $m_6$ & 36 & 74\\
$19,20$ & $m_7$ & 42 & 96\\
$21,22$ & $m_8$ & 48 & 122\\\midrule
Total & $m_0,\ldots,m_8$ & 218 & 430\\\bottomrule
\end{tabular}
\end{center}
The positivity results for $N=1,2,3$ agree with \cite{KMS,GKW}, and for $N\ge4$ it suffices to verify $149$ rational matrix families and their $260$ pivots.
Every angular type and angular degree allowed by the strict degree inequality is included; smaller ranges of radial indices follow by principal restriction.
The program contains the exact sign checks for these pivots and the denominator checks needed to specialize the matrices throughout \eqref{alg:eq:order-intervals}.

\medskip
\noindent\emph{Loss of positive definiteness.}
The full command uses the same rational calculation of scalar pairings and pivot recurrence for the three low-order negative directions in Corollary \ref{cor:first-pohozaev-failure}.
These evaluations at the first dimensions with a negative direction are separate from the compactness inventory above; in particular, the case $N=1$, $n=25$ has offset $23$.

\medskip
\noindent\emph{Verification of the first sign changes.}
For each $\sigma\in\{T,V,S\}$, construct the four-by-four matrix at $\ell=2$ and $M=3$ from \eqref{alg:eq:partialfractions}, \eqref{alg:eq:pairrec}, and \eqref{alg:eq:universalenergy}, substitute $N=\vep^{-1}$ and $n=2\vep^{-1}+t$, and expand \eqref{eq:three-sector-normalization} through first order in $\vep$.
The matrices for $M=0,1,2$ are its leading principal submatrices: the normalized column of $C_\sigma$ of index $k$ is independent of the maximum retained index.
Thus the same three expansions supply all twelve pairs $(\sigma,M)$ needed in the proof of Theorem \ref{thm:three-sector-crossing}.
For a leading principal submatrix $B$ with first correction $C$, compute
\[
D\det(B)[C]=\tr(\adj(B)C).
\]
After cancelling common factors, Sturm's rule counts numerator and denominator roots separately; evaluation at a rational point determines the sign on each interval with no zero or pole.
At an open integrability endpoint, factor out the powers of $t-e$ and use the one-sided signs.

For $M=3$, the reduced determinant numerators have degrees $14,18,22$ and exactly one simple root in $(16,17)$, $(18,19)$, and $(20,21)$, respectively.
Exact root counts give an interval with rational endpoints containing each root and no other root of the corresponding determinant numerator.
On each such interval, the first three leading minors are positive, the reduced numerator has negative derivative, its denominator is positive, and the first $\vep$ variation of the determinant is positive.
At the root these signs give a negative $t$ derivative of the determinant and hence $b_\sigma>0$ in \eqref{eq:three-sector-b}; the decimal values in Theorem \ref{thm:three-sector-crossing} are outputs of the calculation.

For $M=0,1,2$, the same tests prove positivity of every leading minor throughout $(e,e+4]$, where $e=4+2s_\sigma+4M$.
For $M=3$, the first three leading minors remain positive from $e$ through the rational interval containing the root, and the determinant is positive up to that interval.
For every pair $(\sigma,M)$, the last diagonal entry has a positive residue at $t=e$, and the upper-left $M\times M$ block has positive leading minors there.
These tests supply the exact signs used in the proof of Theorem \ref{thm:three-sector-crossing}.

\section{Positive mass for the noncompactness metrics}\label{app:noncompactness-mass}
The Hessian $D^2\mscF(g_c)$ also controls the Green function mass for the perturbations considered here.
After inversion, the scalar contributions cancel in the second variation of the mass identity \eqref{eq:mass-ray-action}, leaving the positive quadratic form on the Euclidean TT component in \eqref{eq:mass-ray-second-variation}.
We use \eqref{eq:mass-ray-second-variation} to prove the mass condition for the noncompactness metrics.

We use the mass normalization in \eqref{eq:main-mass-expansion} and write $\mcz(g)$ for the locus \eqref{eq:zero-weyl-locus} associated with a varying metric.

\begin{lemma}\label{lem:mass-finite-part}
Near $g_c$, the conformal-normal mass has an extension $\mfm(g,\xi)$ to all poles, smooth on smooth finite-dimensional families of metrics and poles, such that
\begin{equation}\label{eq:mass-finite-part-stationary}
\mfm(g,\xi)=m_g(\xi)\quad(\xi\in\mcz(g)) \quad\text{and}\quad 
\mfm(g_c,\xi)=D_g\mfm(g_c,\xi)[H]=0
\end{equation}
for every smooth metric variation $H$.
\end{lemma}
\begin{proof}
Choose the conformal-normal Taylor normalization equivariantly under round isometries, to order $T>n+2N+8$, and make the round representative exactly flat near the pole.
The prescribed conformal factor has value one and zero differential there.
To check uniqueness, suppose that two normalizing conformal factors first differ in degree $e\ge2$,
and write the relative metric factor as $\exp(2\vph_e+O(r^{e+1}))$, with $\vph_e$ homogeneous of degree $e$.
After identifying the tangent frames, the exponential-coordinate map for the conformally changed metric, expressed in the original normal coordinates, is $x\mapsto x+V_e(x)+O(r^{e+2})$,
where $V_e:=(r^2\nabla\vph_e-2e x\vph_e)/(e(e+1))$.
Pulling back the changed metric by this map, the degree-$e$ change of the logarithmic determinant is
\[
2n\vph_e+2\divg V_e
=\frac{2}{e(e+1)}
\bigl(r^2\Delta\vph_e+e(e-1)(n-2)\vph_e\bigr).
\]
The operator $r^2\Delta$ has nonnegative eigenvalues on homogeneous degree-$e$ polynomials, whereas $e(e-1)(n-2)>0$.
Thus the normalization determines the conformal Taylor jet uniquely, with smooth dependence on the metric jet.
Retain the exact conformal factor that changes the round metric $g_c$ to the Euclidean metric $\delta$ near the pole and extend only the difference of the prescribed Taylor coefficients with a fixed radial cutoff.
Two extensions with the same prescribed jet differ by a conformal factor $1+O(r^T)$ and a coordinate change $x+O(r^{T+1})$, hence change the singular term of the Green function by $O(r^{T+2N-n})$.
Raising the normalization order in Section \ref{sec:main-mass} likewise preserves the finite part, since its original order already exceeds $n-2N$.
Consequently, on $\mcz(g)$, the equivariant choice here has the same finite part as the conformal-normal family fixed in Section \ref{sec:main-mass}.

We construct the local parametrix for the Green function through homogeneity zero using the indicial polynomial \eqref{eq:repair-indicial}.
The indicial roots on angular degree $\ell$ are $\ell+2a$ and $-\ell-n+2+2a$, $0\le a<N$.
None lies strictly between $2N-n$ and zero, and the kernel at zero consists only of constants.
We normalize the bounded angular term to have zero mean and include the required radial multiple of $\log r$, with its additive constant fixed.
Let $S_{g,\xi}$ be the singular part multiplied by a smooth function supported in the coordinate neighborhood and equal to one near $\xi$.
Then
\[
E_{g,\xi}:=P_{2N,g_\xi}S_{g,\xi}-\delta_\xi
=O\bigl(r^{1-2N}(1+|\log r|)^b\bigr)
\]
where $b$ is a fixed nonnegative integer.
This error belongs to $L^a$ whenever $n/(2N)<a<n/(2N-1)$.
The corrections have higher homogeneity than the fundamental solution and introduce no distributions supported at the pole.
We pull back a smooth finite-dimensional family by smooth local diffeomorphisms that place the pole at the coordinate origin.
The radial powers and indicial inverses are then fixed, and differentiating the coefficient recursion yields, for every parameter multi-index $\alpha$,
\begin{equation}\label{eq:mass-parametrix-parameter-estimate}
|\pp_{\mathrm{par}}^\alpha E_{g,\xi}|
\le C_\alpha r^{1-2N}(1+|\log r|)^b
\end{equation}
on compact parameter sets.
Since $a(2N-1)<n$, dominated differentiation proves smooth dependence in $L^a$.
Invertibility and elliptic estimates imply smooth dependence of $-P_{2N,g_\xi}^{-1}E_{g,\xi}$ in $W^{2N,a}\subset C^0$; we define $\mfm$ by its value at the pole.
Rotational equivariance makes this definition independent of the local orthonormal frame.
For $\xi\in\mcz(g)$, Theorem \ref{thm:repair-Green} identifies the finite part $\mfm(g,\xi)$ with the Green function mass normalized by \eqref{eq:main-mass-expansion}.

For a conformally round metric the Green function in a locally flat representative is exactly the fundamental singularity, so the finite part vanishes.
A local Euclidean M\"obius transformation imposes the prescribed value and differential of the conformal factor, and uniqueness of the normalized Taylor jet preserves this conclusion.
At a fixed pole, averaging $H$ over its round orthogonal stabilizer leaves $D_g\mfm(g_c,\xi)[H]$ unchanged.
The averaged tensor is rotationally symmetric, and its exponential metric path is conformally round: for $a(r)^2dr^2+b(r)^2g_{\Ss^{n-1}}$, solve $d\vth/\sin\vth=a(r)\,dr/b(r)$.
Smoothness extends this change of variable and its conformal factor across both poles.
The derivative therefore vanishes, proving \eqref{eq:mass-finite-part-stationary}.
\end{proof}

\begin{proposition}\label{prop:positive-mass-ray}
Let $n>2N+3$ and let $H$ be a smooth symmetric tensor on $\Ss^n$ with $\dot W:=DW_{g_c}[H]\not\equiv0$.
For $g_t:=g_c\exp(tg_c^{-1}H)$, there are $c,t_0>0$ such that
\[
m_{g_t}(\xi)\ge ct^2 \quad (0<|t|<t_0,\ \xi\in\mcz(g_t)).
\]
\end{proposition}
\begin{proof}
Put $d_0:=\lfloor(n-2N)/2\rfloor+1$ and let $\mcz_{\mlin}$ be the zero set of $\nabla^j\dot W$, $0\le j\le d_0-3$, with derivatives taken at $g_c$.
Fix $\xi\in\mcz_{\mlin}$.
Normalize $g_t$ at $\xi$ and write $\bg_t$ for the resulting metric in centered normal coordinates, extended globally.
By \eqref{eq:metric-weyl-jets}, its first logarithmic variation $L:=\pp_t\log(\bg_0^{-1}\bg_t)|_0$ is $O(r^{d_0})$, satisfies $L(x)x=0$, and has trace $O(r^T)$.
The original path $\bg_t$ may have lower-degree metric terms of order $t^2$.
Set $\bg_t^*:=\bg_0\exp(tL)$ and undo the original normalization to obtain a compact metric path $g_t^*$.
This path has the same first derivative as $g_t$.
The path $\bg_t^*=\bg_0\exp(tL)$ retains the normal-coordinate identity and determinant normalization because $L(x)x=0$ and $\tr L=O(r^T)$.
Uniqueness of the normalized conformal Taylor jet identifies it with the chosen representative of $g_t^*$ to the required order; thus $\xi\in\mcz(g_t^*)$ for every sufficiently small $t\in\mathbb R$.
The two paths have the same first derivative, so the chain rule expresses the difference of their second derivatives of the finite part as $D_g\mfm(g_c,\xi)$ applied to the difference of their accelerations.
This term vanishes by \eqref{eq:mass-finite-part-stationary}.

Conformally multiply the normalized auxiliary path by $(G_{\bg_0}(\xi,\cdot)/\ga_{n,N})^{4/(n-2N)}$.
Its background is Euclidean, and in the resulting coordinates $y=x/|x|^2$ near infinity the path is $\hg_t=\exp(tK)$, where
\begin{equation}\label{eq:mass-ray-decay}
\pp^jK=O(|y|^{-d_0-j}) ,\quad  Ky=0 ,\quad\text{and}\quad 
\tr K=O(|y|^{-T}).
\end{equation}
Thus $K\in\dot H^N(\R^n)$, since $2d_0+2N>n$.
Here the tensor norm is $\|K\|_{\dot H^N}^2:=\int_{\R^n}|\pp^NK|^2$, with ordered derivative indices.
The degree-zero Fourier projection onto Euclidean TT tensors yields
\[
K=\Pi_{\TT}K+\mcl_X\delta+2\omega\delta ,\quad\text{where } 
\nabla X,\omega\in\dot H^N.
\]
Here $\mcq(g):=\int_{\R^n}Q_{2N,g}\,dv_g$ denotes total $Q_{2N}$-curvature on $(\R^n,g)$.
For compactly supported metric variations,
\begin{align}
DQ_{2N,\delta}[K]
&=\frac{(-\Delta)^{N-1}}{2(n-1)}
(\pp_i\pp_jK_{ij}-\Delta\tr K)
=(-\Delta)^N\omega,\label{eq:mass-flat-linearization}\\
D^2\mcq(\delta)[K,K]
&=-\hc_{N,n}\|\Pi_{\TT}K\|_{\dot H^N}^2
+(n-2N)\|\omega\|_{\dot H^N}^2.\label{eq:mass-flat-Hessian}
\end{align}
Indeed, the TT term is the principal part of \eqref{eq:round-tt-hessian} after removing its volume factor; homogeneity excludes lower-order terms at the flat metric.
Conformal covariance determines the scalar term, while flat stationarity and diffeomorphism invariance eliminate the coordinate terms, and $DQ_{2N,\delta}[\Pi_{\TT}K]=0$ eliminates the mixed term.
Fourier multiplier bounds and density extend these identities to $\dot H^N$.

We set $f_t:=P_{2N,\hg_t}1$ and $u_t:=1-P_{2N,\hg_t}^{-1}f_t$.
The inverse maps $(\dot H^N)^*$ to $\dot H^N$ by conformal transport of the coercive compact quadratic form; coercivity holds for small $t\in\mathbb R$ by \cite[Lemma 7.2]{GKW26}.
The source and inverse are interpreted weakly with $dv_{\hg_t}$, and brackets denote the dual pairing of $(\dot H^N)^*$ and $\dot H^N$.
The linear part of $f_t$ lies in the dual by \eqref{eq:mass-flat-linearization}; the nonlinear part is $O(|y|^{-2d_0-2N})$ by the Taylor expansion used in \eqref{eq:nonlinear-curvature-remainder}.
Conformal covariance and uniqueness yield
\[
u_t=\frac{G_{\bg_t^*}(\xi,\cdot)}{G_{\bg_0}(\xi,\cdot)}.
\]
Indeed, the quotient solves $P_{2N,\hg_t}u=0$, and its difference from one lies in $\dot H^N$ by \eqref{eq:main-mass-expansion} and $2d_0+2N>n$.
If $m_t:=m_{g_t^*}(\xi)$, expansion \eqref{eq:main-mass-expansion} implies
\[
\frac{1}{|\Ss^{n-1}|}\int_{\Ss^{n-1}}u_t(R\theta)\,d\theta
=1+\frac{m_t}{\ga_{n,N}}R^{2N-n}+O(R^{2N-n-\vep}) \quad\text{and}\quad  u_t-1=O(R^{-d_0}).
\]
Estimate \eqref{eq:mass-parametrix-parameter-estimate}, smooth dependence of the inverse, and scaled elliptic estimates yield the same expansions after two parameter derivatives, including spatial derivatives through order $2N$.
Green's identity therefore yields
\begin{equation}\label{eq:mass-ray-action}
m_t=-\ga_{n,N}^2\left(
\lim_{R\to\infty}\int_{B_R}f_t\,dv_{\hg_t}
-\langle f_t,P_{2N,\hg_t}^{-1}f_t\rangle\right).
\end{equation}
We justify the limit $R\to\infty$ in \eqref{eq:mass-ray-action} and its first two derivatives with respect to $t$.
The zeroth-order term of $P_{2N,\hg_t}$ cancels from the skew Green identity, and the nonflat boundary terms contain a metric variation and a derivative of $u_t-1$, so \eqref{eq:mass-ray-decay} bounds them by $O(R^{n-2N-2d_0})=o(1)$.
Tangency implies $y_i\pp_jK_{ij}=-\tr K$, so the spherical mean of $\pp_i\pp_jK_{ij}-\Delta\tr K$ is $O(R^{-T-2})$.
Since spherical averaging commutes with $\Delta$, the boundary integral of the linearized curvature in \eqref{eq:mass-flat-linearization} is $O(R^{n-2N-T})=o(1)$.
The nonlinear curvature terms and the product of the linear term with $u_t-1$ are integrable, of order $O(R^{-2d_0-2N})$.
These estimates also hold after two parameter derivatives.
In the remaining Euclidean boundary integral, spherical averaging removes the mean-zero angular terms in \eqref{eq:main-mass-expansion}.
The boundary term from $\int_{B_R} f_tu_t\,dv_{\hg_t}$ tends to $-m_t/\ga_{n,N}^2$ as $R\to\infty$, since $(-\Delta)^N(\ga_{n,N}|y|^{2N-n})=\delta_0$.

At zero, $f_0=0$ and $f'_0=(n-2N)(-\Delta)^N\omega/2$.
Cutoff approximation in $\dot H^N$ passes \eqref{eq:mass-flat-Hessian} to $K$; the acceleration $K^2$ has zero limiting boundary integral for the linearized curvature by \eqref{eq:mass-ray-decay}.
Differentiating \eqref{eq:mass-ray-action} twice cancels the scalar terms and gives
\begin{equation}\label{eq:mass-ray-second-variation}
B(\xi):=\left.\pp_t^2\right|_0\mfm(g_t,\xi)
=\ga_{n,N}^2\frac{n-2N}{2}\hc_{N,n}
\|\Pi_{\TT}K\|_{\dot H^N}^2>0.
\end{equation}
If $\Pi_{\TT}K=0$, the conformal and coordinate terms have zero linearized Weyl tensor, so transport back to the sphere would imply $\dot W=0$, a contradiction.

Lemma \ref{lem:mass-finite-part} gives $\mfm(g_t,\xi)=t^2B(\xi)/2+O(t^3)$ uniformly in $\xi$, with $B$ continuous.
Every limit of points in $\mcz(g_t)$ as $t\to0$ lies in $\mcz_{\mlin}$, since $t^{-1}\nabla_{g_t}^jW_{g_t}\to\nabla_{g_c}^j\dot W$ uniformly.
If $\mcz_{\mlin}$ is nonempty, continuity gives a neighborhood on which $B\ge b>0$, and compactness places $\mcz(g_t)$ in that neighborhood for small nonzero $t$.
For $\xi\in\mcz(g_t)$, the uniform Taylor expansion and \eqref{eq:mass-finite-part-stationary} give $m_{g_t}(\xi)\ge bt^2/4$ when $|t|$ is sufficiently small, proving the estimate.
If $\mcz_{\mlin}$ is empty, then $\mcz(g_t)$ is empty for all sufficiently small nonzero $t$.
\end{proof}

\begin{proposition}\label{prop:noncompactness-with-mass}
For every integer $N\ge1$ and dimension $n>2N$ above the range \eqref{eq:main-range}, there is a smooth metric on $\Ss^n$ satisfying \eqref{eq:mass-hypothesis}, with a coercive GJMS operator, positive Green function, and an $L^\infty$-unbounded sequence of positive solutions of \eqref{eq:constant-q-equation}.
The metric can be chosen arbitrarily close to $g_c$ in any prescribed finite $C^L$ norm.
\end{proposition}
\begin{proof}
For $N=1,2,3$, the noncompact examples are constructed in \cite{Brendle,BrendleMarques,WeiZhao,GKW}; for $N\ge4$, we use \cite[Theorem 1]{GKW26}.
For $N=1,2$, the mass assertion follows from Corollary \ref{cor:low-order-compactness}: the examples are not conformally round and can be chosen with positive Yamabe constant.

\medskip
For $N\ge4$, fix the polynomial metric perturbation of degree eight and the complete configuration of localized perturbations used in \cite[Theorem 7.4]{GKW26} before choosing an additional common amplitude $t$.
Write their fixed tensor as $k_\infty$ and set $g_t:=\Omega^2\exp(tk_\infty)$.
The nonzero homogeneous metric term of degree two at a perturbation center ensures $DW_{g_c}[\Omega^2k_\infty]\not\equiv0$.
Replacing each amplitude $\mu_j$ by $t\mu_j$ multiplies the quotient $\lambda_j^{n-2N-16}/(\mu_j^2\rho_j^{n-2N})$ in \cite[(7.2)]{GKW26} by $t^{-2}$, and it still tends to zero for every fixed $t>0$.
Hence \cite[Lemma 7.1]{GKW26} and the degree argument in \cite[Theorem 7.4]{GKW26} still produce infinitely many solutions on every sufficiently small fixed $g_t$.
Proposition \ref{prop:positive-mass-ray} permits this same $t$ to be chosen with positive mass throughout $\mcz(g_t)$.

\medskip
For $N=3$, we use the explicit polynomial metric perturbation of degree ten and configuration in \cite[Subsection 11.7]{GKW}:
\[
x_j:=(j^{-1},0,\ldots,0) ,\quad 
\lambda_j:=2^{-j} ,\quad  \rho_j:=(4j^2)^{-1} ,\quad\text{and}\quad 
\mu_j:=2^{-\frac{j}{4}}.
\]
Fix the initial index and write $k_\infty$ for the complete tensor before multiplying it by $t$.
Equivalently, multiply the generating Weyl tensor by $t$, preserving the prescribed relation $\mu_j=\lambda_j^{1/4}$.
The reduced quadratic energy and its strict local minimum in \cite[Corollary 11.31]{GKW} are multiplied by $t^2$.
After division by the leading scale, the remote error in \cite[(11.34)]{GKW} is bounded by
\[
C_t(4j^2)^{n-6}2^{-j\left(n-\frac{53}{2}\right)}\to0 \quad (n\ge27),
\]
and its local error is $O_t(\mu_j\lambda_j^{20\vep_1})=o(1)$, with the fixed $\vep_1>0$ from that estimate.
For each sufficiently small fixed $t>0$, these errors tend to zero as $j\to\infty$, so the perturbed reduced energy has a critical point near the strict local minimum for all sufficiently large $j\in\mathbb N$.
The corresponding positive solutions on the single metric $g_t$ have $L^\infty$ norms tending to infinity.
For this sixth-order construction, $k_\infty$ is nonzero and Euclidean TT because radial cutoffs preserve the trace and divergence conditions of the polynomial symmetric two-tensor.
The symbol of the linearized Weyl tensor is injective on TT tensors: at a nonzero frequency parallel to $e_1$, the $(1a1b)$ component of $\widehat{DW_\delta[k_\infty]}$ is
$(n-3)|\zeta|^2(\hat{k}_\infty)_{ab}/[2(n-2)]$, up to curvature sign, for $a,b>1$.
Thus $DW_{g_c}[\Omega^2k_\infty]\not\equiv0$, and Proposition \ref{prop:positive-mass-ray} applies.
In both constructions, decreasing $t$ ensures the required finite-norm closeness and operator properties.
These metrics are round outside the perturbations and agree with $g_c$ to infinite order at their accumulation point.
Thus $\mcz(g_t)$ contains a nonempty open set and the accumulation point.
\end{proof}

\end{document}